\documentclass[11pt]{article}
\usepackage{amsmath}
\usepackage{amssymb}
\usepackage{amsthm}
\usepackage{enumerate}
\usepackage{amsfonts}
\theoremstyle{plain}
\newtheorem{theorem}{Theorem}[section]
\newtheorem{lemma}[theorem]{Lemma}
\newtheorem{proposition}[theorem]{Proposition}
\newtheorem{corollary}[theorem]{Corollary}
\newtheorem*{theoremA}{Theorem A}
\newtheorem*{theoremB}{Theorem B}
\theoremstyle{remark}
\newtheorem{remark}[theorem]{Remark}
\theoremstyle{definition}
\newtheorem{definition}[theorem]{Definition}
\newtheorem{example}[theorem]{Example}
\newtheorem{convention}[theorem]{Convention}
\usepackage{xcolor}
\definecolor{linkcol}{rgb}{0.10,0.20,0.55}
\definecolor{citecol}{rgb}{0.05,0.40,0.15}
\usepackage[colorlinks,linktoc=all,linkcolor=linkcol,citecolor=citecol,urlcolor=linkcol,bookmarksnumbered]{hyperref}

\newcommand{\Mod}{\mathrm{Mod}}
\newcommand{\T}{\mathcal{T}}
\newcommand{\N}{\mathcal{N}}
\newcommand{\Pp}{\mathcal{P}}
\newcommand{\E}{\mathcal{E}}
\newcommand{\Lst}{L_*}
\newcommand{\red}{\mathrm{red}}
\newcommand{\glue}{\mathrm{gl}}
\newcommand{\cut}{\mathrm{cut}}
\newcommand{\rot}{\mathrm{rot}}
\newcommand{\Stab}{\mathrm{Stab}}
\newcommand{\Nt}{\widetilde{\mathcal{N}}}
\newcommand{\Cos}{\mathrm{Cos}}
\newcommand{\Dt}{\Delta_{\mathrm{top}}}

\begin{document}
\title{Hyperbolic string vertices as geometric chains}
\author{Hao Yu}
\date{}
\maketitle

\begin{abstract}
For $0<L\le 2\,\mathrm{arcsinh}\,1$, the hyperbolic string vertices of Costello and Zwiebach consist of the bordered hyperbolic surfaces with geodesic boundary of length $L$ whose interior simple closed geodesics have length at least $L$. We show that the corresponding regions in moduli space are compact oriented orbifolds with corners whose boundaries are, including orientations and automorphism groups, the spaces of twisted sewings. Hence their generating chain satisfies $(\partial+\Delta)\mathbf V=0$ exactly in a complex of geometric chains. The connected form $\partial V+\Delta V+\frac12\{V,V\}=0$ cannot hold exactly for geometric chains; we determine its exact form, show that it holds modulo gluing and covering relations in a complex computing the homology of the moduli spaces, and that the vertex currents are the unique top-dimensional solution. Finally we compare with Costello's homotopical framework. In a strict model given by hyperbolic surfaces with marked boundary points, every solution of his master equation with his initial condition integrates exactly to the vertex currents, up to the sign $(-1)^{3g-3+n}$; such solutions exist, are unique up to homotopy, and correspond to the hyperbolic vertices. Under a continuity assumption on the comparison map and granting the inputs of Costello's homology computation, the hyperbolic vertices represent his homotopy class of solutions, for a completed notion of homotopy.
\end{abstract}
\tableofcontents

\section{Introduction}\label{sec:intro}
Fix $L>0$. For a compact oriented surface $X$ with a hyperbolic metric for which every boundary component is a geodesic of length $L$, Costello and Zwiebach \cite{CZ} consider the condition
\begin{equation}\label{eq:systole}
\ell_\gamma(X)\ge L\quad\text{for every simple closed geodesic $\gamma$ in the interior of $X$.}
\end{equation}
For $L\le\Lst:=2\,\mathrm{arcsinh}\,1$ they show that the sets $\mathcal V_{g,n}(L)$ of such surfaces satisfy the geometric Batalin--Vilkovisky master equation of closed string field theory \cite{Zw}: the boundary of $\mathcal V_{g,n}(L)$ is the union of the surfaces obtained by sewing two boundary components of a vertex of lower complexity, or one boundary component of each of two vertices, with an arbitrary twist. The same statement for $L$ below $\mathrm{arcsinh}\,1$, with a different proof, is due to Andersen, Borot and Orantin \cite[Theorem on the quantum master equation in the section on string field theory]{ABO}, who attribute the idea that systolic sets satisfy the master equation to Costello. Both proofs compare sets of surfaces: they show that each boundary point arises by sewing and that sewing lands in the boundary. Orientations, multiplicities and the factors coming from automorphisms are not discussed, and the statement at the level of chains is asserted but not proved \cite{ABO}.

This paper gives the chain level statement in the geometric chains of \cite{Yu}. These chains are generated by maps from compact oriented orbifolds with corners, so a compact region of a moduli space with a description of its boundary as an oriented orbifold defines a chain whose boundary can be computed exactly, without choosing a fundamental chain. The work consists in establishing the orbifold structure and in matching the boundary with the sewing spaces, including orientations and automorphism groups.

\textit{The model.} We work in the finite dimensional moduli spaces of hyperbolic surfaces with geodesic boundary of length $L$, as in Mirzakhani \cite{Mir07}, not in the infinite dimensional spaces of surfaces with local coordinates of \cite{CZ}. Sewing two boundary components requires the choice of an orientation reversing isometry between them; these choices form a circle, which plays the role of the twist. The boundary components are unlabelled and the surfaces may be disconnected, with unordered components; the reason is explained below.

\textit{Results.} Let $\N_{\chi,n}$ be the moduli space of possibly disconnected hyperbolic surfaces $X$ with $n$ unlabelled boundary components of length $L$, all components of negative Euler characteristic and $\chi(X)=\chi$, regarded as an orbifold in the wide sense (the generic surface may have non-trivial automorphisms). Let $V_{\chi,n}\subseteq\N_{\chi,n}$ be the region defined by \eqref{eq:systole} for all components.

\begin{theoremA}[Theorem \ref{thm:A}]
For $0<L\le\Lst$, $V_{\chi,n}$ is a compact suborbifold with corners of $\N_{\chi,n}$ of full dimension $-3\chi-n$. A point has depth $k$ if exactly $k$ interior simple closed geodesics have length $L$; these are pairwise disjoint, and near the point $V_{\chi,n}$ is given in Fenchel--Nielsen coordinates by $\ell_{\gamma_1},\dots,\ell_{\gamma_k}\ge L$.
\end{theoremA}

For punctured surfaces and small thresholds this is \cite[Proposition 1.1]{JW}, where it is called well known; we include the proof for bordered surfaces and the explicit range $L\le\Lst$, including the endpoint.

\begin{theoremB}[Theorem \ref{thm:B}]
For $0<L\le\Lst$, there is an orientation reversing isomorphism $\Theta$ of compact orbifolds with corners in the wide sense from the space of sewings of $V_{\chi,n+2}$ along an unordered pair of boundary components onto the boundary $\partial V_{\chi,n}$ in the sense of Joyce \cite{jo}, through which the gluing morphism factors. Here $\N_{\chi,n}$ is oriented by $\prod_i d\ell_i\wedge d\tau_i$ in Fenchel--Nielsen coordinates with left twists, and the space of sewings is oriented by the base followed by the circle of gluing isometries, oriented by rotation.
\end{theoremB}

In particular $\partial\mathbf V_{\chi,n}=-\mathbf S_{\chi,n}$ for the weighted fundamental chain $\mathbf V_{\chi,n}$ of $V_{\chi,n}$ and the weighted chain $\mathbf S_{\chi,n}$ of the sewing space (Corollary \ref{cor:sewing}). In Section \ref{sec:master} we define the sewing operator $\Delta$ on chains over orbifolds in the wide sense, show $\Delta^2=0$ and $\partial\Delta+\Delta\partial=0$, and deduce:

\textit{Theorem C} (Theorem \ref{thm:C}). The sum $\mathbb V$ of the fundamental chains of all $V_{\chi,n}$ satisfies $(\partial+\Delta)\mathbb V=0$ exactly, before passing to the quotient by degenerate chains, in the complex $\mathcal C$ of chains over the orbifolds $\N_{\chi,n}$ in the wide sense. This auxiliary complex does not compute homology; in the complex of \cite{Yu} the content of Theorem C is Corollary \ref{cor:sewing}.

This is the exponential form $(d+\hbar\Delta)e^{S/\hbar}=0$ of the master equation, the power of $\hbar$ of a surface being $-(\chi+n)/2=\sum_i(g_i-1)$. Its connected part (Corollary \ref{cor:connected}) contains, whenever a separating curve cuts a surface into two pieces of the same type $t$, the term $\Delta_\times[\mathrm{Sq}(V_t)]$, the sewing of an unordered pair of vertices of type $t$. Integration turns this term into $\frac12\{T_t,T_t\}$, and the integration currents satisfy the usual form of the master equation (Corollary \ref{cor:currents}).

\textit{Rigidity} (Theorem \ref{thm:D}). Among compactly supported top dimensional currents on the moduli spaces of connected surfaces, the integration currents of the hyperbolic vertices are the only solution of the master equation with the given term on the pair of pants. This is an analogue, in a complex without room for homotopies, of Costello's uniqueness up to homotopy \cite{Cos05}.

\textit{Why disconnected surfaces.} If a separating geodesic cuts a vertex into two pieces of the same topological type $t$, the corresponding boundary face is the quotient of the space of sewings of $V_t\times V_t$ by the involution exchanging the two pieces. In the usual connected form $\partial V+\Delta V+\frac12\{V,V\}=0$ this face is represented by $\frac12\{V_t,V_t\}$, whose source is the space of sewings itself, with coefficient $\frac12$. Geometric chains distinguish a chain $[S,f\circ q]$ from $2[S/\mathbb Z_2,f]$ when $q:S\to S/\mathbb{Z}_2$ is a non-trivial double covering, so the connected form cannot hold exactly as soon as $\dim V_t>0$; this happens for instance for $(g,n)=(2,0)$, cut into two one-holed tori, and for $(0,6)$, cut into two four-holed spheres. With unordered components the two sides of the identity contain the same quotient, and the identity is exact (Corollary \ref{cor:connected}); the failure of the connected form with $\frac12\{V_t,V_t\}$ is Proposition \ref{prop:obstruction}. With the gluing relation of \cite{YuF} in place of $D$, a covering of degree $d$ is equivalent to $d$ times its base, and after reduction the connected form holds in the quotient by gluing relations, a complex which still computes the homology of the moduli spaces (Theorem \ref{thm:connected-G}). With covering relations also in the direction of the kernels, the relations can be put on the auxiliary complex itself; then $\Delta$ descends, the quotient computes homology, $\mathbb V\equiv e^V-1$, and the connected vertices form a Maurer--Cartan element of a differential graded Lie algebra computing the homology of the moduli spaces (Theorem \ref{thm:Tw}, Corollary \ref{cor:exp}).

\textit{Comparison with Costello} (Section \ref{sec:Costello}). Costello \cite{Cos05} defines string vertices homotopically, as the solution, unique up to homotopy, of the master equation in an algebra built from homotopy coinvariants of the circle actions on moduli spaces of surfaces with parametrized boundary; on the strict quotients, where our chains live, his sewing operator is not defined on singular chains \cite[footnote 6]{Cos05}. As a model of his framework we use the moduli spaces $\Nt$ of hyperbolic surfaces with a marked point on each labelled boundary component, on which gluing is strictly associative, and we construct a residue map $\Phi$, which takes the coefficient of $t_1^{-1}\cdots t_n^{-1}$ and projects to the moduli orbifolds, from Costello's algebra for $\Nt$ to the complex of Theorem \ref{thm:Tw}; multiplied by $(-1)^{3g-3+n}$ on the moduli space of type $(g,n)$, it intertwines the BV operators and the brackets.

\textit{Theorem E} (Theorem \ref{thm:E}). Every solution of Costello's master equation for $\Nt$ with his initial condition integrates exactly, up to the sign $(-1)^{3g-3+n}$, which depends on the orientation of Costello's circle operator (Remark \ref{rem:signs}(2)), to the currents $T_{g,n}$ of the hyperbolic vertices. Such solutions exist and are unique up to homotopy, and, if the comparison map $h$ from $\Nt$ to Costello's moduli spaces is continuous, they represent Costello's class. Here and below homotopies and Costello's class refer to the completed notion of homotopy of Section \ref{ssec:Cos-framework}, which may be coarser than Costello's own if one does not use his sketched \cite[Lemma 8.1.1]{Cos05}. The integration statement concerns Maurer--Cartan elements in the model $\Nt$; for Costello's solutions on his own moduli spaces there is no direct integration statement, since there is no map from them to $\N_L$ compatible with gluing.

\textit{Theorem F} (Theorem \ref{thm:F}). The residue map induces a bijection of homotopy classes of Maurer--Cartan elements from Costello's algebra for $\Nt$ to the algebra of Theorem \ref{thm:Tw}. If $h$ is continuous, the comparison map to Costello's moduli spaces also induces a bijection, under which Costello's class corresponds to the class of the hyperbolic vertices $V_\hbar$.

The first statement of Theorem E uses only the rigidity theorem. The existence and uniqueness statement and Theorem F(1) use in addition a computation of the equivariant homology of $\Nt$ by a Koszul resolution, the standard extension of the theorem of Jones and Hood--Jones \cite{J87,HJ87} to tori (Proposition \ref{prop:equivariant}), Theorem \ref{thm:Tw}, Harer's bound on the homology of mapping class groups \cite{Ha} (Lemma \ref{lem:bound}), and an order-by-order argument (Lemma \ref{lem:MC}) in place of the lemma of Goldman--Millson type which Costello sketches. The transfer to Costello's own model is conditional: it assumes that $h$ is continuous, Costello not having specified the topology of his moduli spaces, and it uses his computation of the homology of his algebra, which rests on \cite[Proposition 5.1.1]{Cos05}, whose proof is sketched, and on the rational homotopy equivalence of $M_g(n)/(S^1)^n$ with $M_{g,n}$. Neither result needs $\Phi$ to be a quasi-isomorphism: only the associated graded homology in degrees $-1$ and $-2$ matters, and it vanishes except for the class of the pair of pants (Section \ref{ssec:F}). To our knowledge no explicit geometric string vertex, of minimal area or hyperbolic, had been shown to represent Costello's class before; the combinatorial string vertices of \cite{CCT} represent it by construction, via the homotopy invariance of Maurer--Cartan moduli spaces (they use the framed version of the differential graded Lie algebra of \cite{Cos07}; Section \ref{ssec:lit}).

\textit{Relation to other work.} Moosavian and Pius \cite{MP17} use hyperbolic metrics with constant curvature to construct string vertices satisfying the master equation to first order. F{\i}rat and Valdes-Meller \cite{FVM} derive a Mirzakhani type recursion for the volumes of the hyperbolic vertices. C\u{a}ld\u{a}raru, Costello and Tu \cite{CCT} describe the string vertices of \cite{Cos05} homologically and combinatorially, via ribbon graphs and singular chains, and Tu \cite{Tu} expresses the fundamental classes of the Deligne--Mumford spaces through these vertices. F{\i}rat constructs the local coordinates of the hyperbolic three-string vertex \cite{Fi21}. The open-closed version is \cite{Cho}, see also \cite{HVZ,Fi24}; we treat only closed surfaces here.

\section{Weighted geometric chains}\label{sec:chains}
We use the geometric chain complexes $C^{g,\mathrm{mor}}_*$ and $GC^{\mathrm{mor}}_*$ of \cite[Definition on liftable maps and orbifold morphisms]{Yu}, with rational coefficients. A class $[Q,\hat f]$ consists of a compact oriented effective orbifold with corners $Q$ and an orbifold morphism $\hat f:Q\to X$ to an orbifold with corners $X$, up to isomorphism; $[-Q,\hat f]=-[Q,\hat f]$, a disconnected $Q$ is the sum of its components, and $\partial[Q,\hat f]=\sum_F[F,\hat f\circ i_Q|_F]$ over the components $F$ of the boundary $\partial Q$ in the sense of Joyce, with the outward normal first. Orbifolds with corners in the wide sense are proper \'etale Lie groupoids with corners up to Morita equivalence, not necessarily effective \cite{Yu}. For such an orbifold $Q$ the kernels of the local actions form a locally constant family of finite groups; on a connected component $Q_a$ we write $K_a$ for it, and $Q_{\red}$ for the effective orbifold obtained by dividing the local groups by their kernels \cite{Yu}.

\begin{definition}\label{def:weighted}
Let $Q$ be a compact oriented orbifold with corners in the wide sense, $X$ an effective orbifold with corners, and $\hat f:Q\to X$ an orbifold morphism. Near each point represent $\hat f$ by a pair $(\tilde f,\lambda)$ of a homomorphism $\lambda$ of local groups and a $\lambda$-equivariant smooth map \cite[Lemma on orbifold morphisms]{Yu}. We call $\hat f$ \emph{reducible} if $\lambda$ is trivial on the kernels $K_a$. Then $\hat f$ factors, uniquely up to $2$-isomorphism, as $Q\to Q_{\red}\xrightarrow{\hat f_{\red}}X$, and we put
\begin{equation*}
\langle Q,\hat f\rangle:=\sum_a\frac1{|K_a|}\,[Q_{a,\red},\hat f_{a,\red}]\in C^{g,\mathrm{mor}}_{\dim Q}(X).
\end{equation*}
If $X$ is an orbifold in the wide sense, we compose with the canonical morphism $X\to X_{\red}$ and regard $\langle Q,\hat f\rangle$ as a chain on $X_{\red}$.
\end{definition}
The factorization exists because the groupoid of $Q_{\red}$ is the quotient of that of $Q$ by the bundle of kernels, and $\hat f$ kills this bundle. For $X$ in the wide sense, $\hat f$ is reducible as a morphism to $X_{\red}$ whenever $\lambda$ maps the kernels of $Q$ into those of $X$. This is the convention of \cite{Yu} for non-effective sources, extended to morphisms. If $Q$ is effective, $\langle Q,\hat f\rangle=[Q,\hat f]$.

\begin{lemma}\label{lem:weighted}
\begin{enumerate}
\item If $\hat f$ is reducible, so is $\hat f\circ i_Q:\partial Q\to X$, and $\partial\langle Q,\hat f\rangle=\langle\partial Q,\hat f\circ i_Q\rangle$.
\item If $h:Q'\to Q$ is an isomorphism of orbifolds with corners in the wide sense and $\hat f$ is reducible, then $\langle Q',\hat f\circ h\rangle=\pm\langle Q,\hat f\rangle$, with the sign $+$ if $h$ preserves and $-$ if it reverses orientations.
\end{enumerate}
\end{lemma}
\begin{proof}
(1) Let $F$ be a component of $\partial Q$ lying over $Q_a$. The kernel $K_a$ acts trivially on $\partial Q$, so it is contained in the kernel $K_F$ of $F$. Conversely, let $g$ be an element of a local group $G$, acting linearly on a linear chart $\tilde U\subseteq R^n_k$ as in \cite[Lemma on doubling]{Yu}, which acts trivially near a point of a boundary hyperplane $\{y_i=0\}$. Then $g$ fixes this hyperplane pointwise and preserves the half space $\{y_i\ge0\}$. A finite order linear map with these properties is the identity, so $g$ acts trivially on $\tilde U$ and lies in the kernel. (Concretely, $g=\mathrm{id}+\alpha\otimes w$ with $\alpha$ vanishing on the hyperplane; finite order forces $\alpha(w)\in\{0,-2\}$, and $\alpha(w)=-2$ is a reflection, which does not preserve the half space.) Hence $K_F=K_a$, the morphism $\hat f\circ i_Q$ is reducible, and $(\partial Q)_{\red}=\partial(Q_{\red})$ \cite[convention after the definition of orbifolds]{Yu}. Here the boundary of a groupoid in the wide sense is taken chartwise, as in \cite[Convention on the boundary]{Yu}, the kernel acting trivially. Now (1) follows from the definition of $\partial$.
(2) An isomorphism identifies the kernels and induces an isomorphism $Q'_{\red}\cong Q_{\red}$ of the same orientation behaviour, compatible with the factorizations, and $[-Q,\hat f]=-[Q,\hat f]$.
\end{proof}

\section{Moduli spaces of bordered hyperbolic surfaces}\label{sec:moduli}
\subsection{Teichm\"uller spaces and mapping class groups}
Fix $L>0$. A \emph{type} is a compact oriented surface $S$, possibly disconnected, each of whose components $S_i$ has genus $g_i$, $n_i$ boundary components and $2g_i-2+n_i>0$. We write $\chi(S)=\sum_i(2-2g_i-n_i)<0$ and $n(S)=\sum_in_i$. Only finitely many types have given $(\chi,n)$, up to diffeomorphism.

The Teichm\"uller space $\T_L(S)$ is the set of equivalence classes of pairs $(X,h)$, where $X$ is a hyperbolic surface with totally geodesic boundary, every boundary component of length $L$, and $h:S\to X$ is an orientation preserving diffeomorphism; $(X,h)\sim(X',h')$ if $h'\circ h^{-1}$ is isotopic to an isometry, where isotopies need not fix the boundary. We have $\T_L(S)=\prod_i\T_L(S_i)$. Let $\Mod(S)$ be the group of isotopy classes of orientation preserving diffeomorphisms of $S$, which may permute the components and the boundary components of $S$; it acts on $\T_L(S)$ by $\phi\cdot[X,h]=[X,h\circ\phi^{-1}]$.

Since an orientation preserving isometry of a hyperbolic surface of negative Euler characteristic which is isotopic to the identity is the identity \cite[Ch.~7]{FM}, the isometry $h'\circ h^{-1}$ above is unique. Hence every point $x\in\T_L(S)$ determines a hyperbolic surface $X_x$, well defined up to unique isometry, with a marking $S\to X_x$ up to isotopy, and these surfaces form a smooth bundle over $\T_L(S)$ on which $\Mod(S)$ acts. For a boundary component $b$ of $S$ we write $b_x\subseteq\partial X_x$ for the corresponding boundary geodesic, oriented as the boundary of $X_x$, and $\rot^{b}_t$, $t\in\mathbb R/L\mathbb Z$, for the rotation of $b_x$ by arclength $t$ in this orientation. An \emph{interior curve} of $S$ is an isotopy class of essential simple closed curves in the interior of $S$ which are not boundary parallel; each is represented on $X_x$ by a unique simple closed geodesic in the interior, and every simple closed geodesic in the interior of $X_x$ arises in this way.

\subsection{Fenchel--Nielsen coordinates}\label{ssec:FN}
Let $S$ be connected, and let $P=\{\gamma_1,\dots,\gamma_N\}$, $N=3g-3+n$, be a pants decomposition of $S$ by interior curves. The Fenchel--Nielsen coordinates
\begin{equation*}
(\ell_1,\tau_1,\dots,\ell_N,\tau_N):\T_L(S)\to(\mathbb R_{>0}\times\mathbb R)^N
\end{equation*}
are a real analytic diffeomorphism, and the coordinates for two pants decompositions differ by a real analytic diffeomorphism \cite[Ch.~6]{Bu}, \cite[\S2]{Mir07}. Here $\ell_j$ is the length of the geodesic $\gamma_j$ and $\tau_j$ its twist parameter. We use the following convention for the sign of the twist.

\begin{definition}\label{def:left-twist}
Let $\gamma$ be a simple closed geodesic in the interior of $X$, oriented, so that $X$ has a side $A$ to the left of $\gamma$ and a side $B$ to its right near $\gamma$. The \emph{left twist deformation} of length $s\in\mathbb R$ cuts $X$ along $\gamma$ and reglues after translating the side $B$ by arclength $s$ in the direction of $\gamma$. Reversing the orientation of $\gamma$ exchanges $A$ and $B$ and the direction, so this does not depend on the orientation: each side moves to the left as seen from the other side.
\end{definition}
We choose the twist parameters $\tau_j$ so that the left twist deformation of length $s$ along $\gamma_j$ changes $\tau_j$ by $s$ and leaves the other coordinates unchanged \cite{Wo82}. The Dehn twist about $\gamma_j$ acts by $\tau_j\mapsto\tau_j\pm\ell_j$, the other coordinates being fixed.

\begin{convention}[orientation]\label{conv:orientation}
We orient $\T_L(S)$, $S$ connected, by $d\ell_1\wedge d\tau_1\wedge\dots\wedge d\ell_N\wedge d\tau_N$, and a disconnected type by the product orientation. Since $\dim\T_L(S_i)=6g_i-6+2n_i$ is even, the product orientation does not depend on the order of the components.
\end{convention}
By Wolpert's theorem \cite{Wo83,Wo85}, extended to bordered surfaces with fixed boundary lengths in \cite[\S2]{Mir07}, the form $\sum_jd\ell_j\wedge d\tau_j$ is the Weil--Petersson symplectic form up to a universal sign; in particular it does not depend on $P$. Hence the orientation does not depend on $P$. A mapping class $\phi$ maps Fenchel--Nielsen coordinates for $P$ to Fenchel--Nielsen coordinates for $\phi(P)$, up to changes $\tau_j\mapsto\tau_j+a_j\ell_j+c_j$ of the twists caused by the choice of seams, which do not change $d\ell_j\wedge d\tau_j$; so $\Mod(S)$ preserves the orientation. Only the independence of $P$ is used; the relation to the complex orientation of the moduli space of punctured curves is discussed in Remark \ref{rem:normalization}.

\subsection{The moduli orbifolds}
The group $\Mod(S)$ acts properly discontinuously and real analytically on $\T_L(S)$ \cite[Ch.~12]{FM}, \cite[Ch.~6]{Bu}, so the action groupoid
\begin{equation*}
\N_L(S):=\Mod(S)\ltimes\T_L(S)
\end{equation*}
is an oriented orbifold in the wide sense \cite{Yu}, of dimension $-3\chi(S)-n(S)$. Its underlying space is the set of isometry classes of hyperbolic surfaces of type $S$ with geodesic boundary of length $L$. The kernel of the action, the group $K(S)$ of mapping classes acting trivially on $\T_L(S)$, is the generic automorphism group. It is non-trivial for some types, for instance for the pair of pants, where $\T_L(S)$ is a point and $\Mod(S)\cong S_3$ permutes the boundary components, and for the one-holed torus (the elliptic involution). We put
\begin{equation*}
\N_{\chi,n}:=\bigsqcup_S\N_L(S),
\end{equation*}
the union over the diffeomorphism types $S$ with $\chi(S)=\chi$ and $n(S)=n$. It is an oriented orbifold in the wide sense of even dimension $-3\chi-n$, without boundary, and a finite disjoint union of connected components $\N_L(S)$.

\section{The vertex regions}\label{sec:A}
Put $\Lst:=2\,\mathrm{arcsinh}\,1$, so that $\sinh(\Lst/2)=1$.

\begin{lemma}\label{lem:hyperbolic}
Let $X$ be a hyperbolic surface with geodesic boundary, all components of negative Euler characteristic.
\begin{enumerate}
\item If two distinct simple closed geodesics $\gamma,\delta$ in the interior of $X$ intersect, then $\sinh(\ell_\gamma/2)\sinh(\ell_\delta/2)>1$. In particular, interior simple closed geodesics of length $\le\Lst$ are pairwise disjoint.
\item Every non-simple closed geodesic in $X$ has length $>4\,\mathrm{arcsinh}\,1=2\Lst$.
\item For $x\in\T_L(S)$ and $B>0$, only finitely many interior curves $\gamma$ satisfy $\ell_\gamma(x)\le B$. Every $x$ has a neighbourhood $U$ such that $\ell_\gamma(y)\ge\frac23\ell_\gamma(x)$ for all $y\in U$ and all interior curves $\gamma$.
\item (Bers) There is a constant $B(S,L)$ such that every $x\in\T_L(S)$ admits a pants decomposition all of whose curves have length $\le B(S,L)$.
\end{enumerate}
\end{lemma}
\begin{proof}
For closed surfaces, (1) is the collar lemma \cite[Ch.~4, collar theorem and its corollary on intersecting geodesics]{Bu}, (2) is \cite[Ch.~4]{Bu}, (3) is the discreteness of the length spectrum together with Wolpert's lemma $\ell_\gamma(y)\le K\,\ell_\gamma(x)$ for a $K$-quasiconformal map in the marking class \cite[Lemma 3.1]{Wo79}, For surfaces with geodesic boundary, (1)--(3) follow by applying these statements to the double $DX$ along $\partial X$, a closed hyperbolic surface: interior simple closed geodesics of $X$ are simple closed geodesics of $DX$ with the same intersections, closed geodesics of $X$ are closed geodesics of $DX$, and doubling is continuous, linear in Fenchel--Nielsen coordinates, $\T_L(S)\to\T(DS)$, so that the quasiconformal comparison of \cite{Wo79} applies. Statement (4) is Bers' theorem for compact surfaces with geodesic boundary, with a constant depending on the signature and an upper bound for the boundary lengths \cite[Ch.~5]{Bu}; it does not follow from the closed case by doubling, since a short pants decomposition of $DX$ need not consist of curves in $X$. 
\end{proof}

\begin{definition}\label{def:vertex}
For a type $S$ let $\T^V_L(S)\subseteq\T_L(S)$ be the set of $x$ such that every simple closed geodesic in the interior of $X_x$ has length $\ge L$, and put
\begin{equation*}
V_L(S):=\Mod(S)\ltimes\T^V_L(S)\subseteq\N_L(S),\qquad V_{\chi,n}:=\bigsqcup_SV_L(S)\subseteq\N_{\chi,n}.
\end{equation*}
For $x\in\T^V_L(S)$ let $\Gamma(x)$ be the set of interior simple closed geodesics of $X_x$ of length exactly $L$, regarded as isotopy classes of curves in $S$.
\end{definition}
If $L\le\Lst$, then by Lemma \ref{lem:hyperbolic}(2) the condition also holds for non-simple closed geodesics, and boundary geodesics have length $L$, so $x\in\T^V_L(S)$ if and only if the systole of $X_x$ is $\ge L$, as in \cite{CZ}. Since a disconnected surface satisfies the condition if and only if each component does, $\T^V_L(S)=\prod_i\T^V_L(S_i)$.

\begin{theorem}\label{thm:A}
Let $0<L\le\Lst$ and let $S$ be a type.
\begin{enumerate}
\item $\T^V_L(S)$ is a closed, $\Mod(S)$-invariant submanifold with corners of $\T_L(S)$ of full dimension. For $x\in\T^V_L(S)$ the curves in $\Gamma(x)$ are pairwise disjoint, and for every pants decomposition $P\supseteq\Gamma(x)$ there is a neighbourhood $U$ of $x$ with
\begin{equation*}
U\cap\T^V_L(S)=\{y\in U:\ \ell_\gamma(y)\ge L\ \text{for}\ \gamma\in\Gamma(x)\}.
\end{equation*}
So $x$ has depth $|\Gamma(x)|$, and its local boundary components correspond to the elements of $\Gamma(x)$.
\item $V_L(S)$ is a compact oriented orbifold with corners in the wide sense, of dimension $-3\chi(S)-n(S)$, and the inclusion $V_L(S)\to\N_L(S)$ is a morphism which is a local diffeomorphism onto its image.
\item The kernel of every connected component of $V_L(S)$ and of $\partial V_L(S)$ is $K(S)$. In particular the inclusion, and its composite with $i_{V}:\partial V_L(S)\to V_L(S)$, are reducible as morphisms to $\N_L(S)_{\red}$.
\end{enumerate}
\end{theorem}
\begin{proof}
For disconnected $S$, $\T^V_L(S)$ is the product of the sets for the components, and products of manifolds with corners are manifolds with corners; so in (1) we may assume $S$ connected.

(1) Closedness and invariance are clear, since $\T^V_L(S)$ is the intersection of the closed sets $\{\ell_\gamma\ge L\}$ and mapping classes permute the curves. Let $x\in\T^V_L(S)$. By Lemma \ref{lem:hyperbolic}(1) the curves in $\Gamma(x)$, of length $L\le\Lst$, are pairwise disjoint, so $\Gamma(x)$ extends to a pants decomposition $P$. By Lemma \ref{lem:hyperbolic}(3) with $B=2L$, the set $F$ of curves $\gamma$ with $\ell_\gamma(x)\le2L$ is finite, and on a neighbourhood $U_0$ of $x$ every $\gamma\notin F$ has $\ell_\gamma\ge\frac43L>L$. The curves in $F\setminus\Gamma(x)$ have $\ell_\gamma(x)>L$, so by continuity they have $\ell_\gamma>L$ on a smaller neighbourhood $U$. Hence on $U$ the set $\T^V_L(S)$ is cut out by the inequalities $\ell_\gamma\ge L$, $\gamma\in\Gamma(x)$. These lengths are among the Fenchel--Nielsen coordinates for $P$, so $U\cap\T^V_L(S)$ is, in these coordinates, an open subset of $\{\ell_\gamma\ge L,\ \gamma\in\Gamma(x)\}$, a product of $|\Gamma(x)|$ closed half-lines and a Euclidean space. Transition maps between such charts are restrictions of real analytic diffeomorphisms of $\T_L(S)$, hence smooth up to the corners, and they preserve depth, which is the number of active inequalities. The last statement holds because the local boundary components at $x$ are the hyperplanes $\{\ell_\gamma=L\}$.

(2) $\Mod(S)$ acts properly by diffeomorphisms preserving $\T^V_L(S)$ and its corner structure, so $V_L(S)$ is a proper \'etale groupoid with corners, an orbifold with corners in the wide sense. It is oriented by Convention \ref{conv:orientation}. The inclusion is induced by the equivariant inclusion $\T^V_L(S)\subseteq\T_L(S)$, with the identity of $\Mod(S)$.

For compactness of $|V_L(S)|$ we use Lemma \ref{lem:hyperbolic}(4). There are finitely many pants decompositions of $S$ up to $\Mod(S)$; fix representatives $P^1,\dots,P^m$. Let $x\in\T^V_L(S)$, and let $P$ be a pants decomposition with all lengths $\le B=B(S,L)$; then $P=\phi(P^r)$ for some $\phi$ and $r$. Replacing $x$ by $\phi^{-1}x$ and applying Dehn twists about the curves of $P^r$, which change the twists by multiples of the lengths, we may assume that in the coordinates for $P^r$ we have $\ell_j\in[L,B]$ (the lower bound because $x\in\T^V_L(S)$) and $\tau_j\in[0,\ell_j]$. So $|V_L(S)|$ is the image of the union of $m$ compact sets intersected with the closed set $\T^V_L(S)$, hence compact.

(3) The kernel of a component of $V_L(S)$ consists of the elements of $\Mod(S)$ acting trivially on an open subset of $\T_L(S)$. Such an element $g$ lies in the finite stabiliser of a point, so it has finite order. Averaging a Riemannian metric over $\langle g\rangle$ makes $g$ an isometry of the connected manifold $\T_L(S)$ which is the identity on an open set, hence $g$ acts trivially, i.e. $g\in K(S)$. The converse inclusion is clear. Note that $K(S)$ fixes every interior curve $\gamma$: if $g\gamma\neq\gamma$, then $\ell_\gamma=\ell_{g^{-1}\gamma}$ on $\T_L(S)$, which is false, since $\ell_\gamma$ can be made small with $\ell_{g^{-1}\gamma}$ bounded below if the two curves intersect, and the two lengths are independent Fenchel--Nielsen coordinates if they are disjoint. For a component of $\partial V_L(S)$ the kernel consists of the elements acting trivially on an open subset of a hypersurface $\{\ell_\gamma=L\}$. Such a $g$ fixes a point $y$ with $dg_y$ of finite order, equal to the identity on a hyperplane and of determinant $1$, since $g$ preserves the orientation. So $dg_y=\mathrm{id}$, and $g$, an isometry of an averaged metric, is the identity near $y$, hence on $\T_L(S)$. Conversely $K(S)$ acts trivially on these hypersurfaces. The inclusion maps $K(S)$ to the kernel $K(S)$ of $\N_L(S)$ and is therefore reducible; the same holds after composing with $i_V$.
\end{proof}

\begin{remark}\label{rem:A}
(1) The endpoint $L=\Lst$ is allowed because the inequality in Lemma \ref{lem:hyperbolic}(1) is strict.
(2) Faces need not be embedded: if two curves in $\Gamma(x)$ are exchanged by an element of the stabiliser of $x$, the two local boundary components at $|x|$ belong to the same component of $\partial V_L(S)$. The boundary of \cite{jo}, used in \cite{Yu}, handles this.
(3) For punctured surfaces and $L$ sufficiently small, Theorem \ref{thm:A}(1),(2) is \cite[Proposition 1.1]{JW}.
\end{remark}

\section{Sewing and the boundary of the vertex regions}\label{sec:B}
\subsection{The sewing spaces}\label{ssec:sewing}
Let $S'$ be a type with $n(S')\ge2$, and let $\Pi(S')$ be the set of unordered pairs $\{b,b'\}$ of distinct boundary components of $S'$. For $x\in\T_L(S')$ and $\{b,b'\}\in\Pi(S')$ let $\Phi_{b,b'}(x)$ be the set of orientation reversing isometries $\varphi:b_x\to b'_x$ or $\varphi:b'_x\to b_x$, modulo $\varphi\sim\varphi^{-1}$. The group $\mathbb R/L\mathbb Z$ acts on it by
\begin{equation}\label{eq:rotation}
t\cdot\varphi:=\rot^{b'}_t\circ\varphi\qquad(\varphi:b_x\to b'_x).
\end{equation}
This does not depend on the representative: since $\varphi^{-1}$ reverses orientations, $(\rot^{b'}_t\circ\varphi)^{-1}=\varphi^{-1}\circ\rot^{b'}_{-t}=\rot^{b}_t\circ\varphi^{-1}$, which is \eqref{eq:rotation} with the roles of $b,b'$ exchanged. The action is free and transitive.

The sets $\Phi_{b,b'}(x)$ form a principal $\mathbb R/L\mathbb Z$-bundle $\tilde\E_{b,b'}\to\T_L(S')$. Its smooth structure is given by smooth sections $s_b,s_{b'}$ of the boundary circles over $\T_L(S')$, for instance the feet of the seams used to define Fenchel--Nielsen twists: $\varphi\in\Phi_{b,b'}(x)$ corresponds to the unique $t$ with $\varphi(s_b(x))=\rot^{b'}_ts_{b'}(x)$. A mapping class maps the boundary circles of $X_x$ isometrically and orientation preservingly to those of $X_{\phi x}$, hence commutes with rotations, and acts on $\bigsqcup_{\{b,b'\}}\tilde\E_{b,b'}$ by bundle maps. Put
\begin{equation*}
\Pp(S'):=\Mod(S')\ltimes\big(\T_L(S')\times\Pi(S')\big),\quad\E(S'):=\Mod(S')\ltimes\bigsqcup_{\{b,b'\}\in\Pi(S')}\tilde\E_{b,b'}.
\end{equation*}
Then $p:\Pp(S')\to\N_L(S')$ is an orbifold fibre bundle with finite fibre $\Pi(S')$, and $\pi:\E(S')\to\Pp(S')$ is an orbifold circle bundle whose local groups act on the fibres by rotations, hence preserve the orientation given by \eqref{eq:rotation}. We orient them by these fibre orientations, the finite fibre positively. (In \cite{Yu} oriented orbifold bundles and their pullbacks are defined over effective bases; here the bases are orbifolds in the wide sense, and in this section we only use the total spaces. The corresponding operations on chains are defined in Section \ref{sec:master}.) We write $\Pp_{\chi,n}$ and $\E_{\chi,n}$ for the unions over the types with $(\chi(S'),n(S'))=(\chi,n)$.

\subsection{Cutting and gluing}\label{ssec:cut}
Let $S$ be a type, $\gamma$ an interior essential simple closed curve of $S$ which is not boundary parallel, and $S'$ the surface obtained by cutting $S$ along $\gamma$, with the two new boundary components $b,b'$ and the quotient map $q:S'\to S$. Every component of $S'$ has negative Euler characteristic, $\chi(S')=\chi(S)$ and $n(S')=n(S)+2$. Conversely, given $(S',\{b,b'\})$, gluing $b$ to $b'$ by an orientation reversing diffeomorphism gives $(S,\gamma)$. This is a bijection between the pairs $(S,[\gamma])$, $[\gamma]$ a $\Mod(S)$-orbit of such curves, and the pairs $(S',[\{b,b'\}])$, $[\{b,b'\}]\in\Pi(S')/\Mod(S')$, each up to diffeomorphism. The curve $\gamma$ is separating if and only if $b$ and $b'$ lie in different components of $S'$.

Let $\Stab(\gamma)\subseteq\Mod(S)$ and $\Stab\{b,b'\}\subseteq\Mod(S')$ be the stabilisers, which may reverse the orientation of $\gamma$, resp. exchange $b$ and $b'$. Cutting along $\gamma$ defines a homomorphism
\begin{equation*}
c:\Stab(\gamma)\to\Stab\{b,b'\},
\end{equation*}
which is surjective with kernel the infinite cyclic group generated by the Dehn twist $T_\gamma$. (The cutting homomorphism of \cite[\S3.6]{FM} concerns diffeomorphisms fixing the boundary and oriented curves; we give the argument in our setting.) \emph{Well defined:} a representative fixing $\gamma$ setwise cuts to a diffeomorphism of $S'$ preserving $\{b,b'\}$; two representatives isotopic in $S$ are isotopic through diffeomorphisms preserving $\gamma$, since the space of unoriented embedded curves in an isotopy class is contractible \cite{Gr}, so their cuts are isotopic. \emph{Surjective:} a diffeomorphism of $S'$ preserving $\{b,b'\}$ can be isotoped in a collar of $b'$ so that it is compatible with the gluing. \emph{Kernel:} let $G$ be the group of diffeomorphisms of $S'$ preserving $b$ and $b'$ individually and $G_c\subseteq G$ those compatible with the gluing; restriction to $b'$ gives a fibration $G_c\to G\to\mathrm{Diff}^+(S^1)$ up to homotopy. The identity component of $G$ is contractible \cite{ES}, so $\pi_1G=0$ and the kernel of $\pi_0G_c\to\pi_0G$ is the image of $\pi_1\mathrm{Diff}^+(S^1)\cong\mathbb Z$, generated by a twist in a collar of $b'$, which is $T_\gamma$ in $S$. An element isotopic to the identity in $S'$ preserves $b$ and $b'$ individually, so elements exchanging $b,b'$ do not contribute to the kernel. The subgroup $\langle T_\gamma\rangle$ is normal in $\Stab(\gamma)$, since $\phi T_\gamma\phi^{-1}=T_{\phi(\gamma)}$ for orientation preserving $\phi$.

Let $H_\gamma:=\{x\in\T_L(S):\ \ell_\gamma(x)=L\}$, a real analytic hypersurface. For $x\in H_\gamma$, cutting $X_x$ along the geodesic $\gamma_x$ gives a hyperbolic surface $Y$ with geodesic boundary of length $L$, marked by $S'$, together with the identification of the two copies of $\gamma_x$, an orientation reversing isometry $\varphi$ between $b_Y$ and $b'_Y$. This defines
\begin{equation*}
\cut:H_\gamma\to\tilde\E_{b,b'},\qquad x\mapsto(Y,\varphi).
\end{equation*}

\begin{lemma}\label{lem:cut}
\begin{enumerate}
\item $\cut$ is invariant under $T_\gamma$, equivariant with respect to $c$, and induces a real analytic diffeomorphism $H_\gamma/\langle T_\gamma\rangle\cong\tilde\E_{b,b'}$.
\item Let $P$ be a pants decomposition of $S$ containing $\gamma$, and $P'=q^{-1}(P\setminus\{\gamma\})$ the induced pants decomposition of $S'$. In the Fenchel--Nielsen coordinates $(\tau_\gamma,(\ell_j,\tau_j)_{j\ne\gamma})$ on $H_\gamma$ and $((\ell_j,\tau_j)_{j\ne\gamma},t)$ on $\tilde\E_{b,b'}$, with $t$ as in Section \ref{ssec:sewing} defined by the seam feet, $\cut$ preserves the coordinates $(\ell_j,\tau_j)_{j\neq\gamma}$, and $t\circ\cut=\tau_\gamma+c_0$ modulo $L$ for a constant $c_0$.
\end{enumerate}
\end{lemma}
\begin{proof}
(2) The coordinates $(\ell_j,\tau_j)$, $j\ne\gamma$, are defined by the pants of $P$ and their seams, which are not changed by cutting along $\gamma$; the pants adjacent to $\gamma$ become the pants of $P'$ adjacent to $b$ and $b'$, with the same seams. So these coordinates agree. Now fix $Y$ and vary the gluing. By Definition \ref{def:left-twist}, applied with $\gamma$ oriented as $b'$, the side of $b'$ lies to the left of $\gamma$ and the side of $b$ to its right, since a boundary circle is oriented with the surface on its left. Replacing $\varphi$ by $s\cdot\varphi=\rot^{b'}_s\circ\varphi$ glues each point of $b$ to a point displaced by $s$ in the direction of $\gamma$, i.e. it translates the right side by $s$ in the direction of $\gamma$. So it is the left twist deformation of length $s$, and it changes $\tau_\gamma$ by $s$ and $t$ by $s$. We take $s_b,s_{b'}$ to be the seam feet used in the definition of $\tau_\gamma$, with $\tau_\gamma$ normalized as a length; then $t$ and $\tau_\gamma$ modulo $L$ are both the signed displacement between these feet, and $t-\tau_\gamma$ is constant modulo $L$, also when the other coordinates vary.
(1) $T_\gamma$ changes the marking but not the hyperbolic surface or the curve $\gamma_x$, so $\cut$ is invariant; equivariance is clear from the definitions. In the coordinates of (2), $\cut$ is $(\tau_\gamma,y)\mapsto(y,\tau_\gamma+c_0)$, and $T_\gamma$ acts by $\tau_\gamma\mapsto\tau_\gamma\pm L$, which gives the diffeomorphism.
\end{proof}

\begin{lemma}[collar]\label{lem:collar}
Let $0<L\le\Lst$ and $x\in H_\gamma$, $\cut(x)=(Y,\varphi)$. Then $x\in\T^V_L(S)$ if and only if $Y\in\T^V_L(S')$, and in that case $\Gamma(x)=\{\gamma\}\sqcup q(\Gamma(Y))$.
\end{lemma}
\begin{proof}
The interior simple closed geodesics of $X_x$ other than $\gamma_x$ either are disjoint from $\gamma_x$ or cross it. Those disjoint from $\gamma_x$ are exactly the images of the interior simple closed geodesics of $Y$: a geodesic of $Y$ in the interior is not isotopic to $b$ or $b'$, since otherwise it would be isotopic to $\gamma$ in $X_x$ and hence equal to $\gamma_x$; and it is not isotopic to another boundary component of $Y$, since boundary geodesics of $Y$ other than $b,b'$ are boundary geodesics of $X_x$. Those crossing $\gamma_x$ have length $\ell$ with $\sinh(\ell/2)>1/\sinh(L/2)\ge1$, i.e. $\ell>\Lst\ge L$, by Lemma \ref{lem:hyperbolic}(1). Hence the conditions for $x$ and for $Y$ are equivalent, and the geodesics of length $L$ are $\gamma_x$ and the images of those of $Y$.
\end{proof}

By Lemma \ref{lem:cut}, the action groupoid $\Stab(\gamma)\ltimes H_\gamma$ is Morita equivalent to $\Stab\{b,b'\}\ltimes\tilde\E_{b,b'}$, the normal subgroup $\langle T_\gamma\rangle$ acting freely and properly. Composing the inverse of this Morita equivalence with the inclusion $\Stab(\gamma)\ltimes H_\gamma\to\Mod(S)\ltimes\T_L(S)$, and taking the union over $(S',[\{b,b'\}])$, we obtain the \emph{gluing morphism}
\begin{equation*}
\glue:\E_{\chi,n+2}\to\N_{\chi,n}.
\end{equation*}
On underlying spaces it maps $(Y,\{b,b'\},\varphi)$ to $Y/(z\sim\varphi(z))$. It depends on the choices of representatives $(S',\{b,b'\})$, of the identification $S\cong S'/\!\sim$ and of the Morita inverse only up to $2$-isomorphism.

\subsection{The boundary theorem}\label{ssec:B}
Let
\begin{equation*}
\E^V_{\chi,n+2}:=\pi^{-1}p^{-1}\big(V_{\chi,n+2}\big)\subseteq\E_{\chi,n+2}
\end{equation*}
be the restriction of the sewing bundle to the vertex region, a compact orbifold with corners in the wide sense, locally the product of $V_{\chi,n+2}$ with a finite set and a circle. We orient it base first: by the orientation of $V_{\chi,n+2}$, then the finite fibre positively, then the circle by \eqref{eq:rotation}. This is the base first convention of the pullback in \cite{Yu}. Since $\dim V_{\chi,n+2}$ is even, the fibre first orientation agrees with it.

\begin{theorem}\label{thm:B}
Let $0<L\le\Lst$, $\chi<0$ and $n\ge0$. There is an isomorphism of compact orbifolds with corners in the wide sense
\begin{equation*}
\Theta:\E^V_{\chi,n+2}\to\partial V_{\chi,n}
\end{equation*}
such that $i_{V}\circ\Theta$, followed by the inclusion $V_{\chi,n}\to\N_{\chi,n}$, is $2$-isomorphic to $\glue|_{\E^V}$. It reverses orientations, $\partial V_{\chi,n}$ carrying the boundary orientation (outward normal first).
\end{theorem}
\begin{proof}
\textit{Decomposition.} Fix a type $S$ with $(\chi(S),n(S))=(\chi,n)$. By Theorem \ref{thm:A}(1), the boundary of the manifold with corners $\T^V_L(S)$ in the sense of \cite{jo} is
\begin{equation*}
\partial\T^V_L(S)=\bigsqcup_\gamma\big(H_\gamma\cap\T^V_L(S)\big),
\end{equation*}
the union over all interior essential, non boundary parallel, simple closed curves $\gamma$, with $i$ the inclusion: a point of $\partial\T^V_L(S)$ is a pair $(x,\beta)$ with $\beta$ a local boundary component at $x$, i.e. a curve $\gamma\in\Gamma(x)$. Each $H_\gamma\cap\T^V_L(S)$ is a manifold with corners, with local boundary components at $x$ corresponding to $\Gamma(x)\setminus\{\gamma\}$. Since $\Mod(S)$ permutes the curves,
\begin{equation*}
\partial V_L(S)=\Mod(S)\ltimes\partial\T^V_L(S)\simeq\bigsqcup_{[\gamma]}\Stab(\gamma)\ltimes\big(H_\gamma\cap\T^V_L(S)\big),
\end{equation*}
the union over the $\Mod(S)$-orbits of curves, by the Morita equivalence between an action groupoid on a disjoint union of translates and the stabiliser of one of them. On the other side, $\E^V_{\chi,n+2}$ is the union over the types $S'$ with $(\chi(S'),n(S'))=(\chi,n+2)$ of $\Mod(S')\ltimes\bigsqcup_{\{b,b'\}}\tilde\E_{b,b'}|_{\T^V_L(S')}$, which is Morita equivalent to $\bigsqcup_{[\{b,b'\}]}\Stab\{b,b'\}\ltimes\tilde\E_{b,b'}|_{\T^V_L(S')}$. By Section \ref{ssec:cut} the pieces on the two sides are indexed by the same set.

\textit{Each piece.} By Lemma \ref{lem:collar}, $\cut$ maps $H_\gamma\cap\T^V_L(S)$ onto $\tilde\E_{b,b'}|_{\T^V_L(S')}$, and by Lemma \ref{lem:cut}(1) it induces a diffeomorphism
\begin{equation*}
\big(H_\gamma\cap\T^V_L(S)\big)/\langle T_\gamma\rangle\cong\tilde\E_{b,b'}|_{\T^V_L(S')}.
\end{equation*}
This is a diffeomorphism of manifolds with corners: in the coordinates of Lemma \ref{lem:cut}(2) both sides are cut out by $\ell_\delta\ge L$ for the curves $\delta\in\Gamma(Y)$, which are the same functions on both sides, by Lemma \ref{lem:collar} and Theorem \ref{thm:A}(1) applied to $x$ and to $Y$. It is equivariant with respect to $c$, which is surjective with kernel $\langle T_\gamma\rangle$. Hence it induces an isomorphism of groupoids with corners $\Stab(\gamma)\ltimes(H_\gamma\cap\T^V_L(S))\simeq\Stab\{b,b'\}\ltimes\tilde\E_{b,b'}|_{\T^V_L(S')}$, up to Morita equivalence. Let $\Theta$ be the union of the inverses of these isomorphisms. By construction, $i_V\circ\Theta$ composed with the inclusion into $\N_{\chi,n}$ is the gluing morphism.

\textit{Orientation.} Let $x\in H_\gamma\cap\T^V_L(S)$, and use the coordinates of Lemma \ref{lem:cut}(2). Let $\Omega':=\bigwedge_{j\neq\gamma}(d\ell_j\wedge d\tau_j)$, an even form which orients $\T_L(S')$ by Convention \ref{conv:orientation} (for disconnected $S$ or $S'$ we take the product with the forms of the other components, which commute with everything). The orientation of $\T_L(S)$ is $d\ell_\gamma\wedge d\tau_\gamma\wedge\Omega'$. Near $x$, $\T^V_L(S)\subseteq\{\ell_\gamma\ge L\}$, so the outward normal along $H_\gamma$ is $-\partial/\partial\ell_\gamma$, and the boundary orientation $o_\partial$ is determined by
\begin{equation*}
\big(d\ell_\gamma\wedge d\tau_\gamma\wedge\Omega'\big)\big(-\tfrac{\partial}{\partial\ell_\gamma},v_1,\dots,v_{d-1}\big)=-\big(d\tau_\gamma\wedge\Omega'\big)(v_1,\dots,v_{d-1}),
\end{equation*}
i.e. $o_\partial=-\,d\tau_\gamma\wedge\Omega'=-\,\Omega'\wedge d\tau_\gamma$. The base first orientation of $\E^V$ is $\Omega'\wedge dt$, and $\Theta^*$ maps $d\tau_\gamma$ to $dt$ by Lemma \ref{lem:cut}(2). Hence $\Theta$ reverses orientations. The local groups preserve both orientations, so this computation in the charts suffices.
\end{proof}

\begin{corollary}\label{cor:sewing}
Let $\mathbf V_{\chi,n}:=\langle V_{\chi,n},\iota\rangle\in C^{g,\mathrm{mor}}_{-3\chi-n}\big((\N_{\chi,n})_{\red}\big)$ for the inclusion $\iota$, and $\mathbf S_{\chi,n}:=\langle\E^V_{\chi,n+2},\glue\rangle\in C^{g,\mathrm{mor}}_{-3\chi-n-1}\big((\N_{\chi,n})_{\red}\big)$. Both are defined, and
\begin{equation*}
\partial\mathbf V_{\chi,n}=-\,\mathbf S_{\chi,n}
\end{equation*}
exactly in $C^{g,\mathrm{mor}}_*$, i.e. without passing to the quotient by degenerate chains.
\end{corollary}
\begin{proof}
By Theorem \ref{thm:A}(3), $\iota$ and $\iota\circ i_V$ are reducible. By Theorem \ref{thm:B}, $\glue|_{\E^V}\cong\iota\circ i_V\circ\Theta$ with $\Theta$ an isomorphism, so $\glue|_{\E^V}$ is reducible. Lemma \ref{lem:weighted}(1) gives $\partial\mathbf V_{\chi,n}=\langle\partial V_{\chi,n},\iota\circ i_V\rangle$, and Lemma \ref{lem:weighted}(2), applied to $\Theta$, which reverses orientations, gives $\langle\partial V_{\chi,n},\iota\circ i_V\rangle=-\langle\E^V_{\chi,n+2},\glue\rangle$.
\end{proof}

\begin{example}\label{ex:pants}
Let $\chi=-1$, $n=1$. The only type is the one-holed torus, and $\partial V_{-1,1}$ consists of the surfaces with a non-separating geodesic $\gamma$ of length $L$. Here $\Stab(\gamma)/\langle T_\gamma\rangle\cong\mathbb Z_2$, generated by the elliptic involution, which acts trivially on $\T_L(S)$; so $\partial V_{-1,1}$ is a circle (the twist modulo $L$) with trivial $\mathbb Z_2$-action, and $\partial\mathbf V_{-1,1}=\frac12[S^1_\partial,\cdot]$, where $S^1_\partial$ carries the boundary orientation $-d\tau_\gamma$. On the other side, $\N_{-1,3}=\N_L(\text{pants})$ is a point with group $S_3$; the stabiliser of a pair $\{b,b'\}$ is $\mathbb Z_2$, generated by the isometry $r$ exchanging $b$ and $b'$. Writing $\varphi_t$ for the gluing with $\varphi_t(s_b)=\rot^{b'}_ts_{b'}$, where $s_b,s_{b'}$ are the feet of the seams to the third boundary component, which $r$ exchanges, one finds $r\varphi_t^{-1}r^{-1}=\varphi_t$; so $r$ acts trivially on the fibre, and $\mathbf S_{-1,1}=\frac12[S^1_{\rot},\glue]$, with $S^1_{\rot}$ oriented by $dt=d\tau_\gamma$. Hence $\partial\mathbf V_{-1,1}=-\mathbf S_{-1,1}$, in accordance with Corollary \ref{cor:sewing}. Note that $\mathbf V_{-1,3}=\frac16[\mathrm{pt}]$.
\end{example}

\section{The master equation}\label{sec:master}
Throughout this section $0<L\le\Lst$.

\subsection{Chains over orbifolds in the wide sense}\label{ssec:wide}
The operator $\Delta$ involves pullbacks along the bundles $p$ and $\pi$ of Section \ref{ssec:sewing}, whose bases are orbifolds in the wide sense. It cannot be defined on $C^{g,\mathrm{mor}}_*((\N_{\chi,n})_{\red})$ by reducing these bundles: on reduced orbifolds the finite covering $\Pp_{-1,3}\to\N_{-1,3}$ becomes the identity of a point, whereas its degree is $3$, and a chain on $\N_{\red}$ need not lift to $\N$, nor uniquely. We therefore work with chains over $\N$ and reduce at the end.

\begin{definition}\label{def:wide-chains}
Let $N$ be an orbifold in the wide sense, without boundary. Let $C^w_k(N)$ be the $\mathbb Q$-vector space spanned by the isomorphism classes $[Q,\hat f]$ of pairs of a connected compact oriented $k$-dimensional orbifold with corners in the wide sense $Q$ and an orbifold morphism $\hat f:Q\to N$, modulo $[-Q,\hat f]=-[Q,\hat f]$; a disconnected $Q$ is the sum of its components. Isomorphisms are orientation preserving isomorphisms $h:Q'\to Q$ with $\hat f\circ h\cong\hat f'$. The boundary is $\partial[Q,\hat f]:=[\partial Q,\hat f\circ i_Q]$, with $\partial Q$ taken chartwise as in \cite[Convention on the boundary]{Yu}. A class is \emph{thin} if $\hat f\cong\hat g\circ\hat\pi$ through an orbifold with corners in the wide sense of dimension $<k$; this defines $S^w_*(N)$, $D^w_*(N)$ and $GC^w_*(N)$ as in \cite{Yu}. A class is \emph{kernel compatible} if, in local representations $(\tilde f,\lambda)$, $\lambda$ maps the kernel of $Q$ into the kernel of $N$; these classes span a subcomplex $C^{w,R}_*(N)$.
\end{definition}

The proof of \cite[Lemma on $C^g$]{Yu} applies verbatim: $C^w_*(N)$ is a chain complex, and a class admitting an orientation reversing automorphism $a$ with $\hat f\circ a\cong\hat f$ vanishes. $D^w_*$ is a subcomplex, as in \cite{Yu}. By Lemma \ref{lem:weighted}(1), kernel compatibility passes to faces, and the reduction
\begin{equation*}
\rho:C^{w,R}_*(N)\to C^{g,\mathrm{mor}}_*(N_{\red}),\qquad\rho[Q,\hat f]:=\langle Q,\hat f\rangle
\end{equation*}
of Definition \ref{def:weighted} is a chain map. In this notation $\mathbf V_{\chi,n}=\rho[V_{\chi,n},\iota]$ and $\mathbf S_{\chi,n}=\rho[\E^V_{\chi,n+2},\glue]$.

\begin{remark}\label{rem:wide-not-homology}
$C^w_*$ is an auxiliary complex and does not compute homology: for $N=BG$, $G$ a finite group, the $0$-chains $[BH,\hat\lambda]$ for homomorphisms $\lambda:H\to G$ (up to isomorphism of $H$ and conjugation in $G$) are independent cycles which are not boundaries: a compact connected oriented $1$-dimensional orbifold in the wide sense with non-empty boundary is isomorphic to $[0,1]\times BH$ (orientability excludes reflection points, and gerbes over an interval are trivial), and morphisms from it to $BG$ are constant in the interval direction. Homological statements are obtained by applying $\rho$, which lands in the complex $C^{g,\mathrm{mor}}_*(N_{\red})$ of \cite{Yu}, whose quotient by degenerate chains computes $H_*(|N|;\mathbb Q)$. We do not know whether $\rho$ maps $D^w_*\cap C^{w,R}_*$ into $D^{\mathrm{mor}}_*$; the identities below hold before passing to quotients, so this is not needed.
\end{remark}

\subsection{The operator $\Delta$}\label{ssec:Delta}
\begin{definition}\label{def:pullback-wide}
Let $q:Y\to N$ be an oriented orbifold fibre bundle whose fibre $F$ is a closed oriented manifold (a finite set or a circle below), the local groups acting on $F$ by orientation preserving diffeomorphisms. For $[Q,\hat f]\in C^w_k(N)$ put $q^*[Q,\hat f]:=[Q\times_NY,\widehat{\mathrm{pr}}_Y]$, the fibre product of groupoids, oriented \emph{fibre first}. For a morphism $\hat g:N\to N'$ put $\hat g_*[Q,\hat f]:=[Q,\hat g\circ\hat f]$.
\end{definition}
Locally, where $\hat f$ is represented by $(\tilde f,\lambda)$ on a chart $(\tilde U,G)$, the fibre product is the action groupoid of $G$ on $\tilde U\times F$, with $G$ acting on $F$ through $\lambda$. So it is a compact orbifold with corners in the wide sense, with $\partial(Q\times_NY)=\partial Q\times_NY$, and
\begin{equation}\label{eq:pullback-boundary}
\partial q^*=(-1)^{\dim F}q^*\partial,\qquad\partial\hat g_*=\hat g_*\partial,
\end{equation}
since $\partial(F\times Q)=(-1)^{\dim F}F\times\partial Q$ for the product orientation. As in \cite[pullback over orbifolds]{Yu}, a thin factorization $\hat f\cong\hat g\circ\hat\pi$ through $N'$ induces one of $q^*[Q,\hat f]$ through $N'\times_NY$, of dimension $\dim N'+\dim F$; so $q^*$ and $\hat g_*$ preserve $S^w$ and, by \eqref{eq:pullback-boundary}, $D^w$. If $[Q,\hat f]$ is kernel compatible, so is $q^*[Q,\hat f]$: an element of the kernel of $Q\times_NY$ lies in the kernel of $Q$ and acts trivially on $F$ through $\lambda$, and $\lambda$ of it lies in the kernel of $N$ and acts trivially on $F$, i.e. lies in the kernel of $Y$. Pushforward along $\hat g$ preserves kernel compatibility if $\hat g$ maps kernels into kernels.

\begin{lemma}\label{lem:glue-kernel}
$\glue:\E_{\chi,n+2}\to\N_{\chi,n}$ maps kernels into kernels.
\end{lemma}
\begin{proof}
On the piece indexed by $(S,[\gamma])$, $\glue$ is the Morita inverse of $\Stab(\gamma)\ltimes H_\gamma\simeq\Stab\{b,b'\}\ltimes\tilde\E_{b,b'}$ followed by the inclusion into $\Mod(S)\ltimes\T_L(S)$. An element of the kernel of $\Stab(\gamma)\ltimes H_\gamma$ acts trivially on the hypersurface $H_\gamma$, hence, by the argument in the proof of Theorem \ref{thm:A}(3), which applies to any open subset of $H_\gamma$, on $\T_L(S)$; so it lies in $K(S)$. Here the local groups correspond under the Morita equivalence: $\Stab(\gamma)_x\to\Stab\{b,b'\}_{\cut(x)}$ is injective because $\langle T_\gamma\rangle$ acts freely, and surjective because a lift of an element of the target fixes $x$ up to a power of $T_\gamma$.
\end{proof}

\begin{definition}\label{def:Delta}
The \emph{sewing operator} is
\begin{equation*}
\Delta:=\glue_*\circ\pi^*\circ p^*:\ C^w_k(\N_{\chi,n+2})\to C^w_{k+1}(\N_{\chi,n}),
\end{equation*}
with the fibre first orientations of Definition \ref{def:pullback-wide}, and $\Delta:=0$ on $C^w_*(\N_{\chi,n})$ for $n\le1$. On $\mathcal C:=\prod_{\chi<0,\,n\ge0}C^w_*(\N_{\chi,n})$ (for fixed $\chi$ only finitely many $n$ occur, and $\partial$, $\Delta$ act componentwise) we grade by codimension, $|[Q,\hat f]|:=\dim\N_{\chi,n}-\dim Q$, so that $\partial$ and $\Delta$ have degree $1$.
\end{definition}

\begin{proposition}\label{prop:Delta}
$\Delta$ preserves $S^w$, $D^w$ and $C^{w,R}$, and on $\mathcal C$
\begin{equation*}
\partial\Delta+\Delta\partial=0,\qquad\Delta^2=0 .
\end{equation*}
\end{proposition}
\begin{proof}
The first statement follows from the discussion after Definition \ref{def:pullback-wide} and Lemma \ref{lem:glue-kernel}. By \eqref{eq:pullback-boundary}, $p^*$ and $\glue_*$ commute with $\partial$ and $\pi^*$ anticommutes with it.

For $\Delta^2=0$, let $x=[Q,\hat f]\in C^w_*(\N_{\chi,n+4})$. Let $\E^{(2)}\to\N_{\chi,n+4}$ be the bundle whose fibre over a surface $Y$ consists of an ordered pair $(\pi_1,\pi_2)$ of disjoint unordered pairs of boundary components and gluing data $\varphi_1\in\Phi_{\pi_1}$, $\varphi_2\in\Phi_{\pi_2}$; on each type it is the action groupoid of $\Mod(S'')$ on the corresponding bundle over $\T_L(S'')$, with fibres finite sets times $(\mathbb R/L\mathbb Z)^2$. The boundary components of a sewn surface $Y/\varphi_1$ are those of $Y$ not in $\pi_1$, with the same metrics, so the pairs of boundary components of $Y/\varphi_1$ and their gluing data are those of $Y$ disjoint from $\pi_1$. On the piece indexed by $(S,[\gamma])$ this identifies the pullback of $\E_{\chi,n+2}$ along $\glue:\E_{\chi,n+4}\to\N_{\chi,n+2}$, i.e. along $\Stab\{b,b'\}\ltimes\tilde\E_{b,b'}\simeq\Stab(\gamma)\ltimes H_\gamma\to\Mod(S)\ltimes\T_L(S)$, with the part of $\E^{(2)}$ with $\pi_1=\{b,b'\}$. Hence
\begin{equation*}
\Delta^2x=\big[Q\times_{\N}\E^{(2)},\ \glue^{(2)}\circ\widehat{\mathrm{pr}}\big],
\end{equation*}
where $\glue^{(2)}:\E^{(2)}\to\N_{\chi,n}$ sews both pairs, and the source is oriented by the circle of $\varphi_2$, then that of $\varphi_1$, then $Q$. The involution $\sigma$ of $\E^{(2)}$ exchanging $(\pi_1,\varphi_1)$ and $(\pi_2,\varphi_2)$ is $\Mod$-equivariant and covers the identity of $\N_{\chi,n+4}$, so it induces an automorphism $\mathrm{id}\times_{\N}\sigma$ of the source. It reverses the orientation, since it exchanges the two circle coordinates. Moreover $\glue^{(2)}\circ\sigma\cong\glue^{(2)}$: on Teichm\"uller spaces both sides are induced by the inclusion of $H_{\gamma_1}\cap H_{\gamma_2}$, $\gamma_i$ the images of the pairs, into $\T_L(S)$, composed with the Morita equivalence with the stabiliser of the ordered pair $(\gamma_1,\gamma_2)$, resp. $(\gamma_2,\gamma_1)$, and these stabilisers coincide. (The two step description of $\glue\circ\glue$ agrees with this one step description because the iterated cutting homomorphism from the stabiliser of $(\gamma_1,\gamma_2)$ to that of $(\pi_1,\pi_2)$ is surjective with kernel $\langle T_{\gamma_1},T_{\gamma_2}\rangle$, by the argument of Section \ref{ssec:cut} applied twice; $\gamma_1,\gamma_2$ are disjoint, not isotopic and not boundary parallel, since no component of $S''$ is an annulus. The rotation orientation of the second circle is not affected by the first sewing, since the remaining boundary components keep their boundary orientations.) By the vanishing statement after Definition \ref{def:wide-chains}, applied componentwise ($\sigma$ either exchanges two components, whose classes then cancel, or preserves a component, whose class then vanishes), $\Delta^2x=0$.
\end{proof}

\begin{lemma}\label{lem:glue-sqcup}
Let $S'=S'_1\sqcup S_2$ be a type, with $S'_1$ possibly disconnected and $S_2$ possibly empty ($\N_L(\emptyset)=\N_{0,0}$), let $\{b,b'\}\in\Pi(S'_1)$, and let $S_1$ and $S=S_1\sqcup S_2$ be obtained by gluing $b$ to $b'$, with the image $\gamma$ of $b$; by Section \ref{ssec:cut}, $\gamma$ is an interior curve of $S_1$ which is not boundary parallel, and $S'$ is $S$ cut along $\gamma$. On the piece of $\{b,b'\}$ write $\glue_{S'}:=\iota_S\circ e^{-1}$ and $\glue_{S'_1}:=\iota_{S_1}\circ e_{S_1}^{-1}$, where $e=(c,\cut):\Stab_S(\gamma)\ltimes H_\gamma(S)\to\Stab_{S'}\{b,b'\}\ltimes\tilde\E_{b,b'}(S')$ is the Morita equivalence of Lemma \ref{lem:cut}, $e_{S_1}$ the one for $S_1$, and $\iota_S$, $\iota_{S_1}$ the inclusions of $\Stab(\gamma)\ltimes H_\gamma$ into $\Mod\ltimes\T_L$. Let
\begin{equation*}
j:(\Stab_{S'_1}\{b,b'\}\times\Mod(S_2))\ltimes(\tilde\E_{b,b'}(S'_1)\times\T_L(S_2))\to\Stab_{S'}\{b,b'\}\ltimes\tilde\E_{b,b'}(S')
\end{equation*}
be induced by the inclusion of groups and the equality $\tilde\E_{b,b'}(S')=\tilde\E_{b,b'}(S'_1)\times\T_L(S_2)$; its source is $(\Stab_{S'_1}\{b,b'\}\ltimes\tilde\E_{b,b'}(S'_1))\times(\Mod(S_2)\ltimes\T_L(S_2))$. Then, as morphisms to $\N_L(S)$,
\begin{equation*}
\glue_{S'}\circ j\ \cong\ \sqcup\circ(\glue_{S'_1}\times\mathrm{id}_{\N_L(S_2)}).
\end{equation*}
\end{lemma}
\begin{proof}
The circles $b,b'$ lie in $S'_1$ and $\gamma$ in $S_1$, so $H_\gamma(S)=H_\gamma(S_1)\times\T_L(S_2)$, $\cut_S=\cut_{S_1}\times\mathrm{id}$, $\Stab_{S_1}(\gamma)\times\Mod(S_2)\subseteq\Stab_S(\gamma)$, and the cutting homomorphism of $S$ restricts on this subgroup to $c_{S_1}\times\mathrm{id}$, since cutting along $\gamma$ does not change $S_2$. Put $e_1:=e_{S_1}\times\mathrm{id}_{\Mod(S_2)\ltimes\T_L(S_2)}$, and let $k$ be the inclusion of $(\Stab_{S_1}(\gamma)\times\Mod(S_2))\ltimes(H_\gamma(S_1)\times\T_L(S_2))$ into $\Stab_S(\gamma)\ltimes H_\gamma(S)$. Then $j\circ e_1=e\circ k$ as functors, hence $e^{-1}\circ j\cong e^{-1}\circ j\circ e_1\circ e_1^{-1}=e^{-1}\circ e\circ k\circ e_1^{-1}\cong k\circ e_1^{-1}$. Both $\iota_S\circ k$ and $\sqcup\circ(\iota_{S_1}\times\mathrm{id})$ are induced by the inclusion $\Stab_{S_1}(\gamma)\times\Mod(S_2)\subseteq\Mod(S)$ and the identity of $H_\gamma(S_1)\times\T_L(S_2)$, so they are equal, $\sqcup$ being induced by $\Mod(S_1)\times\Mod(S_2)\subseteq\Mod(S)$ and $\T_L(S)=\T_L(S_1)\times\T_L(S_2)$. Therefore
\begin{equation*}
\glue_{S'}\circ j=\iota_S\circ e^{-1}\circ j\cong\iota_S\circ k\circ e_1^{-1}\cong\sqcup\circ(\iota_{S_1}\times\mathrm{id})\circ(e^{-1}_{S_1}\times\mathrm{id})=\sqcup\circ(\glue_{S'_1}\times\mathrm{id}).
\end{equation*}
The choices in the definition of $\glue$ (representatives, the identification $S\cong S_1\sqcup S_2$, Morita inverses) change both sides only by $2$-isomorphisms. The map $j$ is in general a finite covering onto its image rather than an equivalence: if $S_2$ has a component of the same type as a component of $S'_1$, then $\Stab_{S'}\{b,b'\}$ is larger than $\Stab_{S'_1}\{b,b'\}\times\Mod(S_2)$; it is exactly the map occurring in the base change of Proposition \ref{prop:BV}.
\end{proof}

\subsection{The master equation}\label{ssec:C}
\begin{theorem}\label{thm:C}
Let $0<L\le\Lst$ and put
\begin{equation*}
\mathbb V:=\sum_{\chi<0,\ n\ge0}[V_{\chi,n},\iota]\ \in\ \mathcal C,
\end{equation*}
an element of degree $0$. Then
\begin{equation*}
(\partial+\Delta)\,\mathbb V=0\qquad\text{in }\mathcal C,
\end{equation*}
i.e.\ $\partial[V_{\chi,n},\iota]=-\Delta[V_{\chi,n+2},\iota]$ for all $\chi<0$, $n\ge0$, exactly in $C^w_*$. Moreover $\Delta[V_{\chi,n+2},\iota]=[\E^V_{\chi,n+2},\glue]$, all classes involved are kernel compatible, and applying $\rho$ gives $\partial\mathbf V_{\chi,n}=-\rho\Delta[V_{\chi,n+2},\iota]=-\mathbf S_{\chi,n}$ in $C^{g,\mathrm{mor}}_*((\N_{\chi,n})_{\red})$.
\end{theorem}
\begin{proof}
Since $V_{\chi,n+2}$ is a full suborbifold of $\N_{\chi,n+2}$, $p^*[V_{\chi,n+2},\iota]=[p^{-1}V_{\chi,n+2},\cdot]$ and $\pi^*p^*[V_{\chi,n+2},\iota]=[\E^V_{\chi,n+2},\cdot]$, oriented fibre first, which agrees with base first because $\dim V_{\chi,n+2}$ is even. So $\Delta[V_{\chi,n+2},\iota]=[\E^V_{\chi,n+2},\glue]$. By Theorem \ref{thm:B}, $\Theta:\E^V_{\chi,n+2}\to\partial V_{\chi,n}$ is an orientation reversing isomorphism with $\iota\circ i_V\circ\Theta\cong\glue$, so $[\E^V_{\chi,n+2},\glue]=-[\partial V_{\chi,n},\iota\circ i_V]=-\partial[V_{\chi,n},\iota]$. The degree of $[V_{\chi,n},\iota]$ is $0$ because $V_{\chi,n}$ has full dimension. Kernel compatibility is Theorem \ref{thm:A}(3) and Lemma \ref{lem:glue-kernel}, and the last statement is Corollary \ref{cor:sewing}.
\end{proof}

In the language of Batalin--Vilkovisky algebras, $\mathbb V$ plays the role of $e^{S/\hbar}$, a surface with Euler characteristic $\chi$ and $n$ boundary components carrying the power $\hbar^{-(\chi+n)/2}=\hbar^{\sum_i(g_i-1)}$, which $\partial$ preserves and $\Delta$ raises by one: $\Delta$ sews two boundary components of a possibly disconnected surface, which either lowers the number of components or creates a handle, and $(\partial+\Delta)\mathbb V=0$ is the exponential form $(d+\hbar\Delta)e^{S/\hbar}=0$ of the quantum master equation. The disjoint union of surfaces defines a product on $\mathcal C$, but $\mathbb V$ is not the exponential of its connected part for this product (Section \ref{ssec:connected}); we do not discuss the product further.

\subsection{The connected form}\label{ssec:connected}
For a connected type $t=(g,n)$ write $V_t:=V_L(S_{g,n})$ and $\N_t:=\N_L(S_{g,n})$. A type with two components of types $t_1\neq t_2$ has $\Mod(S_{t_1}\sqcup S_{t_2})=\Mod(S_{t_1})\times\Mod(S_{t_2})$, so its vertex region is $V_{t_1}\times V_{t_2}$; for $t_1=t_2=t$ the mapping class group is $(\Mod(S_t)\times\Mod(S_t))\rtimes\mathbb Z_2$, and the vertex region is
\begin{equation*}
\mathrm{Sq}(V_t):=(V_t\times V_t)\rtimes\mathbb Z_2,
\end{equation*}
$\mathbb Z_2$ exchanging the factors. Write $\Delta=\Delta_0+\Delta_\times$, where $\Delta_0$ sews two boundary components of the same component and $\Delta_\times$ two of different components. For $x\in C^w_*(\N_{t_1})$, $y\in C^w_*(\N_{t_2})$ define the bracket
\begin{equation*}
\{x,y\}:=\Delta_\times\big(\sqcup_*(x\times y)\big),
\end{equation*}
where $\sqcup:\N_{t_1}\times\N_{t_2}\to\N_L(S_{t_1}\sqcup S_{t_2})$ is the disjoint union, induced by the inclusion $\Mod(S_{t_1})\times\Mod(S_{t_2})\subseteq\Mod(S_{t_1}\sqcup S_{t_2})$ and $\T_L(S_{t_1}\sqcup S_{t_2})=\T_L(S_{t_1})\times\T_L(S_{t_2})$, an isomorphism if $t_1\ne t_2$ and the quotient by $\mathbb Z_2$ if $t_1=t_2$, and $x\times y=[Q_1\times Q_2,\hat f_1\times\hat f_2]$ with the product orientation.

\begin{corollary}\label{cor:connected}
For every connected type $(g,n)$ with $2g-2+n>0$, in $C^w_*(\N_L(S_{g,n}))$,
\begin{equation*}
\partial[V_{g,n}]+\Delta_0[V_{g-1,n+2}]+\sum_{\{t_1,t_2\},\ t_1\neq t_2}\{[V_{t_1}],[V_{t_2}]\}+\sum_t\Delta_\times[\mathrm{Sq}(V_t)]=0 .
\end{equation*}
The first sum runs over unordered pairs of distinct connected types $t_i=(g_i,n_i+1)$ with $g_1+g_2=g$, $n_1+n_2=n$, the second over the types $t$ with $(t,t)$ of this form, and terms with unstable or impossible types are zero.
\end{corollary}
\begin{proof}
Take the component of the identity of Theorem \ref{thm:C} in $C^w_*(\N_L(S_{g,n}))$. The pieces of $\E^V$ that glue to a connected surface are those with $S'$ connected and $\{b,b'\}$ in it, contributing $\Delta_0[V_{g-1,n+2}]$, and those with $S'$ having two components and $b,b'$ in different components. For $t_1\ne t_2$, the latter contribution is $\Delta_\times[V_{t_1}\times V_{t_2}]=\{[V_{t_1}],[V_{t_2}]\}$, and $\{[V_{t_1}],[V_{t_2}]\}=\{[V_{t_2}],[V_{t_1}]\}$ via the isomorphism exchanging the factors, which preserves orientations since all dimensions involved are even (in general $\{x,y\}=(-1)^{\dim x\dim y}\{y,x\}$). For $t_1=t_2=t$ it is $\Delta_\times[\mathrm{Sq}(V_t)]$.
\end{proof}

Summed over ordered pairs instead, the terms with $t_1\neq t_2$ acquire the factor $\frac12$, as in the usual form
\begin{equation*}\partial V+\Delta V+\tfrac12\{V,V\}=0\end{equation*}
of the master equation \cite{Zw,CZ}. The terms with $t_1=t_2$ do not have this form.

\begin{proposition}\label{prop:obstruction}
Let $t$ be a connected type with $\dim V_t>0$ for which $(t,t)$ occurs in Corollary \ref{cor:connected}, for instance $t=(1,1)$ for $(g,n)=(2,0)$ or $t=(0,4)$ for $(g,n)=(0,6)$. Then
\begin{equation*}
\rho\,\Delta_\times[\mathrm{Sq}(V_t)]\ \neq\ \tfrac12\,\rho\,\{[V_t],[V_t]\}\qquad\text{in }GC^{\mathrm{mor}}_*\big(\N_L(S_{g,n})_{\red}\big),
\end{equation*}
and the difference is not in $S^{\mathrm{mor}}_*$. Hence the identity of Corollary \ref{cor:connected} with $\frac12\{[V_t],[V_t]\}$ in place of $\Delta_\times[\mathrm{Sq}(V_t)]$ does not hold in $GC^{\mathrm{mor}}_*$, not even modulo $S^{\mathrm{mor}}_*$.
\end{proposition}
\begin{proof}
Let $E_t/\mathbb Z_2$ be the part of $\E^V$ over $\mathrm{Sq}(V_t)$ sewing one boundary component of each factor, and $E_t$ its pullback to $V_t\times V_t$. Let $s$ be the involution of $E_t$ exchanging the factors, $(Y_1,Y_2,b_1,b_2,\varphi)\mapsto(Y_2,Y_1,b_2,b_1,\varphi^{-1})$. It has fixed objects: take $Y_1=Y_2=Y\in V_t$ with trivial isotropy, $b_1=b_2=b$, and $\varphi$ the composite of the identification of the two copies of $b$ with a reflection of $b$; then $\varphi^{-1}=\varphi$ under this identification. Hence $s$ maps the component $C$ of $E_t$ through such an object to itself, and $q|_C:C\to q(C)$ is a connected double covering, $q^{-1}(q(C))=C$. It suffices to compare the coefficients of the basis classes of $C$ and of $q(C)$ in the two chains. Then $\{[V_t],[V_t]\}=[E_t,\glue\circ q]$ and $\Delta_\times[\mathrm{Sq}(V_t)]=[E_t/\mathbb Z_2,\glue]$, where $q:E_t\to E_t/\mathbb Z_2$ is the quotient. Since $\dim V_t>0$, the factor exchange acts effectively on $V_t\times V_t$ and on $E_t$, and it acts freely on a dense open subset; the kernels of $E_t$ and $E_t/\mathbb Z_2$ agree, so on these components $\rho$ gives multiples of two basis classes $[C_{\red},\glue\circ q]$ and $[q(C)_{\red},\glue]$ with coefficients $1/|K|$. By Theorem \ref{thm:B}, $\glue$ on $E_t/\mathbb Z_2$ is $\iota\circ i_V\circ\Theta$, which is injective on a dense open set (points of depth one with trivial isotropy), while $\glue\circ q$ is generically two to one onto its image. Neither class vanishes: if $a$ is an automorphism of $(E_t/\mathbb Z_2)_{\red}$ with $\glue\circ a\cong\glue$, then $|a|$ is the identity on a dense open set, hence $a\cong\mathrm{id}$ by \cite[Lemma on orbifold morphisms, (d)]{Yu}, which preserves orientation; on $C_{\red}$ such an $a$ is the identity or the deck transformation on a dense open set, and the deck transformation preserves orientation. The two classes are not isomorphic, and $\rho\Delta_\times[\mathrm{Sq}(V_t)]-\frac12\rho\{[V_t],[V_t]\}$ is a non-trivial combination of two distinct basis classes. Neither is thin, since both maps are local diffeomorphisms of rank $\dim E_t$ onto a hypersurface on a dense open set, and cannot factor through an orbifold of smaller dimension. As $S^{\mathrm{mor}}_*$ is spanned by thin basis classes, the difference is not in $S^{\mathrm{mor}}_*\supseteq D^{\mathrm{mor}}_*$. All other terms of Corollary \ref{cor:connected} are unchanged, which gives the last statement.
\end{proof}
For $\dim V_t=0$, i.e. $t=(0,3)$ and $(g,n)=(0,4)$, the two sides agree. Here $E_t$ is the action groupoid of $\mathbb Z_2\times\mathbb Z_2$ on a circle, the generators $r_1,r_2$ fixing the chosen boundary components and rotating them by half a turn, so the kernel is $\langle r_1r_2\rangle$; and $E_t/\mathbb Z_2$ is the action groupoid of the dihedral group of order $8$ on the same circle, where, as in Example \ref{ex:pants}, a lift $s$ of the exchange mapping seam feet to seam feet acts trivially, so the kernel $\{1,r_1r_2,s,sr_1r_2\}$ has order $4$. Both reduce to the same circle, and $\rho\,\Delta_\times[\mathrm{Sq}(V_{0,3})]=\frac14[C]=\frac12\rho\{[V_{0,3}],[V_{0,3}]\}$.

\subsection{The connected form with the gluing relation}\label{ssec:connected-G}
The obstruction of Proposition \ref{prop:obstruction} disappears if the relation $D$ of \cite{Yu} is replaced by the gluing relation of \cite{YuF}. For an effective orbifold $X$, possibly with corners and not compact, let $G^{\mathrm{mor}}_*(X)\subseteq C^{g,\mathrm{mor}}_*(X)$ be the span of the relations $[P,\hat f]-\sum_i[Q_i,\hat f\circ\hat\iota_i]$ for orbifold gluing decompositions $(P;(Q_i,\iota_i))$ of compact oriented effective orbifolds with corners, the pieces glued along faces with disjoint interiors \cite[Definition of orbifold gluing decompositions]{YuF}, and their boundaries. The morphisms $\hat\iota_i$ are determined by the decomposition \cite[Lemma on morphisms of embeddings]{YuF}. The quotient $C^{g,\mathrm{mor}}_*(X)/G^{\mathrm{mor}}_*(X)$, without the relation $D$, computes $H_*(|X|;\mathbb Q)$ \cite[Theorem on orbifold targets]{YuF}, and a covering of degree $d$ is equivalent to $d$ times its base \cite[Proposition on finite coverings]{YuF}: if $\pi:P'\to P$ is an orientation preserving orbifold covering of degree $d$ of compact oriented effective orbifolds with corners (locally $\bigsqcup_kG_k\ltimes\tilde U\to G\ltimes\tilde U$ with $\sum_k[G:G_k]=d$), then $[P',\hat f\circ\hat\pi]-d[P,\hat f]\in G^{\mathrm{mor}}_*(X)$.

\begin{theorem}\label{thm:connected-G}
For every connected type $(g,n)$ with $2g-2+n>0$, in $C^{g,\mathrm{mor}}_*(\N_L(S_{g,n})_{\red})/G^{\mathrm{mor}}_*$,
\begin{equation*}
\partial\,\rho[V_{g,n}]+\rho\,\Delta_0[V_{g-1,n+2}]+\frac12\sum_{(t_1,t_2)}\rho\,\{[V_{t_1}],[V_{t_2}]\}\equiv0,
\end{equation*}
the sum running over ordered pairs of connected types $t_i=(g_i,n_i+1)$ with $g_1+g_2=g$, $n_1+n_2=n$. This complex computes $H_*(|\N_L(S_{g,n})|;\mathbb Q)$.
\end{theorem}
\begin{proof}
Apply $\rho$ to Corollary \ref{cor:connected}, which holds in $C^w_*$, all classes involved being kernel compatible (Theorem \ref{thm:C}). The class $\{[V_t],[V_t]\}=[E_t,\glue\circ q]$ is kernel compatible as well, since $q$ maps the kernels of $E_t$ into those of $E_t/\mathbb Z_2$ and $\glue$ is kernel compatible on $E_t/\mathbb Z_2$. The terms with $t_1\ne t_2$ occur twice in the ordered sum, with the factor $\frac12$, and $\{[V_{t_1}],[V_{t_2}]\}=\{[V_{t_2}],[V_{t_1}]\}$. It remains to show $\rho\Delta_\times[\mathrm{Sq}(V_t)]\equiv\frac12\rho\{[V_t],[V_t]\}$. In the notation of the proof of Proposition \ref{prop:obstruction}, $\{[V_t],[V_t]\}=[E_t,\glue\circ q]$ and $\Delta_\times[\mathrm{Sq}(V_t)]=[E_t/\mathbb Z_2,\glue]$. Components of $E_t$ exchanged by $s$ map isomorphically onto the same component of $E_t/\mathbb Z_2$, and contribute exactly twice. If $\dim V_t>0$, the kernels of $E_t$ and $E_t/\mathbb Z_2$ agree, so on a component $C$ with $s(C)=C$, $q$ reduces to an orientation preserving orbifold covering $C_{\red}\to q(C)_{\red}$ of degree $2$ of compact oriented effective orbifolds with corners, both with the weight $1/|K|$, and the covering relation gives $\rho[C,\glue\circ q]\equiv2\rho[q(C),\glue]$. If $\dim V_t=0$, i.e. $t=(0,3)$, the two sides agree already before passing to the quotient (after Proposition \ref{prop:obstruction}). The last statement is \cite[Theorem on orbifold targets]{YuF}.
\end{proof}

\begin{remark}\label{rem:connected-G}
(1) Theorem \ref{thm:connected-G} does not contradict Proposition \ref{prop:obstruction}, which concerns the quotient by $D$: the two quotients compute the same homology, but since gluing relations need not lie in $D$, the identity on generators does not induce a chain map from the quotient by $G^{\mathrm{mor}}$ to the quotient by $D$. We do not know whether $D\subseteq G^{\mathrm{mor}}$.

(2) The identity is stated after $\rho$. A version in $C^w_*$ would need gluing relations for orbifolds in the wide sense, with pieces having the same kernels as the glued source in their interiors; for such relations the component $(g,n)=(0,4)$ would fail: for $t=(0,3)$ the kernels of $E_t$ and $E_t/\mathbb Z_2$ have orders $2$ and $4$, so $q$ is a covering in the direction of the kernels, and the integrals over the parts with kernel of a given order, without the weights, separate the two sides. With labelled boundary components (Remark \ref{rem:Costello}(3)) the same happens for zero dimensional faces, such as the pair of pants sewn to itself in type $(1,1)$; in positive dimension the exchange of the two new boundary components acts effectively, and the covering relation applies as in the proof above.

(3) In Section \ref{sec:wide-G}, with gluing relations whose pieces are injective only on interiors and with covering relations also in the direction of the kernels, $\Delta$ and the bracket descend to a quotient of $C^w_*$ which computes homology, and the identity holds there before $\rho$ (Corollary \ref{cor:exp}).
\end{remark}

\subsection{The de Rham form}\label{ssec:currents}
Integration over geometric chains sees covering degrees, and turns Corollary \ref{cor:connected} into the usual form. For an orbifold $N$ in the wide sense let $\Omega^*_c(N)$ be the compactly supported forms, i.e.\ invariant forms on charts, and $\mathcal D'_*(N)$ the currents, the continuous dual of $\Omega^*_c(N)$, with $(\partial T)(\omega)=T(d\omega)$. For a kernel compatible class put
\begin{equation*}
J[Q,\hat f](\omega):=\int_Q\hat f^*\omega:=\frac1{|K_Q|}\int_{Q_{\red}}\hat f_{\red}^*\omega ,
\end{equation*}
summed over components; this is integration over $\rho[Q,\hat f]$ in the sense of \cite[Remark on injectivity via Stokes]{Yu}. By Stokes' theorem on orbifolds with corners and Lemma \ref{lem:weighted}(1), $J\partial=\partial J$.

For a bundle $q:Y\to N$ as in Definition \ref{def:pullback-wide} let $q_!:\Omega^*_c(Y)\to\Omega^{*-\dim F}_c(N)$ be integration along the fibre, fibre first, computed in charts $H\ltimes(\tilde V\times F)$ (for a finite fibre, the sum over the fibre), and define $q^*T:=T\circ q_!$ on currents; for a proper morphism $\hat g$ put $\hat g_*T:=T\circ\hat g^*$. The maps $\glue$ and $\sqcup$ are proper: over a compact subset of $|\N_{\chi,n}|$ the systole is bounded below, so the cut surfaces lie in a compact set, only finitely many geodesics have length $L$, and the circle fibres are compact; $\sqcup$ is a finite covering. (The definition of $J$ makes sense for any morphism, with weight $1/|G|$ for the whole local group; we use it on kernel compatible classes.) Then $\Delta:=\glue_*\pi^*p^*$ and the bracket $\{T_1,T_2\}:=\Delta_\times\sqcup_*(T_1\times T_2)$ are defined on currents, and by Fubini's theorem in charts
\begin{equation*}
J\circ\Delta=\Delta\circ J,\qquad J\{x,y\}=\{Jx,Jy\}.
\end{equation*}
Call $q:Q\to Q'$ a covering of degree $d$ if it is locally the map $H\ltimes U\to G\ltimes U$ induced by the inclusion of a subgroup of index $d$, up to Morita equivalence; the examples used are $E_t\to E_t/\mathbb Z_2$ and $BK\to BG$. For such $q$, $\int_Qq^*\eta=d\int_{Q'}\eta$, with the weights of Definition \ref{def:weighted}; this includes the case where the deck group acts trivially, so that $q$ changes only the kernel. Applied to $E_t\to E_t/\mathbb Z_2$ this gives $J\Delta_\times[\mathrm{Sq}(V_t)]=\frac12\{JV_t,JV_t\}$ for every $t$.

\begin{corollary}\label{cor:currents}
Let $T_{g,n}:=J[V_{g,n}]\in\mathcal D'_*(\N_L(S_{g,n}))$ be the integration current of the vertex region. Then, for all $(g,n)$ with $2g-2+n>0$,
\begin{equation*}
\partial T_{g,n}+\Delta_0T_{g-1,n+2}+\frac12\sum_{(t_1,t_2)}\{T_{t_1},T_{t_2}\}=0,
\end{equation*}
the sum over ordered pairs as in Corollary \ref{cor:connected}. With $T:=\sum_{g,n}\hbar^gT_{g,n}$ this is $\partial T+\hbar\Delta_0T+\frac12\{T,T\}=0$.
\end{corollary}
\begin{proof}
Apply $J$ to Corollary \ref{cor:connected}, using $J\Delta_\times[\mathrm{Sq}(V_t)]=\frac12\{T_t,T_t\}$ and the symmetry of the bracket for $t_1\ne t_2$. For the last statement, $\Delta_0$ raises the genus by one and the bracket adds genera.
\end{proof}

This has the form of the master equation of \cite{CZ}, for the vertices in the model of bordered surfaces, with orientations given by Convention \ref{conv:orientation} and all factors from automorphisms determined by the orbifold structures. The paper \cite{CZ} does not fix the normalizations of $\Delta$, the bracket and $\mathcal V_{0,3}$ explicitly; the overall normalization there is fixed by requiring the Feynman chain to have degree one over $\overline{\mathcal M}_{g,n}$ \cite[Section 2.4]{CZ}. The following remark gives the dictionary to the conventions with labelled boundary components and to the complex orientation.

\begin{remark}\label{rem:normalization}
(1) \emph{Labelled boundary components.} Let $V^{\mathrm{lab}}_{g,n}$ be the vertex region in the moduli space with labelled boundary components, the quotient of $\T^V_L(S_{g,n})$ by the mapping classes fixing each boundary component, and $T^{\mathrm{lab}}_{g,n}:=J[V^{\mathrm{lab}}_{g,n}]$, an $S_n$-invariant current. The forgetful map $\pi$ is a covering of degree $n!$, so $T_{g,n}=\frac1{n!}\pi_*T^{\mathrm{lab}}_{g,n}$. The proof of Theorem \ref{thm:B} applies verbatim with labels: a non-separating face is the quotient of the sewing space of $V^{\mathrm{lab}}_{g-1,n+2}$ along the boundary components $n+1,n+2$ by the exchange of the two sides, and a separating face is the sewing space of $V^{\mathrm{lab}}_{g_1,I_1\sqcup\{*\}}\times V^{\mathrm{lab}}_{g_2,I_2\sqcup\{*\}}$ along $*$, for an unordered decomposition $\{(g_1,I_1),(g_2,I_2)\}$, divided by the exchange if the two parts coincide ($n=0$, $g_1=g_2$). Applying $J$ and the covering formula gives
\begin{multline*}
\partial T^{\mathrm{lab}}_{g,n}+\tfrac12\,\mathrm{sew}_{n+1,n+2}T^{\mathrm{lab}}_{g-1,n+2}\\
+\tfrac12\sum_{((g_1,I_1),(g_2,I_2))}\mathrm{sew}_*\big(T^{\mathrm{lab}}_{g_1,I_1\sqcup\{*\}}\times T^{\mathrm{lab}}_{g_2,I_2\sqcup\{*\}}\big)=0,
\end{multline*}
the sum over ordered decompositions with $g_1+g_2=g$, $I_1\sqcup I_2=\{1,\dots,n\}$, where $\mathrm{sew}$ sews the indicated boundary components with the twist circle oriented by rotation. In this form the factor $\frac12$ in front of the self-sewing term comes from the exchange of the two sides of a non-separating curve, and the factor $\frac12$ in front of the bracket from ordered decompositions.

(2) \emph{Automorphisms.} The currents $J[\,\cdot\,]$ include the weights $1/|K|$; for instance $T_{0,3}=\frac16\delta$, $T^{\mathrm{lab}}_{0,3}=\delta$, and $T_{1,1}=T^{\mathrm{lab}}_{1,1}$ is half the current of the underlying region, by the elliptic involution. These are the orbifold fundamental classes; in particular the integrals $\int_{V_{g,n}}\omega^{d/2}/(d/2)!$ are systolic Weil--Petersson volumes with the orbifold convention, which for $(1,1)$ differs by a factor $2$ from the convention used in some of the literature.

(3) \emph{Orientation.} As $L\to0$ the Fenchel--Nielsen coordinates identify $\T_L(S)$ with the Teichm\"uller space of the punctured surface, compatibly with $\Mod(S)$. Near a maximally degenerate point, with plumbing parameters $t_j$, $z_jw_j=t_j$, for the curves of $P$, the length $\ell_j$ is an increasing function of $|t_j|$, and increasing $\arg t_j$ glues the $w_j$-side rotated against the boundary orientation of the $z_j$-side, i.e. it is a right twist. So the complex orientation is $\bigwedge_j d\ell_j\wedge d\arg t_j=\bigwedge_j(-d\ell_j\wedge d\tau_j)$ up to positive factors, and our orientation is $(-1)^{3g-3+n}$ times the complex one on $\N_L(S_{g,n})$. Replacing the orientations of all $\N$ by the complex ones and those of the twist circles by $\arg t$, i.e. by right twists, changes $T_{g,n}$ by $(-1)^{3g-3+n}$ and $\Delta_0$ and the bracket by $-1$ each; since $3g-3+n$ drops by one under both operations, all identities of this section keep their form. [Medium confidence: the sign of $\arg t$ relative to the left twist depends on the plumbing convention $z w=t$ with $z,w$ holomorphic, and was checked only in this local model.]
\end{remark}

\section{Rigidity and comparison}\label{sec:comparison}
Throughout this section $0<L\le\Lst$. We use the compactly supported currents $\mathcal D'_{c,*}(N)$, the continuous dual of $\Omega^*(N)$, on which $q^*$, $\hat g_*$, $\Delta$ and the bracket are defined as in Section \ref{ssec:currents}, without properness assumptions. The currents $J[Q,\hat f]$ of chains have compact support.

\subsection{The algebraic structure on currents}\label{ssec:BV-currents}
Let $\mathcal A:=\prod_{\chi\le0,\,n\ge0}\mathcal D'_{c,*}(\N_{\chi,n})$, where $\N_{0,0}$ is a point, the empty surface, graded by codimension. For currents $T_i$ on $N_i$, $T_1\times T_2$ is the current on $N_1\times N_2$ with $(T_1\times T_2)(\alpha\times\beta)=T_1(\alpha)T_2(\beta)$, extended by continuity. Disjoint union of surfaces defines morphisms $\sqcup:\N_{\chi_1,n_1}\times\N_{\chi_2,n_2}\to\N_{\chi_1+\chi_2,n_1+n_2}$, and we put $T_1\cdot T_2:=\sqcup_*(T_1\times T_2)$. The connected part is $\mathfrak g:=\prod_{2g-2+n>0}\mathcal D'_{c,*}(\N_L(S_{g,n}))$.

\begin{proposition}\label{prop:BV}
\begin{enumerate}
\item $(\mathcal A,\cdot,\partial,\Delta)$ is a differential graded Batalin--Vilkovisky algebra in the sense of \cite[Definition of dg-BV algebraic structure]{Yu}, with unit the point $\N_{0,0}$, and its bracket $\{T_1,T_2\}=\Delta(T_1T_2)-(\Delta T_1)T_2-(-1)^{|T_1|}T_1\Delta T_2$ is, for connected $T_1,T_2$, the bracket of Section \ref{ssec:connected}.
\item $\mathfrak g$ is closed under $\partial$, $\Delta$ and $\{\,,\}$, and on $\mathfrak g$ the operator $\Delta$ equals $\Delta_0$. Hence $(\mathfrak g[[\hbar]],\partial+\hbar\Delta_0,[\,,])$, with $[a,b]:=(-1)^{|a|}\{a,b\}$, is a differential graded Lie algebra after the shift $|a|\mapsto|a|+1$, the bracket having degree one in $|\cdot|$.
\item The current $T=\sum_{g,n}\hbar^gT_{g,n}$ of Corollary \ref{cor:currents} is a Maurer--Cartan element: $(\partial+\hbar\Delta_0)T+\frac12\{T,T\}=0$.
\end{enumerate}
\end{proposition}
\begin{proof}
(1) The product is associative and graded commutative, since $\sqcup$ is associative and commutative up to $2$-isomorphism and all moduli spaces are even dimensional, so that the exchange of factors has the sign $(-1)^{|T_1||T_2|}$. $\partial$ is a derivation of the product of currents. The operator $\Delta$ satisfies $\partial\Delta+\Delta\partial=0$ and $\Delta^2=0$ by the arguments of Proposition \ref{prop:Delta}, read on forms: $q_!$ commutes with $d$ up to the sign $(-1)^{\dim F}$, and $\Delta^2$ is dual to fibre integration over $\E^{(2)}$, which is antisymmetric under $\sigma$. For the second order property we use base change: the pullback of $\E\to\Pp\to\N$ along $\sqcup$ decomposes as $(\E_1\times\N_2)\sqcup(\N_1\times\E_2)\sqcup\E_\times$, according to whether the pair lies in $X_1$, in $X_2$, or meets both, and fibre integration commutes with $\sqcup^*$ on each piece (also when $\sqcup$ is a covering); since $\sqcup$ is induced by the inclusion $\Mod(S_1)\times\Mod(S_2)\subseteq\Mod(S_1\sqcup S_2)$ and the identity of $\T_L$, the pullback of $\E$ along it is $(\Mod(S_1)\times\Mod(S_2))\ltimes\bigsqcup\tilde\E$, and on the orbit of a pair in $X_1$ it is Morita equivalent to the source of the map $j$ of Lemma \ref{lem:glue-sqcup}, the map to $\E$ being $j$; so by that lemma the gluing morphism is identified there with $\sqcup\circ(\glue\times\mathrm{id})$, and likewise for $X_2$, using $\sqcup\circ\mathrm{swap}\cong\sqcup$. With the source $S^1\times Q_1\times Q_2$ oriented fibre first, this gives
\begin{equation*}
\Delta(T_1T_2)=(\Delta T_1)T_2+(-1)^{|T_1|}T_1\Delta T_2+\{T_1,T_2\},
\end{equation*}
the sign coming from moving the circle past $Q_1$, with $(-1)^{\dim Q_1}=(-1)^{|T_1|}$ since moduli spaces are even dimensional; the third term is the sewing of one boundary component of each factor, i.e. the bracket of Section \ref{ssec:connected}. For three factors $a,b,c$, $\Delta(abc)$ decomposes into six pieces: $(\Delta a)bc$, $(-1)^{|a|}a(\Delta b)c$, $(-1)^{|a|+|b|}ab\Delta c$, and the cross terms $\{a,b\}c$, $(-1)^{|a|}a\{b,c\}$, $(-1)^{(|a|+1)|b|}b\{a,c\}$, where the last sign comes from moving the circle past $Q_b$ and exchanging $Q_a$ and $Q_b$; each piece is identified with the stated product by Lemma \ref{lem:glue-sqcup}, applied with $S'_1$ the union of the factors containing the pair ($S_a$, $S_a\sqcup S_b$, \dots) and $S_2$ the union of the remaining factors, together with the associativity and commutativity of $\sqcup$ up to $2$-isomorphism. Expanding the right hand side of the seven term identity of \cite[Definition of dg-BV algebraic structure]{Yu} with these formulas, each piece of the first kind occurs twice with sign $+$ and once with sign $-$, and each cross term once, so the identity holds.
(2) Sewing a pair within a connected surface, and sewing one boundary component of each of two connected surfaces, give connected surfaces, and a connected surface has no pairs in different components. In a dg BV algebra the bracket satisfies the shifted Jacobi identity and $\partial$ and $\Delta$ are derivations of it, the latter as a consequence of $\Delta^2=0$ \cite{Ko85}, so the restriction to $\mathfrak g$, which is closed under these operations, is a dg Lie algebra.
(3) This is Corollary \ref{cor:currents}.
\end{proof}

\subsection{Rigidity}\label{ssec:rigidity}
\begin{theorem}\label{thm:D}
Let $T'=\sum_{g,n}\hbar^gT'_{g,n}\in\mathfrak g[[\hbar]]$ be a Maurer--Cartan element such that each $T'_{g,n}$ is a compactly supported current of top dimension $6g-6+2n$ on $\N_L(S_{g,n})$, and $T'_{0,3}=T_{0,3}$. Then $T'=T$.
\end{theorem}
\begin{proof}
Order the connected types by $w:=2g-2+n\ge1$ and then by $g$. The component of the Maurer--Cartan equation on $\N_L(S_{g,n})$ reads
\begin{equation*}
\partial T'_{g,n}=-\Delta_0T'_{g-1,n+2}-\tfrac12\sum_{(t_1,t_2)}\{T'_{t_1},T'_{t_2}\},
\end{equation*}
and the right hand side involves only types $(g-1,n+2)$, with the same $w$ and smaller genus, and $t_1,t_2$ with $w(t_1)+w(t_2)=w$, $w(t_i)\ge1$, hence $w(t_i)<w$. By induction, the right hand sides for $T$ and $T'$ agree, so $D:=T'_{g,n}-T_{g,n}$ is a compactly supported current of top dimension with $\partial D=0$. For $(g,n)\neq(0,3)$, $\N_L(S_{g,n})$ is connected and not compact. On a chart $(\tilde U,G)$, the lift of $D$ is a $G$-invariant distribution times the orientation form, all of whose partial derivatives vanish since $\partial D=0$; so it is constant. Hence $D$ is a constant multiple of the fundamental current on each chart, the constant is locally constant on the connected space $|\N_L(S_{g,n})|$, and it vanishes because $D$ has compact support and $|\N_L(S_{g,n})|$ is not compact. So $D=0$.
\end{proof}

\begin{corollary}\label{cor:rigidity}
Let $V'_{g,n}\subseteq\N_L(S_{g,n})$, $2g-2+n>0$, be compact full dimensional suborbifolds with corners, given by closed $\Mod$-invariant subsets of $\T_L(S_{g,n})$ which are the closures of their interiors, with the induced orientation, and with $V'_{0,3}=V_{0,3}$, whose integration currents satisfy the master equation of Corollary \ref{cor:currents}. Then $V'_{g,n}=V_{g,n}$ for all $(g,n)$.
\end{corollary}
\begin{proof}
By Theorem \ref{thm:D}, $J[V'_{g,n}]=J[V_{g,n}]$, i.e. the indicator functions of $|V'_{g,n}|$ and $|V_{g,n}|$ agree almost everywhere. Two closed sets which are the closures of their interiors and agree up to measure zero are equal.
\end{proof}

So, in the model of bordered surfaces with geodesic boundary of length $L$, the master equation with the initial condition at the pair of pants determines the string vertices completely, among regions and even among currents of top dimension, such as smoothed vertices $f\cdot[\N_L(S_{g,n})]$ with $f$ smooth and compactly supported. The hyperbolic vertices are the solution.

\subsection{Comparison with Costello's uniqueness theorem}\label{ssec:Costello}
\begin{remark}\label{rem:Costello}
(1) In the framework of Costello \cite{Cos05}, which \cite[Section 2]{CZ} use for the homological existence and uniqueness of string vertices, string vertices are chains on spaces rationally homotopy equivalent, as orbifolds, to our $\N_L(S_{g,n})$, and solutions of the master equation with given initial condition are unique up to homotopy. The input is that the obstruction groups vanish: rationally $|\N_L(S_{g,n})|$ has the homology of $\Mod(S_{g,n})$, which vanishes above the virtual cohomological dimension, and this is smaller than $6g-6+2n$ except for $(0,3)$ \cite{Ha}. Theorem \ref{thm:D} is an analogue of this uniqueness in a different complex, and is compared with it in Section \ref{sec:Costello}: compactly supported currents vanish above the top dimension, so there is no room for homotopies, and its proof uses only that closed compactly supported top dimensional currents vanish, neither the dimension bound nor $L\le\Lst$.

(2) Section \ref{sec:Costello} compares the hyperbolic vertices with the homotopy class of Costello's solution. The route is through the moduli spaces $\Nt$ of surfaces with a marked point on each labelled boundary component, a strict model of Costello's framework in which rotating the marked points is a torus action with finite stabilisers whose quotient is $\N_L(S_{g,n})$, together with a residue map $\Phi$ to the complex of Section \ref{sec:wide-G}, compatible with sewing up to the sign $(-1)^{3g-3+n}$. Theorem \ref{thm:E}(1) shows that every solution in this model with Costello's initial condition integrates exactly to the currents $T_{g,n}$, up to this sign, and Theorem \ref{thm:E}(2) and Theorem \ref{thm:F}(1) that such solutions exist, are unique up to homotopy, and correspond to the class of $V_\hbar$ under a bijection of homotopy classes of Maurer--Cartan elements; these statements are unconditional inside the model $\Nt$. The transfer to Costello's own algebra for $C'_*(M^s)$, under which Costello's class corresponds to that of $V_\hbar$ (Theorem \ref{thm:E}(2), second part, and Theorem \ref{thm:F}(2)), assumes that the comparison map $h$ to his moduli spaces is continuous and relies on his computation of the homology of his algebra \cite[Proposition 9.0.2]{Cos05}, which rests on \cite[Proposition 5.1.1]{Cos05}, whose proof is sketched. Homotopies are those of Section \ref{ssec:Cos-framework}. Compactly supported currents still compute $H_*(|\N_L(S_{g,n})|;\mathbb R)$, by de Rham's theorem for orbifolds \cite{Sa56} and duality, so (1) and Theorem \ref{thm:D} are compatible; the difference is that currents vanish above the top dimension, whereas singular and geometric chains do not.

(3) In geometric chains, the operator $\Delta$ is defined on the auxiliary complex $\mathcal C$, which does not compute homology (Remark \ref{rem:wide-not-homology}), and not on the complex of \cite{Yu}, so the homotopy uniqueness of \cite{Cos05} cannot be applied to Theorem \ref{thm:C} directly. What survives in the complex of \cite{Yu} is Corollary \ref{cor:sewing}. This concerns the moduli spaces with unlabelled boundary components used here. With labelled boundary components and coinvariants, as in \cite[Section on applications]{Yu}, the finite coverings $p$ do not occur and the degree loss of Section \ref{ssec:wide} disappears; but the non-separating faces are then quotients by the exchange of the two branches, and the sewing term is a twofold covering of such a face, so an exact chain level identity meets the same obstruction as in Proposition \ref{prop:obstruction}. In either model, an exact chain level master equation in the usual form needs relations identifying a covering of degree $d$ with $d$ times its base, which are not part of \cite{Yu}. The gluing relation of \cite{YuF} provides them after reduction (Theorem \ref{thm:connected-G}), but $\Delta$ is not defined on $C^{g,\mathrm{mor}}_*((\N_{\chi,n})_{\red})$ (Section \ref{ssec:wide}), hence not on that quotient, so the homotopy uniqueness of \cite{Cos05} does not apply there either; Section \ref{sec:wide-G} removes this obstacle by putting the relations on $C^w_*$ itself (Remark \ref{rem:wide-G}(2)), and Section \ref{sec:Costello} carries out the comparison.
\end{remark}

\section{The master equation modulo gluing and covering relations}\label{sec:wide-G}
In this section we put relations on the auxiliary complex $C^w_*$ itself, so that it computes homology and $\Delta$ descends. The gluing relation is that of \cite{YuF}, with injectivity of the pieces required only on interiors, and with morphisms as part of the data; the covering relation identifies a covering of degree $d$ with $d$ times its base, also in the direction of the kernels, which realises in $C^w_*$ the weights $\frac1{|K|}$ of the convention of \cite{Yu} for non-effective sources, for instance $[BG,\hat f]\equiv\frac1{|G|}[\mathrm{pt},\hat f\circ\hat\pi]$.

\begin{definition}\label{def:Rw}
Let $N$ be an orbifold in the wide sense without boundary. A \emph{gluing datum} $\mathfrak d=(P;(Q_i,\hat\iota_i)_i)$ consists of compact oriented orbifolds with corners in the wide sense $P$, $Q_i$ of the same dimension and morphisms $\hat\iota_i:Q_i\to P$ such that
\begin{enumerate}
\item[(W1)] locally $\hat\iota_i$ is represented by pairs $(\tilde\iota,\lambda)$ with $\tilde\iota$ a smooth embedding of manifolds with corners;
\item[(W2)] at interior points $q$ of $Q_i$, $\lambda:H_q\to G_{\iota_i(q)}$ is an isomorphism;
\item[(W3)] $|\hat\iota_i|$ is injective on $\operatorname{int} Q_i$ and preserves orientations;
\item[(G2)] the images of the $|\hat\iota_i|$ cover $|P|$, and $|\hat\iota_i|(\operatorname{int} Q_i)\cap|\hat\iota_j|(Q_j)=\emptyset$ for $i\ne j$;
\item[(G3)] every face $F$ of $Q_i$ maps into the points of positive depth of $P$, or there are $j\neq i$, a face $F'$ of $Q_j$, an isomorphism $\hat\psi:F\to F'$ and a chosen $2$-isomorphism $\hat\iota_j\circ\hat\imath\circ\hat\psi\cong\hat\iota_i\circ\hat\imath$.
\end{enumerate}
A \emph{covering} of degree $d$ is an orientation preserving morphism $\hat\pi:P'\to P$, $P$ connected, such that every point of $P$ has a chart $G\ltimes\tilde U$ over which $\hat\pi$ is isomorphic to $\bigsqcup_kH_k\ltimes\tilde U\to G\ltimes\tilde U$, induced by $\mathrm{id}_{\tilde U}$, with $H_k\subseteq G$ and $\sum_k[G:H_k]=d$; the $H_k$ need not contain the kernel. Let $R^w_*(N)\subseteq C^w_*(N)$ be the span of the classes $[P,\hat f]-\sum_i[Q_i,\hat f\circ\hat\iota_i]$ and $[P',\hat f\circ\hat\pi]-d[P,\hat f]$, and of their boundaries. Let $\mathcal C^0_k:=\prod_{\chi\le0,\,n\ge0}C^w_{\dim\N_{\chi,n}-k}(\N_{\chi,n})$, where $\N_{0,0}$ is a point, the empty surface, and $R^w:=\prod_{\chi,n}R^w_*(\N_{\chi,n})$. Then $\mathcal C\subseteq\mathcal C^0$ is the part with $\chi<0$; Theorem \ref{thm:C} holds in $\mathcal C^0$, the $(0,0)$-component of $\mathbb V$ being $0$ and $\Delta=0$ on $C^w_*(\N_{0,0})$ by Definition \ref{def:Delta}.
\end{definition}

Coverings of degree $d$ in the direction of the kernels occur: for $H\subseteq G$ finite, $BH\to BG$ is one, of degree $[G:H]$, so $[BH,\hat f\circ\hat\pi]\equiv[G:H][BG,\hat f]$; in particular $[BG,\hat f]\equiv\frac1{|G|}[\mathrm{pt},\hat f\circ\hat\pi]$, and the cycles of Remark \ref{rem:wide-not-homology} become multiples of points.

\begin{lemma}\label{lem:Rw-stable}
Let $q:Y\to N$ be an oriented orbifold fibre bundle with closed fibre on which the local groups act preserving orientation, and $\hat g:N\to N'$ a morphism. Then $q^*$, $\hat g_*$ and the product $x\times y$ preserve $R^w$. Hence $R^w$ is preserved by $\partial$, $\Delta$, $\sqcup_*$ and the product of $\mathcal C^0$.
\end{lemma}

\begin{proof}
For a gluing datum, $(P\times_NY;(Q_i\times_NY,\hat\iota_i\times\mathrm{id}))$ is a gluing datum: locally $Q_i\times_NY$ is $H\ltimes(\tilde U\times F)$, with lifts $\tilde\iota\times\mathrm{id}_F$, which are embeddings; at an interior point $q$, $\lambda_q$ is an isomorphism compatible with the actions on the fibre, so the map of fibres $F/H_q\to F/G_{\iota(q)}$ is a bijection, and with (W3) the map is injective on $\operatorname{int}(Q_i\times_NY)=(\operatorname{int} Q_i)\times_NY$, the fibre having no boundary; at boundary points the map of fibres is onto, which gives (G2); outer faces map into the points of positive depth, and glued faces are glued by $\hat\psi\times\mathrm{id}$ with the $2$-isomorphism acted on by $\hat f$; orientations are fibre first. (At boundary points the pulled back piece need not be injective, which is why (W3) only asks for injectivity on interiors.) At an interior point the local group of the pulled back piece is $\mathrm{Stab}_{H_q}(y)\cong\mathrm{Stab}_{G_{\iota(q)}}(y)$, which is (W2). A covering pulls back to a sum of coverings $P'\times_NY\to P\times_NY$ over the components, of the same degree by the double coset formula $\sum_{H_kgS}[S:S\cap g^{-1}H_kg]=[G:H_k]$ for the stabiliser $S$ of a point of the fibre. Push forward does not change the data. Products of data with a source $Q_2$ are data, with pieces $Q_i\times Q_2$, faces $F\times Q_2$ glued by $\hat\psi\times\mathrm{id}$ and outer faces $Q_i\times\partial Q_2$, and $\hat\pi\times\mathrm{id}$ is a covering. The boundary parts follow from \eqref{eq:pullback-boundary} and the Leibniz rule. Since the components of a product in $\mathcal C^0$ are finite sums, $R^w$ is preserved componentwise.
\end{proof}

\begin{proposition}\label{prop:BV-chains}
$(\mathcal C^0,\cdot,\partial,\Delta)$, with $x\cdot y:=\sqcup_*(x\times y)$ and unit $[\mathrm{pt}]\in C^w_0(\N_{0,0})$, is a differential graded BV algebra already at chain level, and $R^w$ is an ideal for $\cdot$, $\partial$ and $\Delta$, so $\mathcal C^0/R^w$ is a differential graded BV algebra. Its bracket is the bracket of Section \ref{ssec:connected} on connected types.
\end{proposition}

\begin{proof}
The proof of Proposition \ref{prop:BV}(1) applies to chains: associativity and graded commutativity hold at chain level since $\sqcup$ is associative and commutative up to $2$-isomorphism and moduli spaces are even dimensional; $\partial$ is a derivation since $\partial(Q_1\times Q_2)=\partial Q_1\times Q_2\sqcup(-1)^{\dim Q_1}Q_1\times\partial Q_2$; $\partial\Delta+\Delta\partial=0$ and $\Delta^2=0$ by Proposition \ref{prop:Delta}. For the second order property, $p^*\sqcup_*$ and $\pi^*$ are computed by base change, which is the pasting of $2$-pullbacks of groupoids; the pullback of $\Pp$ along $\sqcup$ splits into the open and closed parts where the pair lies in $X_1$, in $X_2$ or meets both, so the source of $\Delta(xy)$ is a disjoint union of three parts, and by Lemma \ref{lem:glue-sqcup}, together with $\sqcup\circ\mathrm{swap}\cong\sqcup$ for the part in $X_2$, the first two are $(\Delta x)y$ and $(-1)^{|x|}x\Delta y$; the third is $\{x,y\}$ by definition. The six pieces for three factors are identified in the same way, and the seven term identity holds in $C^w_*$. The second statement is Lemma \ref{lem:Rw-stable}.
\end{proof}

\begin{theorem}\label{thm:Tw}
For every orbifold $N$ in the wide sense without boundary, the comparison map $\sigma\mapsto[\Delta_{\mathrm{top}}^n,\hat\sigma]$ from the complex $C^{\mathrm{sm}}_*(N)$ of morphisms $\hat\sigma:\Delta_{\mathrm{top}}^n\to N$, composed with the inverse of $H_*(C^{\mathrm{sm}}_*(N))\cong H_*(|N|;\mathbb Q)$ (established in the proof), induces an isomorphism $H_*(|N|;\mathbb Q)\cong H_*(C^w_*(N)/R^w_*(N))$. Here $\Delta_{\mathrm{top}}^n$ denotes the standard simplex, to distinguish it from the sewing operator $\Delta$.
\end{theorem}

\begin{proof}
We follow the proof of \cite[Theorem on orbifold targets]{YuF}.

\emph{Singular side.} Let $C^{\mathrm{sm}}_*(N)$ be the complex spanned by the $2$-isomorphism classes of morphisms $\hat\sigma:\Delta_{\mathrm{top}}^n\to N$, with faces $\hat\sigma\circ d_i$; for an open subgroupoid $O$, $C^{\mathrm{sm}}_*(O)\to C^{\mathrm{sm}}_*(N)$ is injective onto the span of a subset of the basis, an open subgroupoid being full. The proof of \cite[Proposition on chains given by liftable maps and orbifold morphisms, part (A)]{Yu} uses only the singular parts of items (3), (6), (7) and the argument of item (8), together with a good cover; these apply to $N$. As good cover take the one constructed for $N_{\red}$ in \cite[proof of the Proposition on good covers of orbifolds]{Yu}, with charts chosen as reductions of local action groupoids $\tilde G_x\ltimes\tilde U_x$ of $N$; over a set $W/H$ the groupoid $N$ is Morita equivalent to $G\ltimes W$, $G$ the preimage of $H$, with kernel $K_N$, and $G$ fixes the centre $c$. Indeed $N$ over $U_x$ is Morita equivalent to $\tilde G_x\ltimes\tilde U_x$, and the preimage of $W/H$ is the disjoint union of the translates of $W$, so restricting to the stabiliser of $W$ gives $G\ltimes W$; effectiveness is used in \cite{Yu} only for $N_{\red}$. Morphisms $\Delta_{\mathrm{top}}^n\to BG$ are trivial for any finite $G$, effective or not, since principal bundles over $\Delta_{\mathrm{top}}^n$ are trivial; this is all that the singular part of item (6) uses. In detail: over $W/H$, a morphism $\Delta_{\mathrm{top}}^n\to G\ltimes W$ is, up to $2$-isomorphism, a smooth simplex in $W$ together with the trivial $G$-bundle, and $2$-isomorphisms act through $G$; so $C^{\mathrm{sm}}_*(G\ltimes W)$ is the complex of $G$-coinvariants of the smooth singular chains of $W$, which over $\mathbb Q$ computes $H_*(W/H;\mathbb Q)=H_*(\mathrm{pt};\mathbb Q)$, since $W$ is contractible and $K_N$ acts trivially on $W$; this is items (3) and (6) for $N$. Item (7), Mayer--Vietoris for $C^{\mathrm{sm}}$, uses only small simplices and subdivision, which commute with the $G$-action, and the de Rham pairing below holds since $\Omega(N)=\Omega(N_{\red})$.

\emph{Integration.} For an $n$-form $\omega$ on $N$ (equivalently on $N_{\red}$) put $\int_{[Q,\hat f]}\omega:=\int_Q\hat f^*\omega$, with the weights $1/|G|$ of the full local groups; on a component, this is $\frac1{|K_Q|}\int_{Q_{\red}}\hat f^*\omega$, since $\hat f^*\omega$ is invariant. Stokes' theorem holds since faces have the same kernel (proof of Lemma \ref{lem:weighted}(1)). On a gluing datum, $\hat\iota_i$ restricts on $\operatorname{int} Q_i$ to an isomorphism onto an open subgroupoid of $P$ by (W1)--(W3), so the weights agree and the argument of \cite[Proposition on the axioms for orbifold gluing, (A5)]{YuF} applies with the interior points whose isotropy is the kernel. On a covering, locally $\sum_k\frac1{|H_k|}\int_{\tilde U}=\frac d{|G|}\int_{\tilde U}$. So integration vanishes on $R^w$, and with the de Rham pairing \cite[Lemma on the de Rham pairing]{YuF}, which holds for $N$ since $\Omega(N)=\Omega(N_{\red})$, the argument of \cite{YuF} gives injectivity, after reducing by compact supports to a finite union of members of the good cover; relations with image in an open subgroupoid are relations there, as morphisms restrict and the data do not involve the map.

\emph{Contractible pieces.} Let $\hat C$ be the strict morphism induced by the $G$-equivariant geodesic contraction of $G\ltimes W$ to $c$, and $z$ a cycle of degree $k\ge1$; then $z\equiv(\hat C_1)_*z$, since the prism operator of $\hat C$ maps $R^w$ into $R^w$: the product of a gluing datum or a covering with $[0,1]$ is again one, the relations do not involve the map, and push forward along $\hat C\circ(\hat f\times\mathrm{id})$ preserves $R^w$ by Lemma \ref{lem:Rw-stable}. For a term $[Q,\hat f]$ let $K_Q$ be the kernel of a component. On convex linear charts $G_Q\ltimes\tilde U$, on which $G_Q/K_Q$ acts linearly, $\hat f$ is represented by $(\tilde f,\lambda)$ \cite[Lemma on orbifold morphisms, (a)]{Yu}, whose proof does not use effectiveness, and $\hat C_1\hat f$ by $(\mathrm{const}_c,\lambda)$. Steps 2 and 3 of \cite[proof of the Proposition on the point for orbifold sources]{YuF} for $G_Q/K_Q$ are gluing relations, with pieces the action groupoids $K_Q\ltimes(\sigma\cap\tilde V)$ and $\hat\iota$ induced by $K_Q\subseteq G_Q$; (W2) holds since the open cones consist of free points of $G_Q/K_Q$, and the faces are glued by lifts of elements of $G_Q/K_Q$, which normalise $K_Q$. On a piece, $1\ltimes\sigma\to K_Q\ltimes\sigma$ is a covering of degree $|K_Q|$, so $[K_Q\ltimes(\sigma\cap\tilde V),\hat C_1\hat f\hat\iota]\equiv\frac1{|K_Q|}\hat c_*[\sigma,\mathrm{const}]$, which lies in $R^w$ by \cite[Corollary on the point]{YuF}. In degree $0$, $[BK,(\lambda,x)]\equiv\frac1{|K|}[\mathrm{pt},x]$, and the weighted count of points $\sum\pm\frac1{|K|}$ vanishes on boundaries (compact oriented $1$-dimensional sources reduce to intervals and circles, with the kernel at the end points) and on $R^w$; so $H_0=\mathbb Q$.

\emph{Mayer--Vietoris and conclusion.} Cutting along a function on $|N|$ preserves the kernels of the pieces, maps gluing data to gluing data and coverings to sums of coverings, by the double coset formula on the components of the slices. Explicitly, for a covering $\pi:P'\to P$ of degree $d$ and a function $h$ with $(h\circ\hat f,c)$ regular, the lifts of $h\circ\hat f$ and $h\circ\hat f\circ\hat\pi$ to a chart $\tilde U$ are the same function, so $c$ is regular for both, and $\cut_\pm([P',\hat f\hat\pi]-d[P,\hat f])=\sum_C([C',\hat f\hat\pi]-d_C[C,\hat f])$ over the components $C$ of $P_\pm$ with preimages $C'$, each a covering of degree $d_C=d$ on the connected part over which it is formed, and likewise for $\Lambda$ with the components of the level set; the lemmas on cut operators and on cuts of relations then hold by linearity. With this, the Mayer--Vietoris step and the conclusion of \cite{YuF} apply.
\end{proof}

\begin{corollary}\label{cor:exp}
Let $V:=\sum_t[V_t,\iota]$, the sum over connected types. In $\mathcal C^0/R^w$, $\mathbb V\equiv e^V-1$, and
\begin{equation*}
\partial V+\Delta V+\tfrac12\{V,V\}\equiv0 .
\end{equation*}
Hence $V_\hbar:=\sum\hbar^g[V_{g,n}]$ is a Maurer--Cartan element of the differential graded Lie algebra $(\mathfrak g^w[[\hbar]],\partial+\hbar\Delta_0,[\,,])$, $\mathfrak g^w:=\prod_{2g-2+n>0}C^w_*(\N_L(S_{g,n}))/R^w$, with $[a,b]:=(-1)^{|a|}\{a,b\}$ (a differential graded Lie algebra after the shift $|a|\mapsto|a|+1$; Maurer--Cartan elements have degree $0$ in the unshifted grading by codimension, and the powers of $\hbar$ are determined by the genus as in Proposition \ref{prop:BV}(2)); by Theorem \ref{thm:Tw}, $\mathfrak g^w$ computes $\prod H_*(|\N_L(S_{g,n})|;\mathbb Q)$. Integration over forms is a morphism of differential graded BV algebras from $\mathcal C^0/R^w$ to the algebra of Section \ref{ssec:BV-currents} mapping $V_\hbar$ to the current $T$ of Proposition \ref{prop:BV}(3).
\end{corollary}

\begin{proof}
For a type with $m_i$ components of type $t_i$, the vertex region is $\prod_i V_{t_i}^{m_i}\rtimes S_{m_i}$, and $\prod_iV_{t_i}^{m_i}$ maps to it by a covering of degree $\prod_im_i!$, which is the identity on charts and hence orientation preserving. In $\mathcal C^0$, $\prod_i[V_{t_i},\iota]^{m_i}$ is the class of this covering composed with the inclusion, the order of the factors being irrelevant since dimensions are even, so the covering relation, applied on each connected component, gives $\prod_i[V_{t_i},\iota]^{m_i}/m_i!\equiv[V(\prod t_i^{m_i}),\iota]$, and $\mathbb V\equiv e^V-1$; here $e^V$ is defined since $V^m$ is supported on $\chi\le-m$. By Theorem \ref{thm:C}, $(\partial+\Delta)(e^V-1)\equiv0$; since $\partial1=\Delta1=0$ and $\Delta$ has order two, $(\partial+\Delta)e^V=(\partial V+\Delta V+\frac12\{V,V\})e^V$, and $e^V$ is invertible. The connected part is closed under the operations as in Proposition \ref{prop:BV}(2), since a relation with disconnected $P$ splits over the components of $P$ (each piece of a gluing datum has image in one component, and the partner faces of (G3) have the same image), so that $R^w$ is the direct sum of its parts over the components $\N_L(S)$; on it $\Delta=\Delta_0$. Integration vanishes on $R^w$ (proof of Theorem \ref{thm:Tw}), satisfies $J(x\cdot y)=Jx\cdot Jy$ and $J[\mathrm{pt}]=1$, and commutes with $\partial$ and $\Delta$ by Stokes' theorem and Fubini's theorem in charts, as in Section \ref{ssec:currents}, with the weights $1/|G|$. The Maurer--Cartan equation holds in the quotient by $R^w$, not in $C^w_*$ (Proposition \ref{prop:obstruction}).
\end{proof}

\begin{remark}\label{rem:wide-G}
(1) Only the covering relation is needed for the identities of Corollary \ref{cor:exp}: the span of the covering relations is itself a differential graded BV ideal. Both relations are used in the proof of Theorem \ref{thm:Tw}. We do not know whether $\rho$ maps $R^w\cap C^{w,R}_*$ into $G^{\mathrm{mor}}_*$; Corollary \ref{cor:exp} and Theorem \ref{thm:connected-G} are proved independently.

(2) With Theorem \ref{thm:Tw}, the obstacle of Remark \ref{rem:Costello}(3) disappears: $\Delta$ and the bracket are defined on complexes computing the homology of the moduli spaces, and $V_\hbar$ is a Maurer--Cartan element there. The comparison with \cite{Cos05} is carried out in Section \ref{sec:Costello}. It does not need a quasi-isomorphism of differential graded Lie algebras: a morphism compatible with sewing, defined on a strict model of Costello's framework, together with the low degree homology given by Theorem \ref{thm:Tw}, suffices (Theorem \ref{thm:F}).
\end{remark}

\section{Comparison with Costello's model}\label{sec:Costello}
Throughout this section $0<L\le\Lst$. Costello \cite{Cos05} constructs a solution of the quantum master equation in a Batalin--Vilkovisky algebra of chains on the moduli spaces of Riemann surfaces with parametrized boundary, unique up to homotopy. His algebra is built from homotopy coinvariants of the circle actions which rotate the parametrizations, because on chains of the strict quotients the sewing operator is not defined: ``pull-back is not well defined on simplicial chains'' \cite[footnote 6, p.~5]{Cos05}. Our chains live on the strict quotients $\N_{\chi,n}$, where pullback along the sewing bundles is defined for geometric chains (Section \ref{ssec:Delta}). We construct a model of Costello's framework on which gluing is strictly associative, the moduli spaces $\Nt$ of hyperbolic surfaces with a marked point on each boundary component (Section \ref{ssec:Ntilde}), and a residue map $\Phi$ from Costello's algebra for $\Nt$ to the complex $\mathcal C^0/R^w$ of Section \ref{sec:wide-G} (Section \ref{ssec:Phi}). The results are Theorem \ref{thm:E}, which identifies the integration currents of all suitably normalized solutions in this model with the vertex currents, and Theorem \ref{thm:F}, which identifies the homotopy classes of Maurer--Cartan elements, Costello's class corresponding to the class of the hyperbolic vertices if the comparison map $h$ to his moduli spaces is continuous.

\subsection{Costello's framework}\label{ssec:Cos-framework}
We recall the parts of \cite[Sections 5--9]{Cos05} that are used, over $\mathbb Q$ and with the $\mathbb Z$-grading by chain degree, differentials having degree $-1$. Here $C_*$ denotes normalized singular chains; normalization is needed in \cite[footnote 4 and Lemma 6.0.2]{Cos05}.

\emph{Functors.} Let $\mathsf A$ be Costello's topological category of infinitely thin annuli \cite[Section 6]{Cos05}, with objects the integers $n\ge0$, $\mathsf A(1,1)=S^1$ and $\mathsf A(n,n)=S^1\wr S_n$. By \cite[Lemma 6.0.2]{Cos05}, $H_*(\mathsf A)\to C_*(\mathsf A)$ is a quasi-isomorphism, and a symmetric monoidal functor $F:H_*(\mathsf A)\to\mathrm{Comp}_{\mathbb Q}$ is given by complexes $F(n)$ with actions of $S_n$ and products $F(n)\otimes F(m)\to F(n+m)$, odd operators $D_i$ on $F(n)$, $1\le i\le n$, the images of the fundamental class of $\mathsf A(1,1)$, and maps $G_{ij}:F(n)\to F(n-2)$, $i<j$, the images of the class of a point of $\mathsf A(2,0)$, compatible with the symmetric groups and the products and subject to
\begin{equation}\label{eq:Cos-axioms}
D_i^2=[D_i,d]=[D_i,D_j]=0,\qquad[G_{ij},d]=G_{ij}\circ(D_i-D_j)=0 .
\end{equation}
For a symmetric monoidal functor $M:\mathsf A\to\mathrm{Top}$ Costello uses the functor
\begin{equation}\label{eq:Cprime}
C'_*(M)(m):=\bigoplus_{k\ge0}\ \bigoplus_{[m]=I_1\sqcup\dots\sqcup I_k}\big(C_*(M^{\mathrm{conn}}(I_1))\otimes\dots\otimes C_*(M^{\mathrm{conn}}(I_k))\big)_{S_k},
\end{equation}
the sum over decompositions of $[m]=\{1,\dots,m\}$ into possibly empty subsets, where $M^{\mathrm{conn}}(I)$ consists of the connected surfaces with boundary components labelled by $I$; $D_i$ acts on the factor containing $i$, and $G_{ij}$ glues within a factor, or, if $i$ and $j$ lie in different factors, after the shuffle map $\mathrm{EZ}:C_*(X)\otimes C_*(Y)\to C_*(X\times Y)$ \cite[Section 6]{Cos05}. His example is the coarse moduli space $M(m)$ of possibly disconnected Riemann surfaces with $m$ parametrized boundary components, where components without boundary are allowed \cite[Section 6 and footnote 5]{Cos05}, and the subfunctor $M^s$ of surfaces all of whose components have negative Euler characteristic.

\emph{Homotopy coinvariants and the BV algebra} \cite[Sections 5 and 7.2]{Cos05}. Put
\begin{equation*}
F_{hS}(n):=\big(F(n)\otimes t_1^{-1}\mathbb Q[t_1^{-1}]\otimes\dots\otimes t_n^{-1}\mathbb Q[t_n^{-1}]\big)_{S_n},
\end{equation*}
with differential $d+\sum_it_iD_i$, where $t_i^{-k}$ has degree $2k-2$, as in \cite[Section 5.1]{Cos05}, where $t^{-k}$ corresponds to a chain of dimension $2k-2$ in $\mathbb{CP}^\infty$; so $t_i^{-1}$ has degree $0$. The sum $\mathcal F(F):=\bigoplus_nF_{hS}(n)$ is a graded commutative algebra, the product being induced by the products of $F$, and
\begin{equation}\label{eq:Cos-Delta}
\Delta^{\Cos}:=\sum_{i<j}\Delta_{ij},\qquad\Delta_{ij}(a\otimes f):=G_{ij}(D_ia)\otimes\mathrm{Res}_{t_i}\mathrm{Res}_{t_j}f,
\end{equation}
where $\mathrm{Res}_{t_i}\mathrm{Res}_{t_j}f$ is the coefficient of $t_i^{-1}t_j^{-1}$ in $f$, the remaining variables being renumbered in increasing order, descends to $\mathcal F(F)$, has order two and satisfies $(\Delta^{\Cos})^2=[d,\Delta^{\Cos}]=0$ \cite[Section 7.2]{Cos05}. A natural transformation of such functors induces a morphism of these BV algebras \cite[Section 7.2]{Cos05}.

\emph{The connected differential graded Lie algebra} \cite[Sections 8 and 9]{Cos05}. Let $F_{g,n}\subseteq F_{hS}(n)$ be the part coming from connected surfaces of genus $g$, for the functors below the summands of \eqref{eq:Cprime} with $k=1$. Costello's algebra $\mathfrak g(F)$ consists of the series $S=\sum S_{g,n}\lambda^{2g-2+n}$ over $2g-2+n>0$, with $S_{g,n}\in F_{g,n}$, and $S\in\mathfrak g(F)^i$ if each $S_{g,n}$ has degree $6g-6+2n+1+i$. Its differential is $\hat d=d+\Delta^{\Cos}$, and its bracket is given by the odd Poisson bracket $\{f,g\}=\hat d(fg)-(-1)^{|f|}f\hat dg-(\hat df)g$; with the sign $[f,g]:=(-1)^{|f|}\{f,g\}$ of Corollary \ref{cor:exp}, $|f|$ the chain degree, it is a differential graded Lie algebra for the grading by $i$, whose parity is that of the shifted grading of Corollary \ref{cor:exp}. On connected elements, $\Delta^{\Cos}$ maps $F_{g-1,n+2}$ to $F_{g,n}$, and the bracket maps $F_{g_1,n_1}\otimes F_{g_2,n_2}$ to $F_{g_1+g_2,n_1+n_2-2}$, sewing one boundary component of each factor; both respect $w:=2g-2+n$, additively for the bracket, so the powers of $\lambda$ only record $w$. A Maurer--Cartan element is an $S\in\mathfrak g(F)^{-1}$, i.e. with $S_{g,n}$ of degree $6g-6+2n$, such that $\hat dS+\frac12\{S,S\}=0$, i.e. $(d+\Delta^{\Cos})e^S=0$; a homotopy is a solution of the master equation over $\mathbb Q[t,dt]$, $dt$ of degree $-1$, with the analogous conditions, and $\pi_0\mathrm{MC}$ denotes the set of homotopy classes \cite[Section 8.1 and Proposition 9.0.2]{Cos05}; we use the completed version of this notion defined after Lemma \ref{lem:GM}. Let $F^k\mathfrak g(F)$ consist of the $S$ with $S_{g,n}=0$ for $2g-2+n<k$. Then $\mathfrak g(F)$ is filtered pro-nilpotent in the sense of \cite[Section 8.1]{Cos05}: complete, $\mathfrak g=F^1\mathfrak g$ and $[F^k,F^l]\subseteq F^{k+l}$. On $\mathrm{Gr}^k\mathfrak g(F)=\prod_{2g-2+n=k}F_{g,n}$ the bracket vanishes and the differential is $d+\Delta^{\Cos}$, which raises the genus.

\begin{lemma}[{\cite[Lemma 8.1.1]{Cos05}}]\label{lem:GM}
Let $f:\mathfrak g\to\mathfrak g'$ be a filtration preserving map of filtered pro-nilpotent differential graded Lie algebras such that $\mathrm{Gr}\,f$ induces isomorphisms on $H_i$ for $i=0,-1,-2$. Then $f$ induces a bijection $\pi_0\mathrm{MC}(\mathfrak g)\to\pi_0\mathrm{MC}(\mathfrak g')$.
\end{lemma}
Costello sketches the proof and refers to Goldman--Millson type results in the literature for it \cite[proof of Lemma 8.1.1]{Cos05}. We do not use Lemma \ref{lem:GM}. In our situation the relevant groups vanish except one (Lemma \ref{lem:low-degrees}), and the direct order-by-order argument of Lemma \ref{lem:MC} below suffices.

\emph{Homotopies.} Let $\Omega:=\mathbb Q[t,dt]$, $t$ of degree $0$ and $dt$ of degree $-1$. For a filtered pro-nilpotent differential graded Lie algebra $\mathfrak g$ we use the completed tensor product
\begin{equation*}
\mathfrak g\hat\otimes\Omega:=\varprojlim_k\,(\mathfrak g/F^k\mathfrak g)\otimes\Omega,
\end{equation*}
with $d(y\otimes\omega)=dy\otimes\omega+(-1)^{|y|}y\otimes d\omega$ and $[y\otimes\omega,y'\otimes\omega']=(-1)^{|\omega||y'|}[y,y']\otimes\omega\omega'$. For $\mathfrak g=\mathfrak g(F)$ or $\mathfrak h^w$ below, $\mathfrak g/F^k\mathfrak g$ is a finite product of the $F_{g,n}$, resp.\ of the $C^w_*(\N_L(S_{g,n}))/R^w$, so an element of $\mathfrak g\hat\otimes\Omega$ is a series over $(g,n)$ whose components are polynomial in $t$ and $dt$, of unbounded degree in $t$. A homotopy is a Maurer--Cartan element of $\mathfrak g\hat\otimes\Omega$, and $\pi_0\mathrm{MC}(\mathfrak g)$ is the quotient of the set $\mathrm{MC}(\mathfrak g)$ of Maurer--Cartan elements by the equivalence relation generated by homotopy, the ends being the evaluations at $t=0,1$, $dt=0$. A homotopy in the sense of \cite[Section 8.1]{Cos05}, with coefficients in $\mathfrak g\otimes\Omega$, is one in this sense; so homotopic elements in his sense are homotopic in ours, and all statements on $\pi_0\mathrm{MC}$ below refer to our notion. Writing $X=x(t)+a(t)\,dt$, the Maurer--Cartan equation for $X$ says that each $x(t)$ is a Maurer--Cartan element and that
\begin{equation}\label{eq:homotopy-ODE}
\partial_tx(t)=d_{x(t)}a(t),\qquad d_xa:=da+[x,a].
\end{equation}
For a Maurer--Cartan element $x$, $dx=-\frac12[x,x]\in F^2\mathfrak g$, so $x$ defines a class $c(x)\in H_{-1}(\mathrm{Gr}^1\mathfrak g)$, $\mathrm{Gr}^1\mathfrak g=\mathfrak g/F^2\mathfrak g$.

\begin{lemma}\label{lem:MC}
Let $\mathfrak g$ be a filtered pro-nilpotent differential graded Lie algebra.
\begin{enumerate}
\item The class $c(x)$ depends only on the image of $x$ in $\pi_0\mathrm{MC}(\mathfrak g)$, and $c(f(x))=(\mathrm{Gr}^1f)_*c(x)$ for each filtration preserving morphism $f$.
\item If $H_{-2}(\mathrm{Gr}^k\mathfrak g)=0$ for all $k\ge2$, every class in $H_{-1}(\mathrm{Gr}^1\mathfrak g)$ is $c(x)$ for some $x\in\mathrm{MC}(\mathfrak g)$.
\item If $H_{-1}(\mathrm{Gr}^k\mathfrak g)=0$ for all $k\ge2$, any $x,x'\in\mathrm{MC}(\mathfrak g)$ with $c(x)=c(x')$ are joined by a single homotopy.
\end{enumerate}
In particular, if both conditions hold, $c$ induces a bijection $\pi_0\mathrm{MC}(\mathfrak g)\to H_{-1}(\mathrm{Gr}^1\mathfrak g)$, and a filtration preserving morphism between two such algebras which induces an isomorphism on $H_{-1}(\mathrm{Gr}^1)$ induces a bijection on $\pi_0\mathrm{MC}$.
\end{lemma}
\begin{proof}
(1) Modulo $F^2\mathfrak g$, \eqref{eq:homotopy-ODE} reads $\partial_tx(t)\equiv da(t)$, the bracket lying in $F^2$; in $(\mathfrak g/F^2\mathfrak g)\otimes\Omega$ everything is polynomial in $t$, so $x(1)-x(0)\equiv d\int_0^1a(t)\,dt$ modulo $F^2$. Naturality is clear.

(2) Represent the class by $x_1\in\mathfrak g^{-1}$ with $dx_1\in F^2$. Suppose $x\in\mathfrak g^{-1}$ satisfies $R(x):=dx+\frac12[x,x]\in F^{k+1}$, $k\ge1$. Then $dR(x)=[dx,x]=[R(x),x]$, since $[[x,x],x]=0$ for $x$ odd by the Jacobi identity, so $dR(x)\in F^{k+2}$ and $R(x)$ is a cycle of degree $-2$ in $\mathrm{Gr}^{k+1}\mathfrak g$. By assumption $R(x)\equiv-dy$ modulo $F^{k+2}$ for some $y\in F^{k+1}\mathfrak g^{-1}$, and $R(x+y)=R(x)+dy+[x,y]+\frac12[y,y]\in F^{k+2}$. Starting with $x_1$, for which $R(x_1)\in F^2$, we obtain $x_1,x_2,\dots$ with $x_{k+1}-x_k\in F^{k+1}$ and $R(x_k)\in F^{k+1}$; by completeness $x:=\lim x_k$ exists, $R(x)=0$, and $c(x)=[x_1]$.

(3) Let $\mathfrak g^+:=\mathfrak g\oplus\mathbb Q\delta$, $\delta$ of degree $-1$, with $[\delta,y]:=dy$ and $[\delta,\delta]:=0$; this is a graded Lie algebra since $d$ is a derivation with $d^2=0$, and $x\in\mathfrak g^{-1}$ is a Maurer--Cartan element if and only if $[x+\delta,x+\delta]=0$. For $A\in\mathfrak g^0$, $\mathrm{ad}_A$ maps $F^k\mathfrak g$ to $F^{k+1}\mathfrak g$ and $\delta$ to $-dA\in F^1\mathfrak g$, so $(\mathrm{ad}_A)^m$ maps $\mathfrak g^+$ to $F^m\mathfrak g$ for $m\ge1$, and $e^{\mathrm{ad}_A}$ converges and is an automorphism of $\mathfrak g^+$ fixing $\delta$ modulo $\mathfrak g$; the gauge action $e^A\cdot x:=e^{\mathrm{ad}_A}(x+\delta)-\delta=e^{\mathrm{ad}_A}x-\frac{e^{\mathrm{ad}_A}-1}{\mathrm{ad}_A}(dA)$ preserves Maurer--Cartan elements, satisfies $e^A\cdot x\equiv x+[A,x]-dA$ modulo terms quadratic in $A$, and $e^{A}\cdot(e^{B}\cdot x)=e^{\mathrm{BCH}(A,B)}\cdot x$ by the Baker--Campbell--Hausdorff formula in the pro-nilpotent Lie algebra $\mathfrak g^0$. The path $x(t):=e^{tA}\cdot x'$ lies in $\mathfrak g\hat\otimes\mathbb Q[t]$, the coefficient of $t^m$ lying in $F^m\mathfrak g$, and $\partial_t(x(t)+\delta)=[A,x(t)+\delta]$, i.e.\ $\partial_tx(t)=[A,x(t)]-dA=-d_{x(t)}A$; so $x(t)-A\,dt$ is a homotopy from $x'$ to $e^A\cdot x'$ by \eqref{eq:homotopy-ODE}. It therefore suffices to find $A\in\mathfrak g^0$ with $e^A\cdot x'=x$.

We construct $A_1,A_2,\dots\in\mathfrak g^0$ with $x-e^{A_k}\cdot x'\in F^{k+1}$ and $A_{k+1}-A_k\in F^{k+1}$. Let $y\in\mathrm{MC}(\mathfrak g)$ with $x-y\in F^k$, $k\ge1$. Then $d(x-y)=-\frac12[x-y,x+y]\in F^{k+1}$, so $x-y$ is a cycle of degree $-1$ in $\mathrm{Gr}^k\mathfrak g$; its class vanishes, for $k=1$ because $c(x)=c(y)$ and for $k\ge2$ by assumption. So $x-y\equiv-da$ modulo $F^{k+1}$ for some $a\in F^k\mathfrak g^0$, and then $e^a\cdot y\equiv y-da\equiv x$ modulo $F^{k+1}$, since $[a,y]$ and the higher terms lie in $F^{k+1}$. Applying this to $y=x'$ and $k=1$ we put $A_1:=a$; applying it to $y=e^{A_k}\cdot x'$, with $k+1$ in place of $k$, we put $A_{k+1}:=\mathrm{BCH}(a,A_k)\equiv A_k+a$ modulo $F^{k+2}$, so that $e^{A_{k+1}}\cdot x'=e^a\cdot y\equiv x$ modulo $F^{k+2}$. The limit $A$ exists by completeness, and $e^A\cdot x'=x$ by continuity.
\end{proof}

\emph{Costello's solution} \cite[Theorem 1 and Proposition 9.0.2]{Cos05}. For $F=C'_*(M^s)$ there is a Maurer--Cartan element $S$ of $\mathfrak g(F)$ such that $S_{0,3}$ is a $0$-chain of degree, i.e. total weight, $1/3!$ in the moduli space of surfaces with three unparametrized unordered boundary components, and it is unique up to homotopy. The proof uses the bound $H_i(M_{g,n}/S_n;\mathbb Q)=0$ for $i\ge6g-7+2n$, $(g,n)\neq(0,3)$ \cite[(9.0.1)]{Cos05}, to show that $H_i(\mathrm{Gr}^k\mathfrak g)=0$ for $i\ge0$ and for $i=-2$, and for $i=-1$ if $k\ge2$, while $H_{-1}(\mathrm{Gr}^1\mathfrak g)=\mathbb Q$; then Lemma \ref{lem:GM}, applied to $\mathfrak g\to\mathfrak g/F^2\mathfrak g$, an abelian quotient, gives $\pi_0\mathrm{MC}(\mathfrak g)\cong H_{-1}(\mathrm{Gr}^1\mathfrak g)=\mathbb Q$; Lemma \ref{lem:MC} gives this directly from the same vanishing, for the completed notion of homotopy defined after Lemma \ref{lem:GM}. The class of a Maurer--Cartan element in $H_{-1}(\mathrm{Gr}^1\mathfrak g)\cong\mathbb Q$ is the degree of $S_{0,3}$, i.e. the sum of the coefficients of the $0$-chain multiplying $t_1^{-1}t_2^{-1}t_3^{-1}$: in degree $0$ only these terms occur in the component $(0,3)$, and this degree is a well defined linear form on the cycles of $\mathrm{Gr}^1$ which induces the isomorphism (Lemma \ref{lem:low-degrees}). In the sense of Section \ref{sec:intro}, ``Costello's class'' is the class of $S$ for this completed notion of homotopy. The identification of $H_*(F_{g,n})$ with $H_*(M_{g,n}/S_n;\mathbb Q)$ used here rests on \cite[Proposition 5.1.1]{Cos05}, whose proof is sketched, and on the rational homotopy equivalence of $M_g(n)/(S^1)^n$ with $M_{g,n}$, stated in the proof of \cite[Proposition 9.0.2]{Cos05}. Costello remarks that the result holds for any sequence of spaces with the same gluing structure \cite[Remark 2 after Theorem 1]{Cos05}.

\emph{The geometric side.} Let $\mathfrak h^w\subseteq\mathfrak g^w[[\hbar]]$ (Corollary \ref{cor:exp}) be the set of series $\sum_{g,n}\hbar^gx_{g,n}$ with $x_{g,n}\in C^w_*(\N_L(S_{g,n}))/R^w$, i.e. those in which the power of $\hbar$ on $\N_L(S_{g,n})$ is the genus. It is closed under $\partial+\hbar\Delta_0$ and the bracket, since $\Delta_0$ raises the genus by one and the bracket adds genera, and
\begin{equation}\label{eq:lambda-hbar}
\textstyle\sum_{g,n}x_{g,n}\lambda^{2g-2+n}\ \longmapsto\ \sum_{g,n}\hbar^gx_{g,n}
\end{equation}
identifies $\mathfrak h^w$, with $\partial+\hbar\Delta_0$ and $\{\,,\}$, with $\prod_{g,n}C^w_*(\N_L(S_{g,n}))/R^w$ in Costello's bookkeeping, with $\partial+\Delta_0$ and $\{\,,\}$; both formal variables only record the type. We grade $\mathfrak h^w$ as Costello does: $x\in(\mathfrak h^w)^i$ if each $x_{g,n}$ has dimension $6g-6+2n+1+i$, i.e. codimension $-1-i$, and we filter it by $2g-2+n$. Then $\mathfrak h^w$ is a filtered pro-nilpotent differential graded Lie algebra, its Maurer--Cartan elements are the series of elements of codimension $0$ with $(\partial+\Delta_0)x+\frac12\{x,x\}=0$, among them $V_\hbar$ (Corollary \ref{cor:exp}), and homotopies and $\pi_0\mathrm{MC}$ are defined as above.

\begin{lemma}\label{lem:bound}
Let $2g-2+n>0$ and $(g,n)\ne(0,3)$. Then the rational homology $H_j(|\N_L(S_{g,n})|;\mathbb Q)$ vanishes for $j\ge6g-7+2n$.
\end{lemma}
\begin{proof}
$\Mod(S_{g,n})$ acts properly, with finite stabilisers, on $\T_L(S_{g,n})$, which is contractible by Section \ref{ssec:FN}, so $H_*(|\N_L(S_{g,n})|;\mathbb Q)\cong H_*(\Mod(S_{g,n});\mathbb Q)$ \cite[Chapter VII]{Bro}. The group $\Mod(S_{g,n})$ is the mapping class group of the surface of genus $g$ with $n$ punctures, permuting the punctures; it contains the pure mapping class group with finite index, so by Harer \cite{Ha} its rational homology vanishes above $4g-5$ for $n=0$, above $4g-4+n$ for $g,n\ge1$, and above $n-3$ for $g=0$. These bounds are smaller than $6g-7+2n$ except for $(g,n)=(1,1)$ and $(0,4)$, where they equal $1=6g-7+2n$. For $(1,1)$, $\Mod(S_{1,1})\cong SL_2(\mathbb Z)$ has finite abelianization, so $H_1=0$. For $(0,4)$, $H_1(|\N_L(S_{0,4})|;\mathbb Q)$ is the space of $S_4$-coinvariants of $H_1(M_{0,4};\mathbb Q)$, which vanishes \cite[Section 2.2.2]{CZ}, \cite[proof of Proposition 9.0.2]{Cos05}.
\end{proof}

\subsection{Hyperbolic surfaces with marked boundary points}\label{ssec:Ntilde}
Let $I$ be a finite set with $n=|I|\ge1$ elements, and $S$ a connected type of genus $g$ with boundary components $b_i$ labelled by $i\in I$. Let $\hat\T_L(S)$ be the set of pairs $(x,p)$ with $x\in\T_L(S)$ and $p=(p_i)_{i\in I}$, $p_i\in(b_i)_x$; via the seam feet of Section \ref{ssec:sewing} it is the manifold $\T_L(S)\times(\mathbb R/L\mathbb Z)^I$. The group $\Mod^I(S)\subseteq\Mod(S)$ of mapping classes preserving every boundary component acts on it by bundle maps, as in Section \ref{ssec:sewing}. Put
\begin{equation*}
\Nt_g(I):=\Mod^I(S)\backslash\hat\T_L(S),\qquad\Nt_g(n):=\Nt_g(\{1,\dots,n\}).
\end{equation*}
Its points are the isometry classes of hyperbolic surfaces of genus $g$ with geodesic boundary of length $L$ and boundary components labelled by $I$, with a marked point on each boundary component, isometries being required to preserve labels and marked points.

\begin{lemma}\label{lem:Ntilde}
\begin{enumerate}
\item $\Mod^I(S)$ acts freely and properly on $\hat\T_L(S)$, so $\Nt_g(I)$ is a connected manifold of dimension $6g-6+3n$.
\item Let the torus $T^n:=(\mathbb R/\mathbb Z)^n$ act on $\Nt_g(n)$ by $\theta\cdot p_i:=\rot^{b_i}_{-L\theta_i}(p_i)$, and $S_n$ by relabelling. The action of $G_n:=T^n\rtimes S_n$ has finite stabilisers, the stabiliser of the class of a surface $X$ being isomorphic to the group of orientation preserving isometries of $X$, and $G_n\ltimes\Nt_g(n)$ is Morita equivalent to $\N_L(S_{g,n})$. Let $q:\Nt_g(n)\to\N_L(S_{g,n})$ be the resulting morphism.
\item $q\circ\gamma\cong q$ for $\gamma\in G_n$, and $q\circ\rot_i\cong q\circ\mathrm{pr}_2$ as morphisms $(\mathbb R/\mathbb Z)\times\Nt_g(n)\to\N_L(S_{g,n})$, where $\rot_i$ is the action of the $i$-th circle of $T^n$.
\end{enumerate}
\end{lemma}
\begin{proof}
(1) If $\phi\in\Mod^I(S)$ fixes $(x,p)$, the isometry $\psi$ of $X_x$ induced by $\phi$ fixes $p_i$ and maps the oriented geodesic $b_i$ to itself, preserving its orientation and the side of $X_x$; so $d\psi$ is the identity at $p_i$, and $\psi=\mathrm{id}$, $X_x$ being connected. Then $\phi$ is isotopic to the identity, by the definition of $\T_L(S)$. Properness follows from the properness of the action on $\T_L(S)$ and the compactness of the fibres of $\hat\T_L(S)\to\T_L(S)$, and $\hat\T_L(S)$ is connected.
(2) Let $\hat\T$ be the set of triples $(x,\beta,p)$ with $x\in\T_L(S_{g,n})$, $\beta$ a bijection from $\{1,\dots,n\}$ to the set of boundary components of $S_{g,n}$, and $p_i$ a point of $(\beta(i))_x$. The groups $\Mod(S_{g,n})$ and $G_n$ act on $\hat\T$ with commuting actions. By (1), and since $\Mod(S_{g,n})$ permutes the labellings transitively with stabiliser $\Mod^{\{1,\dots,n\}}(S_{g,n})$, the first action is free and proper with quotient $\Nt_g(n)$; the second is free and proper with quotient $\T_L(S_{g,n})$, the fibre over $x$ being a $G_n$-torsor. Hence $\Mod(S_{g,n})\ltimes\T_L(S_{g,n})\simeq(\Mod(S_{g,n})\times G_n)\ltimes\hat\T\simeq G_n\ltimes\Nt_g(n)$. The stabiliser in $G_n$ of the class of $(x,\beta,p)$ is the image of the stabiliser of $x$ in $\Mod(S_{g,n})$, the isometry group of $X_x$, which maps injectively since $\Mod(S_{g,n})$ acts freely on $\hat\T$.
(3) In $G_n\ltimes\Nt_g(n)$ the arrows $(\gamma,y)$, resp.\ $(\theta e_i,y)$, which depend smoothly on $y$, resp.\ $(\theta,y)$, are the required $2$-isomorphisms.
\end{proof}

We use the following description of morphisms from simply connected sources, which avoids monodromy questions.

\begin{lemma}\label{lem:sc}
Let $Q$ be a connected and simply connected manifold with corners and $S$ a type.
\begin{enumerate}
\item Every morphism $Q\to\N_L(S)=\Mod(S)\ltimes\T_L(S)$ is $2$-isomorphic to the morphism induced by a smooth map $\tilde f:Q\to\T_L(S)$, and the morphisms induced by $\tilde f$ and $\tilde f'$ are $2$-isomorphic if and only if $\tilde f'=\phi\cdot\tilde f$ for some $\phi\in\Mod(S)$.
\item If $(X_s)_{s\in Q}$ is a smooth family of hyperbolic surfaces of type $S$ and $h_s,h'_s:S\to X_s$ are smooth families of markings, then $s\mapsto[X_s,h_s]$ and $s\mapsto[X_s,h'_s]$ induce $2$-isomorphic morphisms $Q\to\N_L(S)$.
\end{enumerate}
\end{lemma}
\begin{proof}
A morphism from $Q$ to the action groupoid of the discrete group $\Mod(S)$ is a principal $\Mod(S)$-bundle over $Q$ with an equivariant map to $\T_L(S)$, and $2$-isomorphisms are isomorphisms of such pairs. Since $Q$ is simply connected the bundle is trivial; a section gives $\tilde f$, and two sections differ by a locally constant, hence constant, element of $\Mod(S)$. This is (1). In (2), $h_s^{-1}\circ h'_s$ is a continuous family of orientation preserving diffeomorphisms of $S$, so its isotopy class $\phi$ does not depend on $s$, and $[X_s,h'_s]=\phi^{-1}\cdot[X_s,h_s]$.
\end{proof}

\emph{Gluing.} For two labels $i\ne j$ of $y\in\Nt_g(I)$ let $\varphi_0(y):b_i\to b_j$ be the orientation reversing isometry with $\varphi_0(p_i)=p_j$, and let $G_{ij}(y)$ be the surface obtained by gluing $b_i$ to $b_j$ by $\varphi_0(y)$, with the remaining labels and marked points; likewise for $y=(y_1,y_2)\in\Nt_{g_1}(I_1)\times\Nt_{g_2}(I_2)$ with $i\in I_1$, $j\in I_2$. If the glued surface has boundary, this defines maps
\begin{gather*}
G_{ij}:\Nt_g(I)\to\Nt_{g+1}(I\setminus\{i,j\}),\\
G_{ij}:\Nt_{g_1}(I_1)\times\Nt_{g_2}(I_2)\to\Nt_{g_1+g_2}((I_1\sqcup I_2)\setminus\{i,j\}),
\end{gather*}
which are real analytic: by Lemma \ref{lem:cut}, a lift to $\hat\T_L$ is given locally by $\cut^{-1}$ of the gluing datum $\varphi_0$, the remaining marked points lying on the same boundary circles. If the glued surface is closed we put $G_{ij}:=\glue\circ s_{ij}$, a morphism to $\N_L(S_{g+1,0})$, resp.\ $\N_L(S_{g_1+g_2,0})$, where $s_{ij}$ is the section $t=0$ of the bundle of Lemma \ref{lem:pullback-E}(1) below. In both cases we write $\bar G_{ij}$ for $G_{ij}$ followed by $q$ if the glued surface has boundary, and for $G_{ij}$ otherwise.

\begin{lemma}\label{lem:pullback-E}
Let $\tilde Y$ be $\Nt_g(n)$ with $n\ge2$, or $\Nt_{g_1}(I_1)\times\Nt_{g_2}(I_2)$ with $I_1,I_2\ne\emptyset$, and let $q_{\tilde Y}$ be $q$, resp.\ $\sqcup\circ(q\times q)$, a morphism to $\N:=\N_{\chi,n}$. Let $\Pi$ be the set of unordered pairs of labels, resp.\ of pairs $\{i,j\}$ with $i\in I_1$, $j\in I_2$, and $\Pp^\Pi\subseteq\Pp$, $\E^\Pi\subseteq\E$ the corresponding open and closed parts of the pullbacks of $\Pp$, $\E$ along $q_{\tilde Y}$.
\begin{enumerate}
\item $\Pp^\Pi=\tilde Y\times\Pi$ and $\E^\Pi=\bigsqcup_{\{i,j\}\in\Pi}(\mathbb R/L\mathbb Z)\times\tilde Y$, where $t\in\mathbb R/L\mathbb Z$ corresponds to the gluing isometry $t\cdot\varphi_0(y)$ of \eqref{eq:rotation}; the orientations of Definition \ref{def:pullback-wide} are the orientation of $t$ and the positive orientation of $\Pi$.
\item Let $c_{ij}:(\mathbb R/L\mathbb Z)\times\tilde Y\to\N'$ be $\glue$ composed with the projection to $\E$ on the piece $\{i,j\}$, and $\nu(\theta):=-L\theta$. For every connected and simply connected manifold with corners $Q$ and every smooth map $\tau:Q\to(\mathbb R/\mathbb Z)\times\tilde Y$, the morphisms $c_{ij}\circ(\nu\times\mathrm{id})\circ\tau$ and $\bar G_{ij}\circ\rot_i\circ\tau$ are $2$-isomorphic, and so are $c_{ij}\circ(\nu\times\mathrm{id})\circ\tau$ and $\bar G_{ij}\circ\rot_j\circ\tau$.
\end{enumerate}
\end{lemma}
\begin{proof}
(1) By Lemma \ref{lem:Ntilde}(2) the pullbacks along $q$ of $\Pp$ and $\E$ are the $G_n$-equivariant bundles over $\Nt_g(n)$ whose fibres over $y$ are the set of unordered pairs of boundary components of the underlying surface, identified with pairs of labels, and the circles of gluing isometries $\Phi_{b_i,b_j}$ of Section \ref{ssec:sewing}; the latter are trivialized by $t\mapsto t\cdot\varphi_0(y)$, which is smooth in $y$ and compatible with the orientation of \eqref{eq:rotation}. For the product, $\sqcup$ is induced by the inclusion $\Mod(S_1)\times\Mod(S_2)\subseteq\Mod(S_1\sqcup S_2)$ and the identity of $\T_L$, so the same description applies to the pairs meeting both components (proof of Proposition \ref{prop:BV}).
(2) Rotating $p_i$ by $\theta$, i.e. replacing it by $\rot^{b_i}_{-L\theta}(p_i)$, replaces $\varphi_0$ by the isometry $\varphi$ with $\varphi(\rot^{b_i}_{-L\theta}p_i)=p_j$; since $\varphi$ reverses orientations, $\varphi\circ\rot^{b_i}_{s}=\rot^{b_j}_{-s}\circ\varphi$, so $\varphi(p_i)=\rot^{b_j}_{-L\theta}(p_j)$, i.e. $\varphi=(-L\theta)\cdot\varphi_0$. Rotating $p_j$ by $\theta$ gives the same $\varphi$. So both morphisms have as underlying family the surfaces obtained by gluing along $(-L\theta)\cdot\varphi_0$. On the simply connected $Q$ both are induced by smooth families of markings of this family of glued surfaces: for $c_{ij}$ by the Morita inverse of $\cut$ (Lemma \ref{lem:cut}), which on simply connected sources is a lift to $H_\gamma$, and for $\bar G_{ij}$ by a lift of $G_{ij}\circ\rot_i\circ\tau$ to $\hat\T_L$, resp.\ by $\glue\circ s_{ij}$. By Lemma \ref{lem:sc}(2) they are $2$-isomorphic.
\end{proof}

\emph{The functor.} For an orbifold $N$ in the wide sense let $C^{\mathrm{sm}}_*(N)$ be the complex of $2$-isomorphism classes of morphisms $\Dt^k\to N$ of the proof of Theorem \ref{thm:Tw}, modulo the degenerate ones $\hat\sigma\circ s_j$; it is the normalized complex of a simplicial vector space, quasi-isomorphic to the unnormalized one. Define $F^{\Nt}(m)$ by \eqref{eq:Cprime}, with normalized smooth singular chains $C_*(\Nt^{\mathrm{conn}}(I)):=\bigoplus_gC_*(\Nt_g(I))$ for $I\neq\emptyset$ and with $C_*(\Nt^{\mathrm{conn}}(\emptyset)):=\bigoplus_{g\ge2}C^{\mathrm{sm}}_*(\N_L(S_{g,0}))$ for the closed components, the sums over $g$ with $2g-2+|I|>0$. Here $S_m$ acts by relabelling, the products are juxtaposition, $D_i\sigma:=(\rot_i)_*\mathrm{EZ}([S^1]\otimes\sigma)$ on the factor containing $i$, with the Koszul sign for the factors before it, where $[S^1]:\Dt^1\to\mathbb R/\mathbb Z$, $s\mapsto s$, and $G_{ij}$ is induced by the gluing maps, after the shuffle map if $i$ and $j$ lie in different factors.

\begin{lemma}\label{lem:functor}
$F^{\Nt}$ is a symmetric monoidal functor $H_*(\mathsf A)\to\mathrm{Comp}_{\mathbb Q}$ in the sense of \cite[Section 6]{Cos05}, i.e. it satisfies \eqref{eq:Cos-axioms} and the compatibilities with symmetric groups and products. Its connected parts are $F^{\Nt}_{g,n}=\big(C_*(\Nt_g(n))\otimes t_1^{-1}\mathbb Q[t_1^{-1}]\otimes\dots\otimes t_n^{-1}\mathbb Q[t_n^{-1}]\big)_{S_n}$ for $n\ge1$ and $F^{\Nt}_{g,0}=C^{\mathrm{sm}}_*(\N_L(S_{g,0}))$.
\end{lemma}
\begin{proof}
As for $C'_*(M)$ \cite[Section 6]{Cos05}, the axioms reduce to identities between maps and to properties of the shuffle map. Gluing is strictly associative and commutes strictly with the rotations of other boundary components, with relabellings and with disjoint unions, since points of $\Nt$ are isometry classes and the surfaces glued in either order are canonically isometric, compatibly with labels and marked points; for closed targets the corresponding identities hold up to $2$-isomorphism on simplices by Lemma \ref{lem:sc}(2), which suffices, $C^{\mathrm{sm}}_*$ being spanned by $2$-isomorphism classes. The shuffle map is natural, associative and graded commutative. $D_i^2=0$ holds on normalized chains \cite[Lemma 6.0.2]{Cos05}, and $[D_i,D_j]=0$ follows from the graded commutativity of the shuffle map, the actions commuting. Finally $G_{ij}\circ\rot_i=G_{ij}\circ\rot_j$ as maps $(\mathbb R/\mathbb Z)\times\Nt\to\Nt$, by the computation in the proof of Lemma \ref{lem:pullback-E}(2), and for closed targets up to $2$-isomorphism on simplices; hence $G_{ij}D_i=G_{ij}D_j$.
\end{proof}

\emph{The comparison with $M^s$.} For $y\in\Nt_g(I)$ let $h(y)$ be the Riemann surface underlying $X_y$, with the boundary component $b_i$ parametrized by $\theta\mapsto\rot^{b_i}_{-L\theta}(p_i)$, $\theta\in\mathbb R/\mathbb Z$, i.e. starting at $p_i$ and running against the boundary orientation, which is Costello's convention \cite[Section 6]{Cos05}; this parametrization is real analytic, $b_i$ being a geodesic. For closed surfaces $h$ is the uniformization homeomorphism $|\N_L(S_{g,0})|\to M_g(0)$. We use the convention that $\theta\in\mathsf A(1,1)=\mathbb R/\mathbb Z$ acts on parametrizations by $\phi_i\mapsto\phi_i\circ(\theta+\cdot\,)$ and that $D$ is the fundamental class of $\mathbb R/\mathbb Z$ oriented by $\theta$; Costello does not fix this orientation, and the other choice is discussed in Remark \ref{rem:signs}. Then $h$ is equivariant, and it commutes with gluing: gluing $b_i$ to $b_j$ by $\varphi_0$ is the conformal gluing along the parametrizations which identifies their starting points and reverses the direction, i.e. the gluing by the annulus of \cite[proof of Lemma 6.0.2]{Cos05} whose two parametrizations start at the same point. Costello does not specify the topology of $M(m)$ \cite[Section 6]{Cos05}; we assume that it is such that $h$ is continuous, as for the usual topologies, in which conformal structures and analytic boundary parametrizations depend continuously on the hyperbolic metric. This hypothesis enters only the statements about $M^s$, i.e.\ the second half of Theorem \ref{thm:E}(2) and Theorem \ref{thm:F}(2). Then $\sigma\mapsto h\circ\sigma$, and $\hat\sigma\mapsto h\circ|\hat\sigma|$ on closed components, is a natural transformation $F^{\Nt}\to C'_*(M^s)$, and it induces a morphism of BV algebras and of filtered differential graded Lie algebras
\begin{equation*}
h_*:\mathfrak g(F^{\Nt})\to\mathfrak g(C'_*(M^s))
\end{equation*}
\cite[Section 7.2]{Cos05}, which maps $y\otimes t_1^{-1}t_2^{-1}t_3^{-1}$ to $h(y)\otimes t_1^{-1}t_2^{-1}t_3^{-1}$ for a point $y$ of $\Nt_0(3)$, and hence preserves the degree of the component $(0,3)$. Only the statements about $M^s$ below use $h$.

\begin{proposition}\label{prop:equivariant}
For $2g-2+n>0$ there is an isomorphism $H_*(F^{\Nt}_{g,n})\cong H_*(|\N_L(S_{g,n})|;\mathbb Q)$, without shift of degrees. Moreover $H_0(F^{\Nt}_{g,n})=\mathbb Q$, spanned by the class of $y\otimes t_1^{-1}\cdots t_n^{-1}$ for a point $y$.
\end{proposition}
This is the input from equivariant topology on which Theorem \ref{thm:E}(2) and Theorem \ref{thm:F} depend. For $n=1$ its main step, (a) below, is the theorem of Jones and Hood--Jones \cite{J87,HJ87}; for a torus we give the argument.
\begin{proof}
For $n=0$ this is the statement on $C^{\mathrm{sm}}_*$ in the proof of Theorem \ref{thm:Tw}. Let $n\ge1$ and $X:=\Nt_g(n)$, a connected manifold with a smooth action of the torus $T=T^n$ with finite stabilisers (Lemma \ref{lem:Ntilde}).

(a) \emph{The complex $K_X:=C_*(X)\otimes t_1^{-1}\mathbb Q[t_1^{-1}]\otimes\dots\otimes t_n^{-1}\mathbb Q[t_n^{-1}]$, with differential $d+\sum_it_iD_i$ and $t_i^{-k}$ in degree $2k-2$, computes $H_*(ET\times_TX;\mathbb Q)$, naturally in $X$.} For $n=1$ this is the theorem of Jones and Hood--Jones, which identifies the equivariant homology of an $S^1$-space with the cyclic homology of its normalized singular chains with the circle operator \cite{J87,HJ87}, see \cite[Chapter 7]{Lo}. In general we argue as follows, first with continuous normalized chains; the argument is the homological counterpart of the Koszul duality of \cite{GKM}.

(a1) \emph{An exterior algebra.} Normalized chains $C_*(T)$ form a differential graded algebra under the Pontryagin product $x\cdot y:=\mu_*\mathrm{EZ}(x\otimes y)$, $\mu$ the group law, associative by the associativity and naturality of the shuffle map, and $C_*(X)$ is a left differential graded module over it by $x\cdot\sigma:=a_*\mathrm{EZ}(x\otimes\sigma)$, $a$ the action. Let $\Lambda:=\Lambda_{\mathbb Q}[D_1,\dots,D_n]$ be the exterior algebra on generators of degree $1$, with zero differential, and $\iota:\Lambda\to C_*(T)$ the algebra map with $\iota(D_i):=[S^1]_i:=(e_i)_*[S^1]$, where $e_i:\mathbb R/\mathbb Z\to T$ is the $i$-th circle. It is well defined and a map of differential graded algebras: $\partial[S^1]_i=0$; $[S^1]_i\cdot[S^1]_j=-[S^1]_j\cdot[S^1]_i$, since the shuffle map is graded commutative and $\mu$ is commutative; and $[S^1]_i\cdot[S^1]_i=0$, since the two shuffle simplices of $\Dt^1\times\Dt^1$, composed with $(s,s')\mapsto e_i(s+s')$, are the same singular simplex, the affine map with the values $0,1,2$ at the vertices followed by $e_i$, with opposite signs. It is a quasi-isomorphism, since by the K\"unneth theorem $H_*(T;\mathbb Q)$ with the Pontryagin product is the exterior algebra on the classes of the circles. Through $\iota$, $D_i$ acts on $C_*(X)$ by $[S^1]_i\cdot\sigma=(\rot_i)_*\mathrm{EZ}([S^1]\otimes\sigma)$, the operator $D_i$ of Section \ref{ssec:Ntilde}, $\rot_i=a\circ(e_i\times\mathrm{id})$.

(a2) \emph{The bar construction.} Write $M:=C_*(X)$. Since $T$ is commutative and the shuffle map is graded commutative, $C_*(T)$ is graded commutative, and we regard the left module $M$ also as a right module by $\sigma\cdot x:=(-1)^{|\sigma||x|}x\cdot\sigma$, as in (a3) below; correspondingly the face maps of $B_\bullet(X,T,\mathrm{pt})$ use the action $X\times T\to X$, $(x,\theta)\mapsto\theta\cdot x$. The two-sided bar construction $B(M,C_*(T),\mathbb Q)$ computes $M\otimes^{\mathbb L}_{C_*(T)}\mathbb Q$. Bar constructions are normalized, i.e.\ built from the suspended augmentation ideals, except in the comparison with $X\times T^\bullet$ below, which uses the unnormalized bar construction, the total complex of the simplicial chain complex $M\otimes C_*(T)^{\otimes\bullet}$; the normalized one is its quotient by the images of the degeneracies, and the quotient map is a quasi-isomorphism by the normalization theorem for simplicial chain complexes, naturally in the data. The simplicial space $B_\bullet(X,T,\mathrm{pt})=X\times T^\bullet$ is proper (good), its degeneracies inserting the unit of the Lie group $T$ and hence being closed cofibrations, so the homology of its realization is computed by the total complex of the simplicial chain complex of its singular chains. The shuffle maps $C_*(X)\otimes C_*(T)^{\otimes k}\to C_*(X\times T^k)$ are quasi-isomorphisms compatible with the simplicial operators of $B_\bullet(X,T,\mathrm{pt})$, so they identify $B(M,C_*(T),\mathbb Q)$ up to quasi-isomorphism with the total complex of the simplicial chain complex $C_*(X\times T^\bullet)$, which computes the homology of the realization $|B_\bullet(X,T,\mathrm{pt})|$, which is $X\times_T|B_\bullet(T,T,\mathrm{pt})|$, a model of $ET\times_TX$; for $X$ a point this is the approach of Rothenberg and Steenrod to $H_*(BT)$ \cite{RS65}. Since $\iota$ is a quasi-isomorphism of augmented algebras concentrated in non-negative degrees, the induced map $B(M,\Lambda,\mathbb Q)\to B(M,C_*(T),\mathbb Q)$ is a quasi-isomorphism: the suspended augmentation ideals sit in positive degrees, so the filtration by bar degree is finite in each total degree, and on its first page the map is an isomorphism by the K\"unneth theorem over $\mathbb Q$.

(a3) \emph{The Koszul resolution.} Let $\Gamma$ be the graded vector space with basis $u^{(c)}$, $c\in\mathbb Z_{\ge0}^n$, of degree $2|c|$, and $K:=\Lambda\otimes\Gamma$ with the differential $d(\lambda\otimes u^{(c)}):=(-1)^{|\lambda|}\sum_i\lambda D_i\otimes u^{(c-e_i)}$, where $u^{(c-e_i)}:=0$ if $c_i=0$. This differential is compatible with the left $\Lambda$-module structure, squares to zero because the $D_i$ anticommute, and $K$ is the tensor product over $i$ of the complexes $\Lambda[D_i]\otimes\Gamma[u_i]$, each of which has homology $\mathbb Q$ in degree $0$; so $K$ is a semifree resolution of $\mathbb Q$ over $\Lambda$, as is $B(\Lambda,\Lambda,\mathbb Q)$, and $M\otimes_\Lambda K\simeq B(M,\Lambda,\mathbb Q)$. Regarding $C_*(X)$ as a right $\Lambda$-module by $\sigma\cdot\lambda:=(-1)^{|\sigma||\lambda|}\lambda\cdot\sigma$, we have $C_*(X)\otimes_\Lambda K=C_*(X)\otimes\Gamma$ with differential $d(\sigma\otimes u^{(c)})=d\sigma\otimes u^{(c)}+\sum_iD_i\sigma\otimes u^{(c-e_i)}$. Under $u^{(c)}\mapsto t_1^{-c_1-1}\cdots t_n^{-c_n-1}$, which preserves degrees, $u^{(c)}\mapsto u^{(c-e_i)}$ becomes multiplication by $t_i$ (with $t_i\cdot t_i^{-1}=0$), so this is the complex $K_X$. Together with (a1) and (a2), $H_*(K_X)\cong H_*(ET\times_TX;\mathbb Q)$.

(a4) \emph{Naturality.} All constructions in (a1)--(a3) are natural with respect to pairs $(\alpha,f)$ of an automorphism $\alpha$ of $T$ and a map $f:X\to X$ with $f(\theta\cdot x)=\alpha(\theta)\cdot f(x)$: $\alpha_*$ is an automorphism of the differential graded algebra $C_*(T)$, $f_*$ is $\alpha_*$-semilinear, and $(\alpha,f)$ acts on $|B_\bullet(X,T,\mathrm{pt})|$. If $\alpha$ permutes the circles, then $\alpha_*\circ\iota=\iota\circ\alpha_\Lambda$, $\alpha_\Lambda$ permuting the $D_i$ in the same way, and the induced map of $K_X$ is $f_*$ together with the same permutation of the $t_i$. A relabelling $\gamma\in S_n$ of $X=\Nt_g(n)$ satisfies $\gamma(\theta\cdot y)=\gamma(\theta)\cdot\gamma(y)$, $\gamma$ permuting the circles of $T$, so the isomorphism of (a) is $S_n$-equivariant. For the step $M\otimes_\Lambda K\simeq B(M,\Lambda,\mathbb Q)$ of (a3) we note that $K$ is semifree over $\Lambda$, so a $\Lambda$-linear quasi-isomorphism $K\to B(\Lambda,\Lambda,\mathbb Q)$ over the identity of $\mathbb Q$ exists and is unique up to $\Lambda$-linear homotopy; its conjugate by the actions of a permutation $\alpha$ on $K$ and on $B(\Lambda,\Lambda,\mathbb Q)$ is again such a map, hence homotopic to it, and the induced isomorphism $H_*(M\otimes_\Lambda K)\cong H_*(B(M,\Lambda,\mathbb Q))$ is $S_n$-equivariant.

Smooth chains may be used: the inclusion of smooth into continuous normalized chains is a quasi-isomorphism commuting with the $D_i$, $[S^1]$ being smooth, and the filtration by the total power of the $t_i^{-1}$, which is bounded below and exhaustive, shows that it induces a quasi-isomorphism of the complexes $K_X$.

(b) The projection $ET\times_TX\to X/T$ induces an isomorphism on rational homology; for actions of compact Lie groups with finite stabilisers this is standard, see \cite{Hs,AP}, and we recall the argument. Its fibre over the orbit of $x$ is $ET/T_x$ with $T_x$ finite, which is rationally acyclic. Concretely, for the approximations $E_N=(S^{2N+1})^n$ of $ET$ the map $E_N\times_TX\to X/T$ is proper, with fibres $E_N/T_x$, whose rational cohomology is that of $(S^{2N+1})^n$ since $T_x$ is finite and acts through the connected group $T$; so the Vietoris--Begle theorem, in the form of the Leray spectral sequence of a proper map \cite[Chapters II and IV]{Br97}, gives isomorphisms on rational cohomology, hence on rational homology, in degrees $\le2N$, and $N$ is arbitrary. The space $X/T$ is locally contractible, so sheaf and singular cohomology agree.

(c) Over $\mathbb Q$, coinvariants of the finite group $S_n$ commute with homology, and $H_*(X/T;\mathbb Q)_{S_n}=H_*(X/G_n;\mathbb Q)$, with $X/G_n=|\N_L(S_{g,n})|$ by Lemma \ref{lem:Ntilde}(2). The equivariance of (a) and (b) for $S_n$, which permutes the circles, follows from naturality.

The statement on $H_0$ is elementary: in degree $0$ only $C_0(X)\otimes t_1^{-1}\cdots t_n^{-1}$ occurs, in degree $1$ only $C_1(X)\otimes t_1^{-1}\cdots t_n^{-1}$, and $t_iD_i$ vanishes on the latter, so $H_0(F^{\Nt}_{g,n})=H_0(X;\mathbb Q)_{S_n}=\mathbb Q$.
\end{proof}

\subsection{The residue map}\label{ssec:Phi}
\begin{definition}\label{def:Phi}
For a smooth simplex $\sigma:\Dt^k\to\Nt_g(n)$, $n\ge1$, and exponents $a\in\mathbb Z_{\ge1}^n$ put
\begin{equation*}
\Phi\big(\sigma\otimes t_1^{-a_1}\cdots t_n^{-a_n}\big):=\begin{cases}[\Dt^k,q\circ\sigma]&\text{if }a=(1,\dots,1),\\0&\text{otherwise,}\end{cases}
\end{equation*}
in $C^w_k(\N_L(S_{g,n}))/R^w$, and $\Phi(\hat\sigma):=[\Dt^k,\hat\sigma]$ for $\hat\sigma\in C^{\mathrm{sm}}_k(\N_L(S_{g,0}))$. On the summands of \eqref{eq:Cprime} with several factors, $\Phi$ maps $(\sigma_1\otimes\dots\otimes\sigma_r)\otimes t^{-a}$ to $[\prod_l\Dt^{k_l},\sqcup\circ\prod_l\bar\sigma_l]$ if $a=(1,\dots,1)$ and to $0$ otherwise, where $\bar\sigma_l=q\circ\sigma_l$ for factors with boundary and $\bar\sigma_l=\hat\sigma_l$ for closed factors.
\end{definition}
In the terminology of \cite[Section 2.5]{CCT}, $\Phi$ takes the residue, here the coefficient of $t_1^{-1}\cdots t_n^{-1}$ (in the normalization of \cite{CCT} the coefficient of $u^0$), and projects to the quotient by $G_n$; the projection exists for geometric chains, the quotient being taken as an orbifold.

\begin{lemma}\label{lem:Phi-chain}
$\Phi$ is a well defined, multiplicative chain map $(\mathcal F(F^{\Nt}),d+\sum_it_iD_i)\to(\mathcal C^0/R^w,\partial)$, mapping $F^{\Nt}_{g,n}$ to $C^w_*(\N_L(S_{g,n}))/R^w$, chains of degree $k$ to chains of dimension $k$, and the unit to the unit $[\mathrm{pt}]$.
\end{lemma}
\begin{proof}
\emph{Degenerate simplices.} For $\sigma\circ s_j$ with $s_j:\Dt^{k+1}\to\Dt^k$ a degeneracy, the affine involution of $\Dt^{k+1}$ exchanging the vertices $j$ and $j+1$ reverses orientations and satisfies $s_j\circ a=s_j$; so $[\Dt^{k+1},q\circ\sigma\circ s_j]=0$ by the vanishing statement after Definition \ref{def:wide-chains}, and likewise on closed components and on products. \emph{Coinvariants.} The monomial $t_1^{-1}\cdots t_n^{-1}$ is $S_n$-invariant and $q\circ\gamma\cong q$ for $\gamma\in S_n$ (Lemma \ref{lem:Ntilde}(3)). Exchanging two factors of degrees $k$ and $l$ in \eqref{eq:Cprime} has the sign $(-1)^{kl}$, and so has the exchange $\Dt^k\times\Dt^l\to\Dt^l\times\Dt^k$, while $\sqcup\circ\mathrm{swap}\cong\sqcup$. \emph{Products.} $\Phi(xy)=\Phi(x)\cdot\Phi(y)$ by definition, $x\cdot y=\sqcup_*(x\times y)$ in $\mathcal C^0$, $\sqcup$ being associative up to $2$-isomorphism. \emph{Chain map.} The terms $d\sigma\otimes t^{-a}$ are mapped compatibly, since the comparison map $\hat\sigma\mapsto[\Dt^k,\hat\sigma]$ of Theorem \ref{thm:Tw} is a chain map, and by the Leibniz rule for products. It remains to see that $\Phi(t_iD_i\sigma\otimes t^{-a})=0$. This term has exponent $a_i-1$ at $t_i$ and vanishes unless $a_i\ge2$; so it can contribute only if $a_i=2$ and $a_l=1$ for $l\ne i$, and then it equals
\begin{equation*}
\Phi\big(D_i\sigma\otimes t_1^{-1}\cdots t_n^{-1}\big)=\sum_\mu\mathrm{sign}(\mu)\,\big[\Dt^{k+1},q\circ\rot_i\circ([S^1]\times\sigma)\circ\theta_\mu\big]
\end{equation*}
over the shuffle simplices $\theta_\mu:\Dt^{k+1}\to\Dt^1\times\Dt^k$. By Lemma \ref{lem:Ntilde}(3), $q\circ\rot_i\cong q\circ\mathrm{pr}_2$, and $\mathrm{pr}_2\circ\theta_\mu$ is a degeneracy $\Dt^{k+1}\to\Dt^k$. So each term vanishes, as above; factors not containing $i$ are carried along.
\end{proof}

\begin{proposition}\label{prop:Phi-BV}
On $\mathcal F(F^{\Nt})$,
\begin{equation*}
\Phi\circ\Delta^{\Cos}\equiv-\,\Delta\circ\Phi\qquad\text{modulo }R^w.
\end{equation*}
Let $\epsilon$ be the linear automorphism of the graded vector space $\mathcal C^0/R^w$ which is $(-1)^{\dim\N_{\chi,n}/2}$ on $C^w_*(\N_{\chi,n})/R^w$; on connected types it is $\epsilon_{g,n}:=(-1)^{3g-3+n}$. Then $\epsilon\circ\Phi:\mathcal F(F^{\Nt})\to\mathcal C^0/R^w$ is a morphism of differential graded BV algebras, and with \eqref{eq:lambda-hbar} it restricts to a filtration preserving morphism of differential graded Lie algebras
\begin{equation*}
\epsilon\Phi:\mathfrak g(F^{\Nt})\to\mathfrak h^w,
\end{equation*}
which maps $y\otimes t_1^{-1}t_2^{-1}t_3^{-1}$, for a point $y$ of $\Nt_0(3)$, to the point $[\mathrm{pt},q\circ y]$, and hence preserves the degree of the component $(0,3)$, i.e.\ the sum of the coefficients of the $0$-chain multiplying $t_1^{-1}t_2^{-1}t_3^{-1}$, resp.\ the weighted count of points.
\end{proposition}
\begin{proof}
By Lemma \ref{lem:Phi-chain} it suffices to compare, for $x=\sigma\otimes t^{-a}$ and a pair $\{i,j\}$, $\Phi(\Delta_{ij}x)$ with the part of $\Delta\Phi(x)$ which sews $\{i,j\}$, for $x$ supported on one factor of \eqref{eq:Cprime} containing $i$ and $j$, or on two factors containing $i$ and $j$ respectively. Further factors, closed ones included, are carried along on both sides: on Costello's side by definition, on ours by Lemma \ref{lem:glue-sqcup} and the base change in the proof of Proposition \ref{prop:BV-chains}; the Koszul sign for moving the circle past a preceding factor is the same on both sides. If $a\ne(1,\dots,1)$ both sides vanish, since $\Delta_{ij}$ keeps the exponents of the variables other than $t_i,t_j$ and kills the terms with $a_i\ge2$ or $a_j\ge2$. So let $a=(1,\dots,1)$, let $Q$ be $\Dt^k$, resp.\ $\Dt^{k_1}\times\Dt^{k_2}$, and $\sigma:Q\to\tilde Y$ with $\tilde Y$ as in Lemma \ref{lem:pullback-E}, so that $\Phi(x)=[Q,q_{\tilde Y}\circ\sigma]$.

\emph{Our side.} By Lemma \ref{lem:pullback-E}(1) and the pasting of $2$-pullbacks, the part of $\Delta[Q,q_{\tilde Y}\circ\sigma]=\glue_*\pi^*p^*[Q,q_{\tilde Y}\circ\sigma]$ sewing $\{i,j\}$ is
\begin{equation*}
\big[(\mathbb R/L\mathbb Z)\times Q,\ c_{ij}\circ(\mathrm{id}\times\sigma)\big],
\end{equation*}
the circle oriented by $t$ and placed first, exactly in $C^w_*$.

\emph{Costello's side.} By the associativity of the shuffle map, $D_i$ followed by $G_{ij}$, after the shuffle map in the case of two factors, is $(G_{ij}\circ\rot_i)_*\mathrm{EZ}([S^1]\otimes\sigma)$, with $\sigma=\mathrm{EZ}(\sigma_1\otimes\sigma_2)$ in the second case. Hence
\begin{gather*}
\Phi(\Delta_{ij}x)=\sum_\mu\mathrm{sign}(\mu)\,\big[\Dt^{\dim Q+1},\ b_{ij}\circ([S^1]\times\mathrm{id}_Q)\circ\theta_\mu\big],\\
b_{ij}:=\bar G_{ij}\circ\rot_i\circ(\mathrm{id}\times\sigma),
\end{gather*}
the sum over the simplices $\theta_\mu$ of the shuffle triangulation of $\Dt^1\times Q$. The shuffle map is compatible with the product orientation: $\mathrm{sign}(\mu)$ is the orientation sign of $\theta_\mu$. So the $\theta_\mu$, oriented from $\Dt^1\times Q$, form a gluing datum in the sense of Definition \ref{def:Rw}: they are embeddings, cover $\Dt^1\times Q$ with disjoint interiors, and each interior face is shared by two distinct simplices, the other faces lying in the boundary. Therefore $\Phi(\Delta_{ij}x)\equiv[\Dt^1\times Q,b_{ij}\circ([S^1]\times\mathrm{id})]$. Cutting $\Dt^1=[0,1]$ at $\frac12$ is again a gluing datum, so
\begin{equation*}
\Phi(\Delta_{ij}x)\equiv\sum_{l=1,2}\big[A_l\times Q,\ b_{ij}\circ([S^1]\times\mathrm{id})\big],\qquad A_1=[0,\tfrac12],\ A_2=[\tfrac12,1].
\end{equation*}
On the simply connected $A_l\times Q$, Lemma \ref{lem:pullback-E}(2) gives $b_{ij}\circ([S^1]\times\mathrm{id})\cong c_{ij}\circ(\nu\times\mathrm{id})\circ([S^1]\times\sigma)$, and $(\nu\circ[S^1])\times\mathrm{id}_Q$ is an orientation reversing diffeomorphism of $A_l\times Q$ onto $B_l\times Q$, with the arcs $B_1=[-\frac L2,0]$ and $B_2=[-L,-\frac L2]$ of $\mathbb R/L\mathbb Z$. Hence
\begin{equation*}
\Phi(\Delta_{ij}x)\equiv-\sum_{l=1,2}\big[B_l\times Q,\ c_{ij}\circ(\mathrm{id}\times\sigma)\big]\equiv-\big[(\mathbb R/L\mathbb Z)\times Q,\ c_{ij}\circ(\mathrm{id}\times\sigma)\big],
\end{equation*}
the last step by the gluing datum decomposing the circle into the two arcs: each end of $B_1\times Q$ is glued to an end of $B_2\times Q$, as (G3) requires $j\ne i$, so that a single interval with its ends identified would not do. The comparison of the two morphisms is made only on the simply connected pieces; on $(\mathbb R/L\mathbb Z)\times Q$ the families of markings of Lemma \ref{lem:sc} have monodromy, a power of the Dehn twist about the sewn curve. Summing over $i<j$ gives $\Phi\circ\Delta^{\Cos}\equiv-\Delta\circ\Phi$.

The map $\epsilon$ is multiplicative, since dimensions add under $\sqcup$, commutes with $\partial$, and anticommutes with $\Delta$, which maps $C^w_*(\N_{\chi,n+2})$ to $C^w_*(\N_{\chi,n})$ while $\dim\N_{\chi,n+2}=\dim\N_{\chi,n}+2$. Hence $\epsilon\Phi$ is multiplicative, a chain map, and $\epsilon\Phi\circ\Delta^{\Cos}=\Delta\circ\epsilon\Phi$; so it is a morphism of BV algebras, and it preserves the brackets, which are determined by the products and the BV operators. On connected elements it maps $F^{\Nt}_{g,n}$ to $C^w_*(\N_L(S_{g,n}))/R^w$, elements of degree $6g-6+2n+1+i$ to elements of dimension $6g-6+2n+1+i$, and it preserves the filtrations. The last statement holds by definition, since $\epsilon_{0,3}=1$ and $[\mathrm{pt},q\circ y]$ has weighted count $1$.
\end{proof}

\begin{remark}\label{rem:signs}
(1) The linear map $\epsilon$ is not an automorphism of the BV algebra $\mathcal C^0/R^w$: it is multiplicative and commutes with $\partial$, but it is an isomorphism onto the BV algebra with $\Delta$ replaced by $-\Delta$, and on the connected part an isomorphism of differential graded Lie algebras $(\partial,\Delta_0,\{\,,\})\cong(\partial,-\Delta_0,-\{\,,\})$; indeed $\epsilon_{g-1,n+2}=-\epsilon_{g,n}$, $\epsilon_{t_1}\epsilon_{t_2}=-\epsilon_{g,n}$ for the types $t_1,t_2$ of a separating sewing, and $\epsilon_{0,3}=1$, so that the initial condition is preserved.

(2) The sign in Proposition \ref{prop:Phi-BV} comes from the circles. With the conventions of Section \ref{ssec:Ntilde}, Costello's $D_i$ moves the marked points against the boundary orientation, so that the gluings $G_{ij}\circ\rot_i(\theta)$ are right twists in the sense of Definition \ref{def:left-twist}, while $\Delta$ uses the rotation orientation \eqref{eq:rotation}, i.e. left twists (Lemma \ref{lem:cut}(2)). This is consistent with Remark \ref{rem:normalization}(3), where the same factor $(-1)^{3g-3+n}$ relates our orientation to the complex one, which goes with right twists. Costello's \cite[Theorem 9.0.7]{Cos05} states that his solution represents the orbifold fundamental class; that the orientation meant is the complex one is our reading, supported by the use of the canonical orientation of a complex orbifold in his proof. We have not compared the orientation conventions of \cite{Cos05} beyond this. If $D$ is oriented oppositely, $\Delta^{\Cos}$ and Costello's bracket change sign, $\Phi$ itself is the morphism, and the statements below hold with $\Phi$ in place of $\epsilon\Phi$.
\end{remark}

\subsection{The current level}\label{ssec:E}
Let $\Psi:=J\circ\epsilon\circ\Phi$ on $\mathfrak g(F^{\Nt})$, with values in the differential graded Lie algebra $\mathfrak g[[\hbar]]$ of compactly supported currents of Proposition \ref{prop:BV}, powers of $\lambda$ being converted by \eqref{eq:lambda-hbar}. Explicitly, for a form $\omega$ on $\N_L(S_{g,n})$,
\begin{equation}\label{eq:Psi}
\begin{gathered}
\Psi\big(\sigma\otimes t_1^{-a_1}\cdots t_n^{-a_n}\big)(\omega)=\delta_{a,(1,\dots,1)}\,\epsilon_{g,n}\int_{\Dt^k}\sigma^*q^*\omega,\\
\Psi(\hat\sigma)(\omega)=\epsilon_{g,0}\int_{\Dt^k}\hat\sigma^*\omega .
\end{gathered}
\end{equation}
By Corollary \ref{cor:exp} and Proposition \ref{prop:Phi-BV}, $\Psi$ is a morphism of differential graded Lie algebras. We give a direct proof, which uses neither the relations $R^w$ nor Theorem \ref{thm:Tw}; it rests on an identity of forms, in which, unlike for chains, no global comparison of morphisms is needed, pullback of forms being local.

\begin{lemma}\label{lem:formal}
In the situation of Lemma \ref{lem:pullback-E}, for every form $\omega$ on $\N'$,
\begin{equation*}
q_{\tilde Y}^*\big(p_!\,\pi_!\,\glue^*\omega\big)\big|_{\Pi}=-\sum_{\{i,j\}\in\Pi}(\mathrm{pr}_2)_!\big(\bar G_{ij}\circ\rot_i\big)^*\omega
\end{equation*}
as forms on $\tilde Y$, where on the left $p_!\pi_!$ is restricted to the part $\Pp^\Pi$, $i$ denotes the smaller label of the pair in the first case and the label in $I_1$ in the second, and $(\mathrm{pr}_2)_!$ is integration over the fibre of $\mathrm{pr}_2:(\mathbb R/\mathbb Z)\times\tilde Y\to\tilde Y$, fibre first, $\mathbb R/\mathbb Z$ oriented by $\theta$.
\end{lemma}
\begin{proof}
Integration along the fibres commutes with pullback along base changes, fibrewise and with the fibre orientations, since both are computed in charts. By Lemma \ref{lem:pullback-E}(1), the base change of $\Pp^\Pi$ and $\E^\Pi$ along $q_{\tilde Y}$ is $\bigsqcup_{\{i,j\}\in\Pi}(\mathbb R/L\mathbb Z)\times\tilde Y$, with $\glue$ pulled back to $c_{ij}$; so the left hand side is $\sum_{\{i,j\}}(\mathrm{pr}_2)_!c_{ij}^*\omega$, integrated over $\mathbb R/L\mathbb Z$ oriented by $t$. Locally on $(\mathbb R/\mathbb Z)\times\tilde Y$, $c_{ij}\circ(\nu\times\mathrm{id})\cong\bar G_{ij}\circ\rot_i$ by Lemma \ref{lem:pullback-E}(2), so $(\nu\times\mathrm{id})^*c_{ij}^*\omega=(\bar G_{ij}\circ\rot_i)^*\omega$, and $\nu$ reverses the orientation of the fibre.
\end{proof}

\begin{proposition}\label{prop:Psi}
$\Psi$ is a filtration preserving morphism of differential graded Lie algebras from $\mathfrak g(F^{\Nt})$ to the algebra $(\mathfrak g[[\hbar]],\partial+\hbar\Delta_0,[\,,])$ of Proposition \ref{prop:BV}(2), restricted to the series with $\hbar^g$ on $\N_L(S_{g,n})$.
\end{proposition}
\begin{proof}
\emph{Differential.} For $a=(1,\dots,1)$, $\Psi(d\sigma\otimes t^{-a})=\partial\Psi(\sigma\otimes t^{-a})$ by Stokes' theorem on $\Dt^k$. The terms $t_iD_i\sigma\otimes t^{-a}$ with $a_i=2$, $a_l=1$ otherwise, give $\int_{\Dt^{k+1}}\theta_\mu^*([S^1]\times\sigma)^*\rot_i^*q^*\omega$, and $\rot_i^*q^*\omega=\mathrm{pr}_2^*q^*\omega$ by Lemma \ref{lem:Ntilde}(3), so the integrand is pulled back along the degeneracy $\mathrm{pr}_2\circ\theta_\mu$ from a $(k+1)$-form on $\Dt^k$, and vanishes. \emph{Sewing.} For $x=\sigma\otimes t_1^{-1}\cdots t_n^{-1}$ on $\Nt_g(n)$ and a form $\omega$ on the target of $\Delta_{ij}$, the shuffle map being compatible with the product orientation and $[S^1]$ a diffeomorphism from the interior of $\Dt^1$ onto the circle minus a point,
\begin{align*}
\Psi(\Delta_{ij}x)(\omega)&=\epsilon_{g+1,n-2}\int_{\Dt^1\times\Dt^k}([S^1]\times\sigma)^*(\bar G_{ij}\circ\rot_i)^*\omega\\
&=\epsilon_{g+1,n-2}\int_{\Dt^k}\sigma^*(\mathrm{pr}_2)_!(\bar G_{ij}\circ\rot_i)^*\omega
\end{align*}
by Fubini's theorem, fibre first. Summing over $i<j$ and using Lemma \ref{lem:formal} and $\epsilon_{g+1,n-2}=-\epsilon_{g,n}$, this is $\epsilon_{g,n}\int_{\Dt^k}\sigma^*q^*(p_!\pi_!\glue^*\omega)=(\Delta_0\Psi(x))(\omega)$, since $(\Delta T)(\omega)=T(p_!\pi_!\glue^*\omega)$ for currents (Section \ref{ssec:currents}). The same computation with $\tilde Y=\Nt_{g_1}(n_1)\times\Nt_{g_2}(n_2)$, the second case of Lemma \ref{lem:formal}, $\epsilon_{t_1}\epsilon_{t_2}=-\epsilon_{g,n}$ and $\Psi(\mathrm{EZ}(\sigma_1\otimes\sigma_2)\otimes\cdots)=\Psi(\sigma_1\otimes\cdots)\times\Psi(\sigma_2\otimes\cdots)$ (Fubini) gives the compatibility with the brackets, both brackets being the sums of the sewings of one boundary component of each factor, $D_i$ acting on the first factor. Terms with some exponent $\ge2$ are mapped to $0$ on both sides. The filtrations are preserved since $\Psi$ preserves the types.
\end{proof}

\begin{theorem}\label{thm:E}
Let $0<L\le\Lst$, let $D$ be oriented as in Section \ref{ssec:Ntilde}, and let homotopies be as in Section \ref{ssec:Cos-framework}.
\begin{enumerate}
\item Let $S=\sum S_{g,n}\lambda^{2g-2+n}$ be a Maurer--Cartan element of $\mathfrak g(F^{\Nt})$ whose component $S_{0,3}$ has degree $1/3!$. For $n\ge1$ let $S^\circ_{g,n}\in C_{6g-6+2n}(\Nt_g(n))_{S_n}$ be the coefficient of $t_1^{-1}\cdots t_n^{-1}$ in $S_{g,n}$, and put $S^\circ_{g,0}:=S_{g,0}$. Then $\Psi(S)=T$, i.e. for all $(g,n)$ with $2g-2+n>0$ and all $(6g-6+2n)$-forms $\omega$ on $\N_L(S_{g,n})$,
\begin{equation*}
(-1)^{3g-3+n}\int_{S^\circ_{g,n}}q^*\omega=\int_{V_{g,n}}\omega,
\end{equation*}
the right hand side with the weights of Definition \ref{def:weighted} (for $n=0$ the left hand side is $\int_{S_{g,0}}\omega$ with the sign). The sign $(-1)^{3g-3+n}$ depends on the orientation of $D$; for the opposite orientation it disappears (Remark \ref{rem:signs}(2)).
\item Maurer--Cartan elements $S$ as in (1) exist and are unique up to homotopy. Assume that $h$ is continuous (Section \ref{ssec:Ntilde}). Then $h_*$ maps their homotopy class to the class of Costello's solution in $\mathfrak g(C'_*(M^s))$; this relies on Costello's computation of the homology of $\mathfrak g(C'_*(M^s))$ in \cite[Proposition 9.0.2]{Cos05}, which rests on \cite[Proposition 5.1.1]{Cos05}, whose proof is sketched, and on the rational homotopy equivalence of $M_g(n)/(S^1)^n$ with $M_{g,n}$ (Section \ref{ssec:Cos-framework}).
\end{enumerate}
\end{theorem}
\begin{proof}
(1) By Proposition \ref{prop:Psi}, $\Psi(S)=\sum_{g,n}\hbar^g\Psi(S_{g,n})$ is a Maurer--Cartan element. Its component on $\N_L(S_{g,n})$ is the current $\epsilon_{g,n}\int_{S^\circ_{g,n}}q^*(\cdot)$, which has compact support and top dimension $6g-6+2n$, since only the part of $S_{g,n}$ of chain degree $6g-6+2n$ multiplying $t_1^{-1}\cdots t_n^{-1}$ contributes. Its component on $\N_{0,3}$, a point with group $S_3$, is $f\mapsto\frac16f(\mathrm{pt})$, which is $T_{0,3}$ (Remark \ref{rem:normalization}(2)). By Theorem \ref{thm:D}, $\Psi(S)=T$.
(2) By Lemma \ref{lem:low-degrees} below, $\mathfrak g(F^{\Nt})$ satisfies both hypotheses of Lemma \ref{lem:MC}, and the degree of the component $(0,3)$, i.e.\ the sum of the coefficients of the $0$-chain multiplying $t_1^{-1}t_2^{-1}t_3^{-1}$, is a well defined linear form on the cycles of degree $-1$ of $\mathrm{Gr}^1$ which induces an isomorphism $H_{-1}(\mathrm{Gr}^1)\cong\mathbb Q$; a generator is $y\otimes t_1^{-1}t_2^{-1}t_3^{-1}+b$ with $b$ in the component $(1,1)$. Since $\Delta^{\Cos}$ maps the component $(0,3)$ to the component $(1,1)$, $c(S)$ is the class of the cycle $S_{0,3}\lambda+S_{1,1}\lambda$, not of $S_{0,3}\lambda$ alone; by the linear form, $c(S)$ is the degree of $S_{0,3}$ times this generator, as in Section \ref{ssec:Cos-framework}. So Lemma \ref{lem:MC} gives existence and uniqueness up to homotopy. If $h$ is continuous, $h_*$ is a filtration preserving morphism and maps such an $S$ to a Maurer--Cartan element of $\mathfrak g(C'_*(M^s))$ whose component $(0,3)$ has degree $1/3!$. By Lemma \ref{lem:low-degrees} for $C'_*(M^s)$, which uses the input of \cite[Proposition 9.0.2]{Cos05}, Lemma \ref{lem:MC}(3) applies to $\mathfrak g(C'_*(M^s))$, and the classes $c$ of $h_*(S)$ and of Costello's solution are both $1/3!$ times the generator, so $h_*(S)$ is homotopic to Costello's solution.
\end{proof}
Part (1) is an exact statement for every normalized solution, not only up to homotopy; homotopic solutions have the same current, as currents of top dimension leave no room for homotopies. It depends only on Proposition \ref{prop:Psi} and Theorem \ref{thm:D}, not on Proposition \ref{prop:equivariant}, Theorem \ref{thm:Tw}, or the homotopy type of the moduli spaces. Part (1) concerns Maurer--Cartan elements in the model $\Nt$. For Costello's solutions in $C'_*(M^s)$ themselves there is no direct integration statement, since there is no map from $M^s$ to the moduli spaces $\N_L$ compatible with gluing along which their chains could be pushed forward; they are related to the vertices only through homotopy classes, by Part (2) and Theorem \ref{thm:F}. The first sentence of Part (2), existence and uniqueness up to homotopy in $\mathfrak g(F^{\Nt})$, depends on Proposition \ref{prop:equivariant}, whose case $n=0$ uses the proof of Theorem \ref{thm:Tw}, on Lemma \ref{lem:bound}, which uses Harer's bound \cite{Ha}, the identification of the rational homology of $|\N_L(S_{g,n})|$ with that of $\Mod(S_{g,n})$ \cite{Bro} and, for the type $(0,4)$, the vanishing of the $S_4$-coinvariants of $H_1(M_{0,4};\mathbb Q)$ recalled from \cite{CZ,Cos05}, and on Lemma \ref{lem:MC}; it does not use the sketched Lemma \ref{lem:GM}. The second sentence, that these solutions represent Costello's class, is conditional: it assumes that $h$ is continuous, Costello not having specified the topology of his moduli spaces, and it uses his computation of the homology of $\mathfrak g(C'_*(M^s))$, which rests on \cite[Proposition 5.1.1]{Cos05}, whose proof is sketched, and on the rational homotopy equivalence $M_g(n)/(S^1)^n\simeq_{\mathbb Q}M_{g,n}$ stated in the proof of \cite[Proposition 9.0.2]{Cos05}; Costello's solution itself exists under the same input.

\subsection{Homotopy classes of Maurer--Cartan elements}\label{ssec:F}
\begin{lemma}\label{lem:low-degrees}
Let $\mathfrak L$ be $\mathfrak g(F^{\Nt})$ or $\mathfrak h^w$. Then $H_i(\mathrm{Gr}^k\mathfrak L)=0$ for $i\in\{0,-1,-2\}$ and all $k\ge1$, except $H_{-1}(\mathrm{Gr}^1\mathfrak L)\cong\mathbb Q$. The degree of the component $(0,3)$, i.e.\ the sum of the coefficients of the $0$-chain multiplying $t_1^{-1}t_2^{-1}t_3^{-1}$, resp.\ the weighted count of points, is a well defined linear form on the cycles of degree $-1$ of $\mathrm{Gr}^1\mathfrak L$, vanishing on boundaries, and it induces this isomorphism. A generator is $y\otimes t_1^{-1}t_2^{-1}t_3^{-1}+b$, resp.\ $[\mathrm{pt},\hat y]+\hat b$, for any point $y$ of $\Nt_0(3)$, resp.\ $\hat y$ of $\N_L(S_{0,3})$, with $b$, resp.\ $\hat b$, in the component $(1,1)$. The same holds for $\mathfrak L=\mathfrak g(C'_*(M^s))$, with a generator $y\otimes t_1^{-1}t_2^{-1}t_3^{-1}+b$ for any point $y$ of $M_0(3)$, by \cite[proof of Proposition 9.0.2]{Cos05}, granting the input recalled in Section \ref{ssec:Cos-framework}.
\end{lemma}
\begin{proof}
$\mathrm{Gr}^k\mathfrak L=\prod_{2g-2+n=k}\mathfrak L_{g,n}$, with the differential $d+\Delta$, resp.\ $\partial+\Delta_0$, where $\Delta$ raises the genus. Filtering by the genus, a finite filtration, the first page is $\bigoplus_{2g-2+n=k}H_*(\mathfrak L_{g,n})$, and $H_*(\mathfrak L_{g,n})\cong H_*(|\N_L(S_{g,n})|;\mathbb Q)$ by Proposition \ref{prop:equivariant}, resp.\ Theorem \ref{thm:Tw}. The degrees $i=0,-1,-2$ correspond to the dimensions $6g-5+2n$, $6g-6+2n$ and $6g-7+2n$, where these groups vanish by Lemma \ref{lem:bound}, except for $(g,n)=(0,3)$, where $H_0=\mathbb Q$ is spanned by a point (Proposition \ref{prop:equivariant}, resp.\ the comparison map of Theorem \ref{thm:Tw}) and $H_1=H_{-1}=0$. The only type with $k=1$ other than $(0,3)$ is $(1,1)$, whose groups in the relevant dimensions $1,2,3$ vanish. So the spectral sequence has no room for differentials in the relevant degrees, and the vanishing statements follow. For $k=1$, since $\Delta$ raises the genus, the component $(1,1)$ is a subcomplex of $\mathrm{Gr}^1\mathfrak L$ with quotient the component $(0,3)$, and $\Delta$ maps the component $(0,3)$ to the component $(1,1)$; so a point of the component $(0,3)$ alone is not a cycle of $\mathrm{Gr}^1\mathfrak L$. In the long exact sequence of this extension the groups of the component $(1,1)$ in degrees $-1$ and $-2$, i.e.\ in dimensions $2$ and $1$, vanish, so the projection to the component $(0,3)$ induces $H_{-1}(\mathrm{Gr}^1\mathfrak L)\cong H_0$ of the component $(0,3)$, which is $\mathbb Q$ via the degree (Proposition \ref{prop:equivariant}, resp.\ the weighted count of points of the proof of Theorem \ref{thm:Tw}), and every point $y$ lifts to a cycle $y\otimes t_1^{-1}t_2^{-1}t_3^{-1}+b$, resp.\ $[\mathrm{pt},\hat y]+\hat b$: the obstruction, the class of $\Delta$ of the point in the component $(1,1)$, lies in a vanishing group. The degree vanishes on boundaries, since boundaries of degree $-1$ have as component $(0,3)$ the boundary of a $1$-chain multiplying $t_1^{-1}t_2^{-1}t_3^{-1}$, the operators $t_iD_i$ vanishing on these monomials.
\end{proof}

\begin{theorem}\label{thm:F}
Let $0<L\le\Lst$, and let homotopies be as in Section \ref{ssec:Cos-framework}.
\begin{enumerate}
\item The morphism $\epsilon\Phi$ induces a bijection
\begin{equation*}
\pi_0\mathrm{MC}\big(\mathfrak g(F^{\Nt})\big)\ \xrightarrow{\ \epsilon\Phi\ }\ \pi_0\mathrm{MC}(\mathfrak h^w).
\end{equation*}
Both sets are identified with $\mathbb Q$ by the degree of the component $(0,3)$, i.e.\ by the class in $H_{-1}(\mathrm{Gr}^1)$, and the map respects these identifications. In particular $\epsilon\Phi(S)$ is homotopic to the hyperbolic vertices $V_\hbar=\sum_{g,n}\hbar^g[V_{g,n}]$ in $\mathfrak h^w$ for every Maurer--Cartan element $S$ of $\mathfrak g(F^{\Nt})$ whose component $S_{0,3}$ has degree $1/3!$.
\item Assume that $h$ is continuous (Section \ref{ssec:Ntilde}). Then $h_*$ induces a bijection
\begin{equation*}
\pi_0\mathrm{MC}\big(\mathfrak g(F^{\Nt})\big)\ \xrightarrow{\ h_*\ }\ \pi_0\mathrm{MC}\big(\mathfrak g(C'_*(M^s))\big),
\end{equation*}
compatible with the identifications with $\mathbb Q$, and under the bijections of (1) and (2) the class of Costello's solution corresponds to the class of $V_\hbar$. This part relies on Costello's computation of the homology of $\mathfrak g(C'_*(M^s))$ in \cite[Proposition 9.0.2]{Cos05}, which rests on \cite[Proposition 5.1.1]{Cos05}, whose proof is sketched, and on the rational homotopy equivalence of $M_g(n)/(S^1)^n$ with $M_{g,n}$ (Section \ref{ssec:Cos-framework}).
\end{enumerate}
\end{theorem}
\begin{proof}
The algebras are filtered pro-nilpotent, and $\epsilon\Phi$ and, if $h$ is continuous, $h_*$ are filtration preserving morphisms (Proposition \ref{prop:Phi-BV} and Section \ref{ssec:Ntilde}). By Lemma \ref{lem:low-degrees}, which for $\mathfrak g(C'_*(M^s))$ uses the input of \cite[Proposition 9.0.2]{Cos05}, each algebra satisfies both hypotheses of Lemma \ref{lem:MC}, and $H_{-1}(\mathrm{Gr}^1)\cong\mathbb Q$ via the degree of the component $(0,3)$, with generators $y\otimes t_1^{-1}t_2^{-1}t_3^{-1}+b$, resp.\ $[\mathrm{pt},\hat y]+\hat b$; both morphisms induce chain maps on $\mathrm{Gr}^1$ which preserve this degree (Proposition \ref{prop:Phi-BV} and Section \ref{ssec:Ntilde}), hence map a generator to a generator. So Lemma \ref{lem:MC} identifies each $\pi_0\mathrm{MC}$ with $\mathbb Q$, compatibly with the maps, and gives the bijections; the identification with the degree of the component $(0,3)$ is as in the proof of Theorem \ref{thm:E}(2). Costello's solution has degree $1/3!$ by definition. In $\mathfrak h^w$, $[V_{0,3}]=[BS_3,\iota]\equiv\frac16[\mathrm{pt},\iota\circ\hat\pi]$ by the covering relation for the covering $\mathrm{pt}\to BS_3$ of degree $6$, so the class of $V_\hbar$ is $\frac16$ as well.
\end{proof}

Theorem \ref{thm:F} does not use that $\Phi$ or $h$ is a quasi-isomorphism: Lemma \ref{lem:MC} only needs the vanishing of $H_{-1}$ and $H_{-2}$ of $\mathrm{Gr}^k$ for $k\ge2$ and the comparison of $H_{-1}(\mathrm{Gr}^1)$, where the degrees of the components $(0,3)$ match. Part (1) is unconditional. Under the hypotheses of Part (2), $V_\hbar$, the Maurer--Cartan element of Corollary \ref{cor:exp}, represents Costello's homotopy class, in the sense that it corresponds to it under the bijections induced by the morphisms $h_*$ and $\epsilon\Phi$ of Costello's algebra for the model $\Nt$.

\begin{remark}\label{rem:full-qi}
On $\mathrm{Gr}$, $\Phi$ factors through the residue map $F^{\Nt}_{g,n}\to C^{\mathrm{sm}}_*(\N_L(S_{g,n}))$, $\sigma\otimes t^{-a}\mapsto\delta_{a,(1,\dots,1)}\,q\circ\sigma$, followed by the comparison map of Theorem \ref{thm:Tw}, which induces an isomorphism on homology. We expect the residue map to be a rational quasi-isomorphism inducing the isomorphism of Proposition \ref{prop:equivariant}; this is the analogue for locally free actions of the assertion of \cite[Section 2.5]{CCT} for free actions, which is stated there without proof. Then $\epsilon\Phi$ would be a filtered quasi-isomorphism, hence a quasi-isomorphism for $\hat d$ by the spectral sequence argument of \cite[Lemma 9.0.5]{Cos05}, and likewise $h_*$ if $h$ is a rational weak equivalence, as one expects, both sides being rationally classifying spaces of mapping class groups for $n\ge1$. This would identify the full BV algebras up to quasi-isomorphism, not only the homotopy classes of Maurer--Cartan elements. We do not prove this and do not use it.
\end{remark}

\begin{remark}\label{rem:n0}
\emph{Closed surfaces.} Costello's spaces include components without boundary, with chains on the coarse moduli spaces $M_g(0)$ \cite[Section 6]{Cos05}. Singular simplices in $|\N_L(S_{g,0})|$ do not lift to orbifold morphisms, so $\Phi$ could not be defined on them; this is why $F^{\Nt}$ uses $C^{\mathrm{sm}}_*$ on closed components, which computes the rational homology of $|\N_L(S_{g,0})|\cong M_g(0)$ (proof of Theorem \ref{thm:Tw}) and on which $h$ is still defined. Costello's constructions in \cite[Sections 7--9]{Cos05} only use a functor $H_*(\mathsf A)\to\mathrm{Comp}_{\mathbb Q}$. Alternatively, the components with $n=0$ of each of the three algebras form a dg ideal $\mathfrak i$ on which the bracket vanishes identically, closed surfaces having no boundary to sew; since $H_j(M_g(0);\mathbb Q)\cong H_j(|\N_L(S_{g,0})|;\mathbb Q)=0$ for $j>4g-5$ (proof of Lemma \ref{lem:bound}), and the degrees $0,-1,-2,-3$ correspond to the dimensions $6g-5,\dots,6g-8>4g-5$, its associated graded has no homology in these degrees, so by the long exact sequence $\mathrm{Gr}(\mathfrak L/\mathfrak i)$ has the same $H_{-1}$ and $H_{-2}$ as $\mathrm{Gr}\,\mathfrak L$, while $\mathrm{Gr}^1\mathfrak i=0$, and $\mathfrak L\to\mathfrak L/\mathfrak i$ induces a bijection on $\pi_0\mathrm{MC}$ by Lemma \ref{lem:MC}. Hence Theorem \ref{thm:F} holds equally for the quotients by $\mathfrak i$, which do not involve closed surfaces at all; and in Theorem \ref{thm:E}(1) the components $T_{g,0}$ are determined by the others, by the argument of Theorem \ref{thm:D}.
\end{remark}

\begin{remark}\label{rem:sketches}
\emph{What is used from \cite{Cos05}.} Our construction does not use \cite[Lemma 9.0.4]{Cos05}, which concerns the comparison with Deligne--Mumford spaces in \cite[Theorem 9.0.3]{Cos05} only. Of the geometric statements of \cite{Cos05} whose proofs are sketched there, the topology and homotopy type of $M_g(n)$ and \cite[Proposition 5.1.1]{Cos05}, Theorem \ref{thm:E}(1) uses none. For $\Nt$, Theorem \ref{thm:E}(2) and Theorem \ref{thm:F} replace them by Lemma \ref{lem:Ntilde}, where gluing is strictly associative and the homotopy type is explicit, and by Proposition \ref{prop:equivariant}, proved with the Koszul resolution and the Vietoris--Begle theorem. They enter only through the transfer to Costello's own model $C'_*(M^s)$, i.e.\ the second sentence of Theorem \ref{thm:E}(2) and Theorem \ref{thm:F}(2), which are stated under the assumption that $h$ is continuous and use Costello's computation of the homology of $\mathfrak g(C'_*(M^s))$. From the algebraic part of \cite{Cos05} we use the construction of $\mathcal F(F)$ and $\mathfrak g(F)$; instead of Lemma \ref{lem:GM}, which Costello sketches with references, we use Lemma \ref{lem:MC}, proved here, with homotopies over the completed algebra $\mathfrak g\hat\otimes\mathbb Q[t,dt]$.
\end{remark}

\subsection{Relation to the literature}\label{ssec:lit}
We are not aware of a previous proof that explicit geometric string vertices, hyperbolic or of minimal area, represent the homotopy class of solutions of \cite{Cos05}; the combinatorial string vertices of \cite{CCT} discussed below represent it by construction, via the homotopy invariance of Maurer--Cartan moduli spaces (they use the framed version of the differential graded Lie algebra of \cite{Cos07}). Costello's definition is homotopical, sewing being defined on homotopy coinvariants because it is not defined on singular chains of the strict quotients \cite[footnote 6, p.~5]{Cos05}. Costello and Zwiebach \cite[Section 2]{CZ} review the homological existence and uniqueness argument of \cite{Cos05} for string vertices which are $S_n$-invariant singular chains with real coefficients on the spaces $\hat{\mathcal P}_{g,n}$ of surfaces with local coordinates defined up to phase, defined in \cite[Section 1, (1.3)]{CZ} and used in \cite[Section 2.1]{CZ}, i.e. on strict quotients by the circle actions, and they note that this material is not required for the construction of the hyperbolic vertices, which are then constructed as sets. The review of Sen and Zwiebach \cite[Section 7.5]{SZ24} attributes the existence of string vertices to Costello, with full details of the argument in \cite{CZ}. C\u{a}ld\u{a}raru, Costello and Tu \cite{CCT} work with homotopy coinvariants throughout and construct combinatorial string vertices from ribbon graphs, as the solution, unique up to homotopy, of the master equation in a resolution of the differential graded Lie algebra of \cite{Cos05}, so that they represent Costello's class by construction, as above; they state without proof that taking the residue and projecting, $C_*(X)_{(S^1)^n}\to C_*(X/(S^1)^n)$, is a quasi-isomorphism for free actions \cite[Section 2.5]{CCT}. Tu \cite{Tu} expresses the fundamental classes of the Deligne--Mumford spaces in terms of these homotopical string vertices. Harrelson, Voronov and Z\'u\~niga \cite{HVZ} set up a quantum master equation for open-closed moduli spaces, solved by their fundamental chains. The work on the hyperbolic vertices themselves, for instance on local coordinates for the three-string vertex \cite{Fi21}, on open-closed vertices in the large $N$ limit \cite{Fi24} and on volumes \cite{FVM}, is at the level of sets and of integrals. The model $\Nt$ and the residue map $\Phi$ connect the strict quotients, on which the hyperbolic vertices are geometric chains, with Costello's homotopy coinvariants.

\section*{Declaration on the use of AI}
The author uses AI to check the validity of the proofs and, under the author's direction, to sketch some new proofs. The overall ideas and methods are from the author.

\end{document}